\documentclass{jgokova}
\usepackage{pstricks}
\usepackage{pst-node}
\usepackage{amssymb}

\newcommand{\ben}{\begin{enumerate}}
\newcommand{\een}{\end{enumerate}}
\newcommand{\be}{\begin{equation}}
\newcommand{\ee}{\end{equation}}
\newcommand{\bea}{\begin{eqnarray}}
\newcommand{\eea}{\end{eqnarray}}
\newcommand{\bc}{\begin{center}}
\newcommand{\ec}{\end{center}}

\newtheorem{thm}{Theorem}[section]

\newtheorem{lem}[thm]{Lemma}
\newtheorem{prop}[thm]{Proposition}

\theoremstyle{definition}
\newtheorem{defn}[thm]{Definition}

\theoremstyle{remark}
\newtheorem{rem}[thm]{\rm\bfseries{Remark}}

\newtheorem{exm}[thm]{\rm\bfseries{Example}}

\theoremstyle{definition}
\newtheorem{qdefi}[thm]{``Definition''\!}
\newtheorem*{caveat*}{Caveat}

\numberwithin{equation}{section}

\theoremstyle{remark}
\newtheorem*{remark*}{\rm\bfseries{Remark}}

\newcommand{\Z}{\mathbb{Z}}
\newcommand{\R}{\mathbb{R}}
\newcommand{\C}{\mathbb{C}}

\newcommand{\Sym}{\mathrm{Sym}}
\renewcommand{\Re}{\mathrm{Re}\,}
\renewcommand{\Im}{\mathrm{Im}\,}
\newcommand{\F}{\mathcal{F}}
\newcommand{\A}{\mathcal{A}}
\title{Fukaya categories of symmetric products and bordered 
Heegaard-Floer homology}
\author{Denis Auroux}
\address{
Department of Mathematics, M.I.T., Cambridge MA 02139-4307, USA\newline
\indent Department of Mathematics, UC Berkeley, Berkeley CA 94720-3840,
USA}
\email{auroux@math.berkeley.edu, auroux@math.mit.edu}
\thanks{This work was partially supported by NSF grants DMS-0600148 and 
DMS-0652630.}

\begin{document}
\begin{abstract}
The main goal of this paper is to discuss a symplectic interpretation of
Lipshitz, Ozsv\'ath and Thurston's bordered Heegaard-Floer homology
\cite{LOT} in terms of Fukaya categories of symmetric products and Lagrangian
correspondences. More specifically, we give a description of the algebra
$\A(F)$ which appears in the work of Lipshitz, Ozsv\'ath and Thurston
in terms of (partially wrapped) Floer
homology for product Lagrangians in the symmetric product, and outline 
how bordered Heegaard-Floer homology itself can conjecturally be understood
in this language. 
\end{abstract}


\maketitle

\section{Introduction}

Lipshitz, Ozsv\'ath and Thurston's {\it bordered Heegaard-Floer homology}
\cite{LOT} extends the hat version of Heegaard-Floer homology to an invariant for 3-manifolds
with parametrized boundary. Their construction associates to
a (marked and parametrized) surface $F$ a certain algebra
$\A(F)$, and to a 3-manifold with boundary $F$ a pair of 
($A_\infty$-)\,modules over $\A(F)$, which satisfy a TQFT-like gluing
theorem. On the other hand, recent work of Lekili and Perutz \cite{LP}
suggests another construction, whereby a 3-manifold with boundary yields
an object in (a variant of) the Fukaya category of the symmetric product
of $F$. 


\subsection{Lagrangian correspondences and Heegaard-Floer homology}
\label{ss:1.1}

Given a closed 3-manifold $Y$, the Heegaard-Floer homology group
$\widehat{HF}(Y)$ is classically constructed by Ozsv\'ath and
Szab\'o from a Heegaard decomposition by considering
the Lagrangian Floer homology of two product tori in the symmetric
product of the punctured Heegaard surface. Here is an alternative
description of this invariant.

Equip $Y$ with a Morse function 
(with only one minimum and one maximum, and with distinct critical 
values). Then the complement $Y'$ of a ball in $Y$ (obtained by deleting
a neighborhood of a Morse trajectory from the maximum to the minimum)
can be decomposed into a succession of elementary cobordisms $Y'_i$
($i=1,\dots,r$) between connected Riemann surfaces with boundary 
$\Sigma_0,\Sigma_1,\dots,\Sigma_r$ (where $\Sigma_0=\Sigma_r=D^2$, and
the genus increases or decreases by 1 at each step). By a construction
of Perutz \cite{Perutz}, each $Y'_i$ determines a Lagrangian correspondence
$L_i\subset \Sym^{g_{i-1}}(\Sigma_{i-1})\times \Sym^{g_i}(\Sigma_i)$ between
symmetric products. The {\it quilted Floer homology} of the sequence
$(L_1,\dots,L_r)$, as defined by Wehrheim and Woodward 
\cite{WW,WW2}, is then isomorphic to $\widehat{HF}(Y)$. 
(This relies on two
results from the work in progress of Lekili and Perutz \cite{LP}: the
first one concerns the invariance of this quilted Floer homology
under exchanges of critical points, which allows one to reduce to the
case where the genus first increases from $0$ to $g$ then decreases back
to $0$; the second one states that the composition
of the Lagrangian correspondences from $\Sym^0(D^2)$ to $\Sym^g(\Sigma_g)$
is then Hamiltonian isotopic to the product torus considered by Ozsv\'ath and
Szab\'o.)

Given a 3-manifold $Y$ with boundary $\partial Y
\simeq F\cup_{S^1} D^2$ (where $F$ is a connected genus $g$ surface
with one boundary component), we can similarly view $Y$ as a succession of
elementary cobordisms (from $D^2$ to $F$), and 
hence associate to it a sequence of Lagrangian correspondences
$(L_1,\dots,L_r)$. This
defines an object $\mathbb{T}_Y$ of the {\it extended Fukaya category}
$\mathcal{F}^\sharp(\Sym^g(F))$, as defined by Ma'u, Wehrheim and Woodward 
\cite{MWW} (see \cite{WW,WW2} for the cohomology level version).

More generally, we can consider a cobordism between two connected surfaces 
$F_1$ and $F_2$ (each with one boundary component), i.e.,
a 3-manifold $Y_{12}$ with connected boundary, together with
a decomposition $\partial Y_{12}\simeq -F_1\cup_{S^1} F_2$. 
The same construction associates to such $Y$ a generalized Lagrangian 
correspondence 
(i.e., a sequence of correspondences) from $\Sym^{k_1}(F_1)$ to 
$\Sym^{k_2}(F_2)$, whenever $k_2-k_1=g(F_2)-g(F_1)$; by Ma'u, Wehrheim
and Woodward's formalism, such a correspondence defines an 
$A_\infty$-functor from
$\F^\sharp(\Sym^{k_1}(F_1))$ to $\F^\sharp(\Sym^{k_2}(F_2))$. 

To summarize, this suggests that we should associate:

\begin{itemize}
\item to a genus $g$ surface $F$ (with one boundary),
the collection of extended Fukaya categories of its symmetric products,
$\F^\sharp(\Sym^k(F))$ for $0\le k\le 2g$;\smallskip
\item to a 3-manifold $Y$ with boundary $\partial Y\simeq F\cup_{S^1} D^2$,
an object of $\F^\sharp(\Sym^g(F))$ (namely, the generalized Lagrangian
$\mathbb{T}_Y$);\smallskip
\item to a cobordism $Y_{12}$ with boundary $\partial Y_{12}\simeq
-F_1\cup_{S^1} F_2$, a collection of $A_\infty$-functors from
$\F^\sharp(\Sym^{k_1}(F_1))$ to $\F^\sharp(\Sym^{k_2}(F_2))$.
\end{itemize}

These objects behave naturally under gluing: for example, if a closed
$3$-manifold decomposes as $Y=Y_1\cup_{F\cup D^2} Y_2$, where $\partial
Y_1=F\cup D^2=-\partial Y_2$, then we have a quasi-isomorphism
\begin{equation}\label{eq:Tpairing}
\hom_{\F^\sharp(\Sym^g(F))}(\mathbb{T}_{Y_1},\mathbb{T}_{-Y_2})\simeq
\widehat{CF}(Y).\smallskip\end{equation}

Our main goal is to relate this construction to
bordered Heegaard-Floer homology. More precisely, our main results
concern the relation between the algebra $\A(F)$ introduced in 
\cite{LOT} and the Fukaya category of $\Sym^g(F)$. For 3-manifolds with
boundary, we also propose
(without complete proofs) a dictionary between the $A_\infty$-module
$\widehat{CFA}(Y)$ of \cite{LOT} and the generalized Lagrangian submanifold
$\mathbb{T}_Y$ introduced above.

\begin{remark*} The cautious reader should be aware of the following issue
concerning the choice of a symplectic form on $\Sym^g(F)$.
We can equip $F$ with an exact area form, and choose exact Lagrangian
representatives of all the simple closed curves that appear in
Heegaard diagrams. By Corollary 7.2 in \cite{PerHH}, the symmetric 
product $\Sym^g(F)$ carries an 
exact K\"ahler form for which the relevant product tori are exact
Lagrangian. Accordingly, a sizeable portion of this paper, namely all the
results which do not involve correspondences, can be understood in the
exact setting. However, Perutz's construction of Lagrangian correspondences
requires the K\"ahler form to be deformed by
a negative multiple of the first Chern class (cf.\ Theorem A of
\cite{Perutz}). Bubbling is not an issue in any case, because
the symmetric product of $F$ does not contain any closed holomorphic 
curves (also, we can arrange for all Lagrangian submanifolds and
correspondences to be {\it balanced} and in particular monotone).
Still, we will occasionally need to ensure that our results hold for
the perturbed K\"ahler form on $\Sym^g(F)$ and not just in the exact case. 
\end{remark*}

\subsection{Fukaya categories of symmetric products}

Let $\Sigma$ be a double cover of the complex plane branched at $n$ points.
In Section~\ref{s:lf}, we describe the symmetric product $\Sym^k(\Sigma)$
as the total space of a Lefschetz fibration $f_{n,k}$, for any integer
$k\in\{1,\dots,n\}$. The fibration $f_{n,k}$ has $\binom{n}{k}$ critical 
points, and
the Lefschetz thimbles $D_s$ ($s\subseteq \{1,\dots,n\}$, $|s|=k$)
can be understood explicitly as products of arcs on $\Sigma$.

For the purposes of understanding bordered Heegaard-Floer homology,
it is natural to apply these considerations to the case of
the once punctured genus $g$ surface $F$, viewed as a double cover of the
complex plane branched at $2g+1$ points.
However, the algebra $\mathcal{A}(F,k)$ considered by Lipshitz, Ozsv\'ath
and Thurston only has $\binom{2g}{k}$ primitive idempotents \cite{LOT},
whereas our Lefschetz fibration has $\binom{2g+1}{k}$ critical points.

In Section \ref{s:A12}, we consider a somewhat easier case, namely
that of a twice punctured genus $g-1$ surface $F'$, viewed as a
double cover of the complex plane branched at $2g$ points.
We also introduce a subalgebra $\mathcal{A}_{1/2}(F',k)$
of $\mathcal{A}(F,k)$, consisting of collections of Reeb chords on a
matched {\it pair} of pointed circles,
and show that it has a natural interpretation in terms of the
Fukaya category of the Lefschetz fibration $f_{2g,k}$
as defined by Seidel \cite{SeVCM,SeBook}:

\begin{thm}\label{thm:A12}
$\mathcal{A}_{1/2}(F',k)$ is isomorphic to the endomorphism algebra of the
exceptional collection $\{D_s,\ s\subseteq \{1,\dots,2g\},\ |s|=k\}$
in the Fukaya category $\F(f_{2g,k})$. 
\end{thm}

By work of Seidel \cite{SeBook}, the thimbles $D_s$ generate the Fukaya
category $\F(f_{2g,k})$; hence we obtain a derived equivalence between
$\mathcal{A}_{1/2}(F',k)$ and $\F(f_{2g,k})$.

Next, in Section \ref{s:A} we turn to the case of the genus $g$ surface $F$,
which we now regard as a surface with boundary,
and associate a {\it partially wrapped}
Fukaya category $\F_z$ to the pair
$(\Sym^k(F),\{z\}\times \Sym^{k-1}(F))$ where $z$ is a marked
point on the boundary of $F$ (see Definition \ref{def:Fz}). Viewing $F'$ as a subsurface of $F$,
we specifically consider the same collection of $\binom{2g}{k}$ product 
Lagrangians $D_s$, $s\subseteq \{1,\dots,2g\}$, $|s|=k$ as in Theorem \ref{thm:A12}.
Then we have:

\begin{thm}\label{thm:A}
$\mathcal{A}(F,k)\simeq
\bigoplus\limits_{s,s'}\hom_{\F_z}(D_s,D_{s'}).$
\end{thm}

\noindent
As we will explain in Section \ref{ss:moreA}, a similar result
also holds when the algebra $\A(F,k)$ is defined using a different matching
than the one used throughout the paper.

Our next result concerns the structure of the $A_\infty$-category $\F_z$.

\begin{thm}\label{thm:generate}
The partially wrapped Fukaya category $\F_z$ is generated by the $\binom{2g}{k}$
objects $D_s$, $s\subseteq \{1,\dots,2g\}$, $|s|=k$. In particular,
the natural functor from the category of
$A_\infty$-modules over $\F_z$ to that of $\mathcal{A}(F,k)$-modules
is an equivalence.
\end{thm}

Moreover, the same result still holds if we enlarge the category $\F_z$ to
include compact closed ``generalized Lagrangians'' (i.e., sequences of
Lagrangian correspondences) of the sort that arose in the previous section.

\begin{caveat*}
As we will see in Section \ref{s:generate}, this result relies on the
existence of a ``partial wrapping'' (or ``acceleration'') $A_\infty$-functor 
from the Fukaya category of $f_{2g+1,k}$ to $\F_z$, and requires a detailed 
understanding of the relations between various flavors of Fukaya categories.
This would be best achieved in the context of a more systematic study 
of partially wrapped Floer theory, as opposed to the {\it ad hoc} approach
used in this paper (where, in particular, transversality issues are not 
addressed in full generality).
In~\S \ref{s:generate} we sketch a construction of the
acceleration functor in our setting, but do not give full details; we also do
not show that the functor is well-defined and cohomologically unital.
These properties should follow without major difficulty from the 
techniques introduced by Abouzaid and Seidel, but a careful argument
would require a lengthy technical discussion
which is beyond the scope of this paper; thus, the cautious reader should
be warned that the proof of Theorem~\ref{thm:generate} given here is not
quite complete.
\end{caveat*}

\subsection{Yoneda embedding and $\widehat{CFA}$}
Let $Y$ be a 3-manifold with parameterized boundary 
$\partial Y\simeq F\cup_{S^1} D^2$. Following \cite{LOT}, the
manifold $Y$ can be described by a bordered Heegaard diagram, i.e.\
a surface $\Sigma$ of genus $\bar{g}\ge g$ with one boundary
component, carrying:
\begin{itemize}
\item $\bar{g}-g$ simple closed curves
$\alpha_1^c,\dots,\alpha^c_{\bar{g}-g}$, and
$2g$ arcs $\alpha_1^a,\dots,\alpha_{2g}^a$; 
\item $\bar{g}$ simple closed curves $\beta_1,\dots,\beta_{\bar{g}}$;
\item a marked point $z\in\partial\Sigma$.
\end{itemize}
As usual, the $\beta$-curves determine a product torus
$T_\beta=\beta_1\times\dots\times\beta_{\bar{g}}$ inside
$\Sym^{\bar{g}}(\Sigma)$. As to the closed $\alpha$-curves, using Perutz's
construction they determine a Lagrangian correspondence
$T_{\alpha}$ from $\Sym^g(F)$ to $\Sym^{\bar{g}}(\Sigma)$ (or, equivalently,
$\bar{T}_\alpha$ from $\Sym^{\bar{g}}(\Sigma)$ to $\Sym^g(F)$).
The object $\mathbb{T}_Y$ of the extended Fukaya category $\F^\sharp(\Sym^g(F))$
introduced in \S \ref{ss:1.1} is then isomorphic to the formal composition
of $T_\beta$ and $\bar{T}_\alpha$.

There is a contravariant Yoneda-type $A_\infty$-functor $\mathcal{Y}$ from 
the extended Fukaya category of $\Sym^g(F)$ to the category of right 
$A_\infty$-modules over $\A(F,g)$. Indeed, 
$\F^\sharp(\Sym^g(F))$ can be enlarged into a partially wrapped
$A_\infty$-category $\F^\sharp_z$ by adding to it the same non-compact objects
(products of properly embedded arcs) as in $\F_z$.
This allows us to associate to a generalized Lagrangian
$\mathbb{L}$ the $A_\infty$-module $$\mathcal{Y}(\mathbb{L})=
{\textstyle \bigoplus_s} \hom_{\F^\sharp_z}(\mathbb{L},D_s),$$ where the module
maps are given by products in the partially wrapped Fukaya category.
With this understood,
the right $A_\infty$-module constructed
by Lipshitz, Ozsv\'ath and Thurston \cite{LOT} is simply the image of
$\mathbb{T}_Y$ under the Yoneda functor $\mathcal{Y}$:

\begin{thm}\label{th:CFA}
$\widehat{CFA}(Y)\simeq \mathcal{Y}(\mathbb{T}_Y)$.
\end{thm}

\noindent
Since the Lagrangian correspondence $T_\alpha$ maps $D_s$ to
$$T_\alpha(D_s):=\alpha^c_1\times\dots\times\alpha^c_{\bar{g}-g}\times
\prod_{i\in s}\alpha^a_i\subset \Sym^{\bar{g}}(\Sigma),$$
a more down-to-earth formulation of Theorem \ref{th:CFA} is:
$$\widehat{CFA}(Y)\simeq {\textstyle \bigoplus_s
CF^*(T_\beta,T_\alpha(D_s)).}$$
However the module structure is less apparent in this formulation.

Consider now a closed 3-manifold $Y$ which decomposes as the union $Y_1\cup_{F\cup
D^2} Y_2$ of
two manifolds with $\partial Y_1=F\cup_{S^1} D^2=-\partial Y_2$. Then 
we have:

\begin{thm}\label{th:pairing}
$\hom_{\A(F,g)\text{-mod}}(
\widehat{CFA}(-Y_2),\widehat{CFA}(Y_1))$ is quasi-isomorphic to
$\widehat{CF}(Y)$.
\end{thm}

\noindent
This statement is equivalent to the pairing theorem in \cite{LOT} via
a duality property relating $\widehat{CFA}(-Y_2)$ to $\widehat{CFD}(Y_2)$
which is known to Lipshitz, Ozsv\'ath and Thurston. Thus, it should be
viewed not as a new result, but rather as a different insight into the 
main result in \cite{LOT} (see also \cite{AuICM} and \cite{LOTnew} for
recent developments). Observe that the formulation given here does 
not involve $\widehat{CFD}$; this is advantageous since, even though the 
two types of modules contain equivalent information, $\widehat{CFA}$ 
is much more natural from our perspective.

\begin{caveat*}
While the main ingredients in the proofs of Theorems \ref{th:CFA} and 
\ref{th:pairing} are presented in Section~\ref{s:CFApairing},
much of the technology on which the arguments rely is still being
developed; therefore, full proofs are well beyond the scope of this paper.
In particular, the argument for Theorem \ref{th:CFA} relies heavily on
Lekili and Perutz's recent work \cite{LP}, and on the properties
of $A_\infty$-functors associated to Lagrangian correspondences \cite{MWW},
neither of which have been fully written up yet. The cautious reader should
also note that the argument given for Proposition \ref{prop:CFA}
uses a description of the
degeneration of strip-like ends to Morse trajectories as the Hamiltonian 
perturbations tend to zero which, to our knowledge, has not been written 
up in detail anywhere in the form needed here. Finally, we point out that,
while in our approach Theorem \ref{th:pairing} is obtained
as a corollary of Theorems \ref{thm:generate} and \ref{th:CFA}, a direct
proof of this result has recently been obtained by Lipshitz, Ozsv\'ath and
Thurston \cite{LOTnew} purely within the framework of bordered
Heegaard-Floer theory.
\end{caveat*}

%

\subsection*{Acknowledgements}
I am very grateful to Mohammed Abouzaid, Robert Lipshitz,  Peter Ozsv\'ath,
Tim Perutz, Paul Seidel and Dylan Thurston, whose many helpful suggestions and 
comments influenced this work in decisive ways. In particular, I am heavily
indebted to Mohammed Abouzaid for his patient explanations of wrapped
Fukaya categories and for suggesting the approach outlined in the
appendix. Finally, I would like to thank the referees for valuable
comments.
This work was partially supported by
NSF grants DMS-0600148 and DMS-0652630.

\section{A Lefschetz fibration on $\Sym^k(\Sigma)$}\label{s:lf}

Fix an ordered sequence of $n$ real numbers $\theta_1<\theta_2<\dots<\theta_n$,
and consider the points $p_j=i\theta_j$ on the imaginary axis in the 
complex plane. Let $\Sigma$ be the
double cover of $\C$ branched at $p_1,\dots,p_n$: hence $\Sigma$ is
a Riemann surface of genus $\lfloor \frac{n-1}{2}\rfloor $ with one 
(resp.\ two) puncture(s) if $n$ is odd (resp.\ even). We denote by 
$\pi:\Sigma\to\C$ the covering map, and let $q_j=\pi^{-1}(p_j)\in\Sigma$.

We consider the $k$-fold symmetric product of the Riemann surface
$\Sigma$ ($1\le k\le n$), with the product complex structure $J$, and the
holomorphic map $f_{n,k}:\Sym^k(\Sigma)\to \C$ defined by
$f_{n,k}([z_1,\dots,z_k])=\pi(z_1)+\dots+\pi(z_k)$.

\begin{prop}\label{prop:lf}
$f_{n,k}:\Sym^k(\Sigma)\to \C$ is a Lefschetz fibration,
whose $\binom{n}{k}$ critical points are the tuples consisting
of $k$ distinct points in $\{q_1,\dots,q_{n}\}$.
\end{prop}

\proof
Given $\underline{z}\in\Sym^k(\Sigma)$, denote by
$z_1,\dots,z_r$ the distinct elements in the $k$-tuple
$\underline{z}$, and by $k_1,\dots,k_r$ the multiplicities with which
they appear. The tangent space $T_{\underline{z}}\Sym^k(\Sigma)$ decomposes
into the direct sum of the $T_{[z_i,\dots,z_i]}\Sym^{k_i}(\Sigma)$,
and $df_{n,k}(\underline{z})$ splits into the direct sum of the
differentials $df_{n,k_i}([z_i,\dots,z_i])$. Thus $\underline{z}$ is a
critical point of $f_{n,k}$ if and only if $[z_i,\dots,z_i]$ is a critical
point of $f_{n,k_i}$ for each $i\in \{1,\dots,r\}$.

By considering the 
restriction of $f_{n,k_i}$ to the diagonal stratum, we see that
$[z_i,\dots,z_i]$ cannot be a critical point of $f_{n,k_i}$ unless $z_i$
is a critical point of $\pi$. Assume now that $z_i$ is a critical point
of $\pi$, and pick a local complex coordinate $w$ on $\Sigma$ near $z_i$,
in which $\pi(w)=w^2+\mathrm{constant}$. Then a neighborhood
of $[z_i,\dots,z_i]$ in $\Sym^{k_i}(\Sigma)$ identifies with a 
neighborhood of the origin in $\Sym^{k_i}(\C)$, with coordinates
given by the elementary symmetric functions $\sigma_1,\dots,\sigma_{k_i}$.
The local model for $f_{n,k_i}$
is then
$$f_{n,k_i}([w_1,\dots,w_{k_i}])=w_1^2+\dots+w_{k_i}^2+\mathrm{constant}=
\sigma_1^2-2\sigma_2+\mathrm{constant}.$$ Thus, for $k_i\ge 2$ the point $[z_i,\dots,z_i]$
is never a critical point of $f_{n,k_i}$. We conclude that the only critical
points of $f_{n,k}$ are tuples of distinct critical points of $\pi$; moreover
these critical points are clearly non-degenerate.
\endproof

We denote by $\mathcal{S}^n_k$ the set of all $k$-element subsets of
$\{1,\dots,n\}$, and for $s\in \mathcal{S}^n_k$ we call $\vec{q}_s$
the critical point $\{q_j,\ j\in s\}$ of $f_{n,k}$.

We equip $\Sigma$ with an area form $\sigma$, and equip $\Sym^k(\Sigma)$ 
with an exact K\"ahler form~$\omega$ that coincides with the product
K\"ahler form on $\Sigma^k$ away from the diagonal strata (see e.g.\
Corollary 7.2 in \cite{PerHH}). 
The K\"ahler form $\omega$ defines a symplectic horizontal distribution 
on the fibration $f_{n,k}$ away from its critical points, given
by the symplectic orthogonal to the fibers. Because $f_{n,k}$ is holomorphic,
this horizontal distribution is spanned by the gradient vector fields
for $\Re f_{n,k}$ and $\Im f_{n,k}$ with respect to the K\"ahler metric
$g=\omega(\cdot,J\cdot)$.

Given a critical point $\vec{q}_s$ of $f_{n,k}$ and an embedded arc 
$\gamma$ in $\C$ connecting $f_{n,k}(\vec{q}_s)$ to infinity,
the {\em Lefschetz thimble} associated to $\vec{q}_s$ and $\gamma$ is
the properly embedded Lagrangian disc consisting of all points in
$f_{n,k}^{-1}(\gamma)$ whose parallel transport along $\gamma$ converges
to the critical point $\vec{q}_s$ \cite{SeVCM,SeBook}. In our case,
we take $\gamma$ to be the straight line 
$\gamma(\theta_s)=\R_{\ge 0} + i\theta_s$, where 
$\theta_s=\Im f_{n,k}(\vec{q}_s)=\sum_{j\in s} \theta_j$,
and we denote by $D_s\subset \Sym^k(\Sigma)$ the corresponding Lefschetz thimble.

The thimbles $D_s$ have a simple description in terms of the
disjoint properly embedded arcs $\alpha_j=\pi^{-1}(\R_{\ge 0} + i\theta_j)
\subset \Sigma$. Namely:

\begin{lem}
$D_s=\prod\limits_{j\in s} \alpha_j$.
\end{lem}

\proof
Since $\gamma_s$ is parallel to the real axis, parallel transport is given
by the gradient flow of $\Re f_{n,k}$ with respect to the K\"ahler metric $g$.
Away from the diagonal strata, $g$ is a product metric, and so
the components of the gradient vector of $\Re f_{n,k}$ at $[z_1,\dots,z_k]$
are $\nabla\Re \pi(z_1),\dots,\nabla\Re \pi(z_k)$. Thus parallel transport
along $\gamma_s$ decomposes into the product of the parallel transports along
the arcs $\R_{\ge 0}+i\theta_j$.
\endproof

In the subsequent discussion, we will also need to consider perturbed
versions of the thimbles $D_s$. Fix a positive real number $\epsilon$.
Given $\theta\in \R$, we consider the arc
$\gamma^\pm(\theta)=\{i\theta+(1\mp i\epsilon)t,\ t\ge 0\}$ in the complex
plane, connecting $i\theta$ to infinity. For $s\in \mathcal{S}^n_k$ we denote 
by $D^\pm_s\subset \Sym^k(\Sigma)$ the thimble associated to the arc 
$\gamma^\pm(\theta_s)$, and for $j\in \{1,\dots,n\}$ we set
$\alpha^\pm_j=\pi^{-1}(\gamma^\pm(\theta_j))\subset \Sigma$ (see
Figure \ref{fig:alphas}). The same argument as above then gives:

\begin{figure}[t]
\setlength{\unitlength}{7.5mm}
\begin{picture}(5.5,3.1)(-1,-1.5)
\psset{unit=\unitlength}
\psellipse[linewidth=0.5pt,linestyle=dashed,dash=2pt 2pt](3,0.85)(0.2,0.7)
\psellipticarc(3,0.85)(0.2,0.7){-90}{90}
\psellipse[linewidth=0.5pt,linestyle=dashed,dash=2pt 2pt](3,-0.85)(0.2,0.7)
\psellipticarc(3,-0.85)(0.2,0.7){-90}{90}
\psellipticarc(3,0)(0.8,0.18){90}{270}
\psellipse(1,0)(0.5,0.2)
\psline[linearc=1.6](3,1.55)(-0.5,1.55)(-0.5,-1.55)(3,-1.55)
\psellipticarc(3.1,0)(0.9,0.35){90}{180}
\psellipticarc(3.15,0)(1.66,0.6){90}{180}
\psellipticarc(3.18,0)(2.7,0.9){90}{180}
\psellipticarc(3.13,0)(3.63,1.2){90}{180}
\psellipticarc[linestyle=dashed,dash=2pt 2pt](2.88,0)(0.68,0.35){180}{270}
\psellipticarc[linestyle=dashed,dash=2pt 2pt](2.83,0)(1.34,0.6){180}{270}
\psellipticarc[linestyle=dashed,dash=2pt 2pt](2.80,0)(2.32,0.9){180}{270}
\psellipticarc[linestyle=dashed,dash=2pt 2pt](2.85,0)(3.35,1.2){180}{270}
\put(3.29,0.2){\tiny $\alpha_1$}
\put(3.25,1.15){\tiny $\alpha_n$}
\psline[linestyle=dotted](3.4,0.95)(3.4,0.55)
\psline{->}(4,0)(4.8,0)
\put(4.25,0.2){\tiny $\pi$}
\put(4.1,-0.4){\tiny 2:1}
\end{picture}
\qquad
\begin{picture}(5,3)(-1,-1.5)
\psset{unit=\unitlength}
\psframe(-1,-1.6)(3.2,1.6)
\pscircle*(0,-0.6){0.05}
\pscircle*(0,-0.2){0.05}
\pscircle*(0,0.2){0.05}
\pscircle*(0,0.6){0.05}
\psline(0,-0.6)(3.14,-0.6)
\psline(0,-0.2)(3.16,-0.2)
\psline(0,0.2)(3.16,0.2)
\psline(0,0.6)(3.14,0.6)
\put(-0.6,-0.7){\tiny $p_1$}
\put(-0.6,0.55){\tiny $p_{n}$}
\put(3.25,-0.7){\tiny $\pi(\alpha_1)$}
\put(3.25,0.6){\tiny $\pi(\alpha_{n})$}
\psline[linestyle=dotted](3.6,0.2)(3.6,-0.2)
\psline[linestyle=dotted](-0.5,-0.2)(-0.5,0.2)
\end{picture}
\qquad
\begin{picture}(5,3)(-1,-1.5)
\psset{unit=\unitlength}
\psframe(-1,-1.6)(3,1.6)
\pscircle*(0,-0.6){0.05}
\pscircle*(0,-0.2){0.05}
\pscircle*(0,0.2){0.05}
\pscircle*(0,0.6){0.05}
\psline(0,-0.6)(3,-1.5)
\psline(0,-0.2)(3,-1.1)
\psline(0,0.2)(3,-0.7)
\psline(0,0.6)(3,-0.3)
\psline(0,0.6)(3,1.5)
\psline(0,0.2)(3,1.1)
\psline(0,-0.2)(3,0.7)
\psline(0,-0.6)(3,0.3)
\put(-0.6,-0.7){\tiny $p_1$}
\put(-0.6,0.55){\tiny $p_{n}$}
\put(3.1,-1.55){\tiny $\pi(\alpha^+_1)$}
\put(3.1,-0.35){\tiny $\pi(\alpha^+_{n})$}
\put(3.1,0.2){\tiny $\pi(\alpha^-_1)$}
\put(3.1,1.4){\tiny $\pi(\alpha^-_{n})$}
\psline[linestyle=dotted](3.4,-1.05)(3.4,-0.65)
\psline[linestyle=dotted](3.4,0.7)(3.4,1.1)
\psline[linestyle=dotted](-0.5,-0.2)(-0.5,0.2)
\end{picture}
\caption{The arcs $\alpha_j$ and $\alpha^\pm_j$}\label{fig:alphas}
\end{figure}
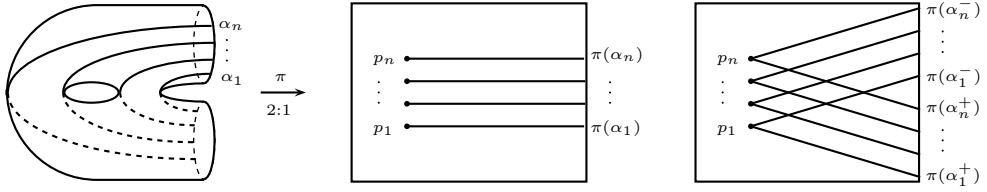

\begin{lem}\label{l:UT} 
$D_s^\pm=\prod\limits_{j\in s}\alpha^\pm_j.$
\end{lem}

\section{The algebra $\A_{1/2}(F',k)$ and the Fukaya category 
of $f_{2g,k}$}\label{s:A12}

\subsection{The algebra $\A_{1/2}(F',k)$}\label{ss:strands}
We start by briefly recalling the definition of the differential
algebra $\A(F,k)$ associated to a genus $g$ surface $F$ with one boundary;
the reader is referred to \cite[\S 3]{LOT} for details.
Consider $4g$ points $a_1,\dots,a_{4g}$ along an oriented segment 
(thought of as the complement of a marked point in an oriented circle),
carrying the labels $1,\dots,2g,1,\dots,2g$ (we fix this specific matching
throughout). The generators of $\A(F,k)$
are unordered $k$-tuples consisting of two types of items:
\begin{itemize}
\item ordered pairs $(i,j)$ with $1\le i<j\le 4g$, corresponding to Reeb
chords connecting pairs of points on the marked circle; in the notation
of \cite{LOT} these are denoted by a column
$\left[\begin{smallmatrix}i\\j\end{smallmatrix}\right]$, or graphically by an 
upwards strand connecting the $i$-th point to the $j$-th point;
\item unordered pairs $\{i,j\}$ such that $a_i$ and $a_j$ carry the same
label (i.e., in our case, $i$ and $j$ differ by $2g$), denoted by a
column $\left[\begin{smallmatrix}i\\ \vphantom{j}\end{smallmatrix}\right]$, or
graphically by two horizontal dotted lines.
\end{itemize}
The $k$ source labels (i.e., the labels of the initial points) are
moreover required to be all distinct, and similarly for the $k$ target labels.
We will think of $\mathcal{A}(F,k)$ as a finite category with objects
indexed by $k$-element subsets of $\{1,\dots,2g\}$,
where, given $s,t\in \mathcal{S}_k:=\mathcal{S}^{2g}_k$,
$\hom(s,t)$ is the linear span of the generators with source labels the
elements of $s$ and target labels the elements of $t$. For instance, taking
$g=k=2$, the generator 
\begin{center}
\setlength{\unitlength}{3mm}
\begin{picture}(7,7)(-5,1)
\put(-1.8,4.5){\makebox(0,0)[rc]{$\left[\begin{smallmatrix}5 & 2 \\[2pt]
8 & \end{smallmatrix}\right]=$}}
\psset{unit=\unitlength}
\multiput(0,1)(0,1){8}{\circle*{0.2}}
\multiput(2,1)(0,1){8}{\circle*{0.2}}
\put(-0.45,1){\makebox(0,0)[rc]{\tiny 1}}
\put(-0.45,2){\makebox(0,0)[rc]{\tiny 2}}
\put(-0.45,3){\makebox(0,0)[rc]{\tiny 3}}
\put(-0.45,4){\makebox(0,0)[rc]{\tiny 4}}
\put(-0.45,5){\makebox(0,0)[rc]{\tiny 5}}
\put(-0.45,6){\makebox(0,0)[rc]{\tiny 6}}
\put(-0.45,7){\makebox(0,0)[rc]{\tiny 7}}
\put(-0.45,8){\makebox(0,0)[rc]{\tiny 8}}
\put(2.45,1){\makebox(0,0)[lc]{\tiny 1}}
\put(2.45,2){\makebox(0,0)[lc]{\tiny 2}}
\put(2.45,3){\makebox(0,0)[lc]{\tiny 3}}
\put(2.45,4){\makebox(0,0)[lc]{\tiny 4}}
\put(2.45,5){\makebox(0,0)[lc]{\tiny 5}}
\put(2.45,6){\makebox(0,0)[lc]{\tiny 6}}
\put(2.45,7){\makebox(0,0)[lc]{\tiny 7}}
\put(2.45,8){\makebox(0,0)[lc]{\tiny 8}}
\psline[linestyle=dashed,dash=2pt 2pt](0,2)(2,2)
\psline[linestyle=dashed,dash=2pt 2pt](0,6)(2,6)
\pscurve(0,5)(0.2,5)(1.8,8)(2,8)
\end{picture}
\end{center}
is viewed as a morphism from $\{1,2\}$ to $\{2,4\}$.

Composition in $\A(F,k)$ is given by concatenation of strand diagrams,
provided that no two strands of the concatenated diagram
cross more than once; otherwise the product is zero \cite{LOT}.
(Of course, the product also vanishes if the target and source labels fail
to match up.) The primitive idempotents of $\A(F,k)$ correspond
to diagrams consisting only of dotted lines, which are the identity
endomorphisms of the various objects. Finally, the differential in $\A(F,k)$
is described graphically as the sum of all the ways of resolving one
crossing of the strand diagram (again excluding resolutions
in which two strands intersect twice). In these operations, a pair of
dotted lines should be treated as the sum of the corresponding arcs.
For example,
\begin{equation}
\label{eq:exdiff}
\partial \left[\begin{smallmatrix}5 & 2 \\ 8 & \end{smallmatrix}\right]=
\left[\begin{smallmatrix}5 & 6 \\ 6 & 8\end{smallmatrix}\right].
\end{equation}

\begin{defn}
We define $\A_{1/2}(F',k)$ to be the subalgebra of $\A(F,k)$ generated by
the strand diagrams for which no strand crosses the interval $[2g,2g+1]$.
\end{defn}

\noindent
(This definition makes sense, as $\A_{1/2}(F',k)$ is clearly
closed under both the differential and the product of $\A(F,k)$.)

\begin{rem}
It is useful to think of $\A_{1/2}(F',k)$ as the algebra associated to a {\it pair}
of pointed circles, one of them carrying the $2g$ points
$a_1,\dots,a_{2g}$ while the other carries $a_{2g+1},\dots,a_{4g}$;
in addition, each of the two circles is equipped with a marked point 
through which Reeb chords are not allowed to pass. 
Connecting two annuli by $2g$ bands in the manner prescribed by
the labels and further attaching a pair of discs yields a
twice punctured genus $g-1$ surface, which we denote by $F'$; as we will
see in the rest of this section, the algebra $\A_{1/2}(F',k)$ can be
understood in terms of the symplectic geometry of this surface
and its symmetric products.
\end{rem}

The algebra $\A_{1/2}(F',k)$ is significantly smaller than $\A(F,k)$: for 
instance, every object of
$\A_{1/2}(F',k)$ is exceptional, i.e.\
$\hom(s,s)=\Z_2\,\mathrm{id}_s$, while there are many more endomorphisms in
$\A(F,k)$. Another feature distinguishing $\A_{1/2}(F',k)$ from $\A(F,k)$ is
directedness. In fact, as will be clear from the rest of this paper, the
relation between $\A_{1/2}(F',k)$ and $\A(F,k)$ is analogous to that between
the directed Fukaya category of a Lefschetz fibration and a partially
wrapped counterpart.

\subsection{The Fukaya category of $f_{n,k}$}\label{ss:fukayalf}
The Fukaya category of the Lefschetz fibration $f_{n,k}$ is a variant of the
Fukaya category of $\Sym^k(\Sigma)$ which allows potentially non-compact
Lagrangian submanifolds as long as they are {\it admissible}, i.e.\
invariant under the gradient flow of $\Re f_{n,k}$ outside of a compact subset.
While the construction finds its roots in ideas of Kontsevich about homological
mirror symmetry for Fano varieties, it has been most extensively studied
by Seidel; see in particular \cite{SeVCM,SeBook}.
In order to make intersection theory for admissible non-compact Lagrangians
well-defined, one needs to choose Hamiltonian perturbations that behave
in a consistent manner near infinity. The description we give here is
slightly different from that in Seidel's work, but can easily be checked
to be equivalent; it is also closely related to the viewpoint given by
Abouzaid in Section 2 of \cite{Ab1}, except we place the base point at
infinity.

Given a real number $\nu$, we say that an exact Lagrangian submanifold 
$L$ of $\Sym^k(\Sigma)$ is
{\em admissible with slope $\nu=\nu(L)$} if the restriction of $f_{n,k}$ to $L$ is
proper and, outside of a compact set, takes values in the half-line
$i\theta+(1+i\nu)\R_+$ for some $\theta\in\R$.
A pair of admissible exact Lagrangians $(L_1,L_2)$ is said to be 
{\em positive} if their slopes satisfy $\nu(L_1)>\nu(L_2)$.

Given two admissible Lagrangians $L_1$ and $L_2$, we can always deform
them by Hamiltonian isotopies (among admissible Lagrangians) to a positive
pair $(\tilde{L}_1,\tilde{L}_2)$. We define
$\hom_{\F(f_{n,k})}(L_1,L_2)=CF^*(\tilde{L}_1,\tilde{L}_2)$,
the Floer complex of the pair $(\tilde{L}_1,\tilde{L}_2)$,
equipped with the Floer differential.
Positivity ensures that the intersections of $\tilde{L}_1$ and
$\tilde{L}_2$ remain in a bounded subset, and the maximum principle
applied to $\Re f_{n,k}$ prevents sequences of holomorphic discs from 
escaping to infinity. Moreover, the Floer cohomology defined in this
manner does not depend on the chosen Hamiltonian isotopies.
The composition $\hom_{\F(f_{n,k})}(L_1,L_2)
\otimes \hom_{\F(f_{n,k})}(L_2,L_3)\to \hom_{\F(f_{n,k})}(L_1,L_3)$ is
similarly defined using the pair-of-pants product in Floer
theory, after replacing each $L_i$ by a Hamiltonian isotopic
admissible Lagrangian $\tilde{L}_i$
in such a way that the pairs $(\tilde{L}_1,\tilde{L}_2)$ and
$(\tilde{L}_2,\tilde{L}_3)$ are both positive; likewise for the
higher compositions.

In order for this construction to be well-defined at the chain
level, in general one needs to specify a procedure for perturbing 
Lagrangians towards positive position. If one considers a collection
of Lefschetz thimbles as will be the case here, then there is a natural
choice, for which the morphisms and $A_\infty$ operations can
be described in terms of Floer theory for the vanishing cycles inside
the fiber of $f_{n,k}$ \cite{SeVCM,SeBook}. (This dimensional reduction is one of the key
features that make Seidel's construction computationally powerful;
however, in the present case it is more efficient to consider
the thimbles rather than the vanishing cycles).

\begin{rem}
We will work over $\Z_2$ coefficients to avoid getting into sign
considerations, and to match with the construction in \cite{LOT};
however, the Lefschetz thimbles $D_s$ are contractible and hence
carry canonical spin structures, which can be used to orient all
the moduli spaces. Keeping track of orientations should give a
procedure for defining the algebras $\mathcal{A}_{1/2}(F',k)$
and $\mathcal{A}(F,k)$ over $\Z$.
\end{rem}

\subsection{Proof of Theorem \ref{thm:A12}}\label{ss:proof12}

We now specialize to the case $n=2g$, and
consider a twice punctured genus $g-1$ surface $F'$ 
(viewed as a double cover of $\C$ branched at $2g$ points),
and the Lefschetz fibration $f_{2g,k}:\Sym^k(F')\to\C$.
Consider two $k$-element subsets $s,t\in \mathcal{S}_k=
\smash{\mathcal{S}^{2g}_k}$,
and the thimbles $D_s,D_t\subset \Sym^k(F')$ defined in Section~\ref{s:lf}. Positivity
can be achieved in a number of manners, e.g.\ we may consider
any of the pairs $(D_s^-,D_t^+)$, $(D_s,D_t^+)$, or $(D_s^-,D_t)$.
We pick the first possibility. By Lemma
\ref{l:UT}, $$D_s^-\cap D_t^+=\biggl(\prod_{i\in s}\alpha_i^-\biggr)\cap
\biggl(\prod_{j\in t}\alpha_j^+\biggr).$$

\begin{prop}\label{prop:complex12}
The chain complexes $\hom_{\F(f_{2g,k})}(D_s,D_t)$ and $\hom_{\A_{1/2}(F',k)}(s,t)$
are isomorphic.
\end{prop}

\proof
The intersections of $D_s^-$ with $D_t^+$ consist of $k$-tuples of 
intersections between the arcs $\alpha_i^-$, $i\in s$ and $\alpha_j^+$,
$j\in t$. These can be determined by looking at Figure~\ref{fig:alphas}. 
Namely, $\alpha_i^-\cap \alpha_j^+$ is
empty if $i>j$, a single point (the branch point $q_i$) if $i=j$, and
a pair of points if $i<j$. 
The preimage $\pi^{-1}(\{\Re z>0\})$ consists of two distinct
components, which we call $V$ and $V'$; then for $i<j$ we call
$q_{i^-j^+}$ (resp.\ $q'_{i^-j^+}$) the point of 
$\alpha_i^-\cap \alpha_j^+$ which lies in $V$ (resp.\ $V'$).

The dictionary between intersection points and generators of
$\hom(s,t)$ is as follows:
\begin{itemize}
\item the point $q_i$ corresponds to the column
$\left[\begin{smallmatrix} i\\\vphantom{j}\end{smallmatrix}\right]$;
\item the point $q_{i^-j^+}$ corresponds to the column
$\left[\begin{smallmatrix} i\\j\end{smallmatrix}\right]$;
\item the point $q'_{i^-j^+}$ corresponds to the  column
$\left[\begin{smallmatrix} 2g+i\\2g+j\end{smallmatrix}\right]$.
\end{itemize}
In both cases, we consider $k$-tuples of such items with the property
that the labels in $s$ and $t$ each appear exactly once; thus we have
a bijection between the generators of $\hom_{\F(f_{2g,k})}(D_s,D_t)$ and 
those of $\hom_{\A_{1/2}(F',k)}(s,t)$.

Next, we consider the Floer differential on
$\hom_{\F(f_{2g,k})}(D_s,D_t)=CF^*(D_s^-,D_t^+)$.
Since the thimbles $D_s^-=\prod_{i\in s} \alpha_i^-$
and $D_t^+=\prod_{j\in t}\alpha_j^+$ are products of arcs in $F'$, 
results from Heegaard-Floer theory can be used in this setting. 
The key observation is that the arcs $\alpha_i^-$ and $\alpha_j^+$ form
a {\em nice} diagram on $F'$, in the sense that the bounded
regions of $F'$ delimited by the arcs $\alpha_i^-$ and $\alpha_j^+$ 
are all rectangles (namely, the preimages of the bounded regions depicted
on Figure \ref{fig:alphas} right). As observed by Sarkar and Wang,
this implies that the Floer differential on $CF^*(D_s^-,D_t^+)$ counts
empty embedded rectangles \cite[Theorems 3.3 and~3.4]{SW}.

(Recall that an {\it embedded rectangle} connecting 
$x\in D_s^-\cap D_t^+$ to $y\in D_s^-\cap D_t^+$ is an
embedded rectangular domain $R$ in the Riemann surface $F'$, satisfying
a local convexity condition, and with boundary
on the arcs that make
up the product Lagrangians $D_s^-$ and $D_t^+$; the two corners where the 
boundary of $R$ jumps from some $\alpha_i^-$ to some $\alpha_j^+$ 
are two of the components of the $k$-tuple $x$,
while the two other corners are components of~$y$. 
The embedded rectangle $R$ is said to be {\it empty} if the other
intersection points which make up the generators $x$ and $y$ all lie outside
of $R$. Because the Maslov index of a holomorphic strip in the 
symmetric product is given by its intersection number with the diagonal,
any index 1 strip
must project to an empty embedded rectangle in $F'$. See
\cite[Section 3]{SW}.)

The embedded rectangles we need to consider lie either within the closure
of $V$ or within that of $V'$; thus they can be understood by looking
at Figure \ref{fig:alphas} right. If a rectangle in $V$ has its sides on
$\alpha^-_i,\alpha^-_j,\alpha^+_l,\alpha^+_m$ ($i<j\le l<m$), then its
``incoming'' vertices are $q_{i^-m^+}$ and either $q_{j^-l^+}$ (if $j<l$)
or $q_j$ (if $j=l$), and its ``outgoing'' vertices are $q_{i^-l^+}$ and
$q_{j^-m^+}$. Via the above dictionary, this corresponds precisely to
resolving the crossing between a strand that connects $a_i$ to $a_m$ and
a strand that connects $a_j$ to $a_l$ (the latter possibly dotted if $j=l$).

The rectangle bounded by $\alpha^-_i,\alpha^-_j,\alpha^+_l,\alpha^+_m$
in $V$ is empty if and only if
the generators under consideration do not include any of the
intersection points $q_{v^-w^+}$ (or equivalently, strands connecting
$a_v$ to $a_w$) with $i<v<j$ and $l<w<m$; this forbidden configuration is
precisely the case in which resolving the crossing would create a double
crossing, which is excluded by the definition of the differential on
$\A_{1/2}(F',k)$.

Empty rectangles in $V'$ can be described similarly
in terms of resolving crossings between strands that connect pairs of points in
$\{a_{2g+1},\dots,a_{4g}\}$. Thus the differential on $CF^*(D_s^-,D_t^+)$
agrees with that on $\hom_{\A_{1/2}(F',k)}(s,t)$.
\endproof

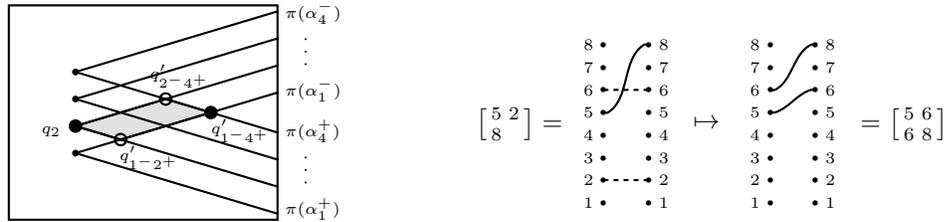
\begin{figure}[b]
\setlength{\unitlength}{9mm}
\begin{picture}(7,3.1)(-1,-1.5)
\psset{unit=\unitlength}
\newrgbcolor{lt2}{0.88 0.88 0.88}
\pspolygon[fillstyle=solid,fillcolor=lt2](0,-0.2)(1.333,0.2)(2,0)(0.667,-0.4)
\psframe(-1,-1.6)(3,1.6)
\pscircle*(0,-0.6){0.05}
\pscircle*(0,-0.2){0.1}
\pscircle*(0,0.2){0.05}
\pscircle*(0,0.6){0.05}
\psline(0,-0.6)(3,-1.5)
\psline(0,-0.2)(3,-1.1)
\psline(0,0.2)(3,-0.7)
\psline(0,0.6)(3,-0.3)
\psline(0,0.6)(3,1.5)
\psline(0,0.2)(3,1.1)
\psline(0,-0.2)(3,0.7)
\psline(0,-0.6)(3,0.3)
\put(3.1,-1.5){\tiny $\pi(\alpha^+_1)$}
\put(3.1,-0.35){\tiny $\pi(\alpha^+_{4})$}
\put(3.1,0.25){\tiny $\pi(\alpha^-_1)$}
\put(3.1,1.4){\tiny $\pi(\alpha^-_{4})$}
\psline[linestyle=dotted](3.4,-1.05)(3.4,-0.65)
\psline[linestyle=dotted](3.4,0.7)(3.4,1.1)
\pscircle*(2,0){0.1}
\pscircle(1.333,0.2){0.1}
\pscircle(0.667,-0.4){0.1}
\put(-0.5,-0.3){\tiny $q_2$}
\put(0.65,-0.7){\tiny $q'_{1^-2^+}$}
\put(2,-0.3){\tiny $q'_{1^-4^+}$}
\put(1.1,0.5){\tiny $q'_{2^-4^+}$}
\end{picture}
\setlength{\unitlength}{3mm}
\begin{picture}(7,7.5)(-5,0.5)
\put(-1.8,4.5){\makebox(0,0)[rc]{$\left[\begin{smallmatrix}5 & 2 \\[2pt]
8 & \end{smallmatrix}\right]=$}}
\psset{unit=\unitlength}
\multiput(0,1)(0,1){8}{\circle*{0.2}}
\multiput(2,1)(0,1){8}{\circle*{0.2}}
\put(-0.45,1){\makebox(0,0)[rc]{\tiny 1}}
\put(-0.45,2){\makebox(0,0)[rc]{\tiny 2}}
\put(-0.45,3){\makebox(0,0)[rc]{\tiny 3}}
\put(-0.45,4){\makebox(0,0)[rc]{\tiny 4}}
\put(-0.45,5){\makebox(0,0)[rc]{\tiny 5}}
\put(-0.45,6){\makebox(0,0)[rc]{\tiny 6}}
\put(-0.45,7){\makebox(0,0)[rc]{\tiny 7}}
\put(-0.45,8){\makebox(0,0)[rc]{\tiny 8}}
\put(2.45,1){\makebox(0,0)[lc]{\tiny 1}}
\put(2.45,2){\makebox(0,0)[lc]{\tiny 2}}
\put(2.45,3){\makebox(0,0)[lc]{\tiny 3}}
\put(2.45,4){\makebox(0,0)[lc]{\tiny 4}}
\put(2.45,5){\makebox(0,0)[lc]{\tiny 5}}
\put(2.45,6){\makebox(0,0)[lc]{\tiny 6}}
\put(2.45,7){\makebox(0,0)[lc]{\tiny 7}}
\put(2.45,8){\makebox(0,0)[lc]{\tiny 8}}
\psline[linestyle=dashed,dash=2pt 2pt](0,2)(2,2)
\psline[linestyle=dashed,dash=2pt 2pt](0,6)(2,6)
\pscurve(0,5)(0.2,5)(1.8,8)(2,8)
\put(4,4.5){\makebox(0,0)[lc]{$\mapsto$}}
\end{picture}
\begin{picture}(12,7.5)(-5,0.5)
\put(4,4.5){\makebox(0,0)[lc]{$=\left[\begin{smallmatrix}5 & 6 \\[2pt]
6 & 8 \end{smallmatrix}\right]$}}
\psset{unit=\unitlength}
\multiput(0,1)(0,1){8}{\circle*{0.2}}
\multiput(2,1)(0,1){8}{\circle*{0.2}}
\put(-0.45,1){\makebox(0,0)[rc]{\tiny 1}}
\put(-0.45,2){\makebox(0,0)[rc]{\tiny 2}}
\put(-0.45,3){\makebox(0,0)[rc]{\tiny 3}}
\put(-0.45,4){\makebox(0,0)[rc]{\tiny 4}}
\put(-0.45,5){\makebox(0,0)[rc]{\tiny 5}}
\put(-0.45,6){\makebox(0,0)[rc]{\tiny 6}}
\put(-0.45,7){\makebox(0,0)[rc]{\tiny 7}}
\put(-0.45,8){\makebox(0,0)[rc]{\tiny 8}}
\put(2.45,1){\makebox(0,0)[lc]{\tiny 1}}
\put(2.45,2){\makebox(0,0)[lc]{\tiny 2}}
\put(2.45,3){\makebox(0,0)[lc]{\tiny 3}}
\put(2.45,4){\makebox(0,0)[lc]{\tiny 4}}
\put(2.45,5){\makebox(0,0)[lc]{\tiny 5}}
\put(2.45,6){\makebox(0,0)[lc]{\tiny 6}}
\put(2.45,7){\makebox(0,0)[lc]{\tiny 7}}
\put(2.45,8){\makebox(0,0)[lc]{\tiny 8}}
\pscurve(0,5)(0.2,5)(1.8,6)(2,6)
\pscurve(0,6)(0.2,6)(1.8,8)(2,8)
\end{picture}

\caption{An empty rectangle and the corresponding differential.}
\label{fig:emptyrec}
\end{figure}

To illustrate the above construction, Figure \ref{fig:emptyrec} shows the
image under $\pi$ of the
empty rectangle (contained in $V'$) which determines (\ref{eq:exdiff}).

Next we need to compare the products in $\F(f_{2g,k})$ and $\A_{1/2}(F',k)$.
Given $s,t,u\in \mathcal{S}_k$, the composition $\hom(D_s,D_t)\otimes \hom(D_t,D_u)
\to \hom(D_s,D_u)$ in $\F(f_{2g,k})$ is defined in terms of perturbations of the
thimbles for which positivity holds: namely, we can consider the Floer
pair-of-pants product $$CF^*(D_s^-,D_t)\otimes CF^*(D_t,D_u^+)\to
CF^*(D_s^-,D_u^+).$$

\begin{prop}\label{prop:prod12}
The isomorphism of Proposition \ref{prop:complex12} intertwines the product
structures of $\F(f_{2g,k})$ and $\A_{1/2}(F',k)$.
\end{prop}

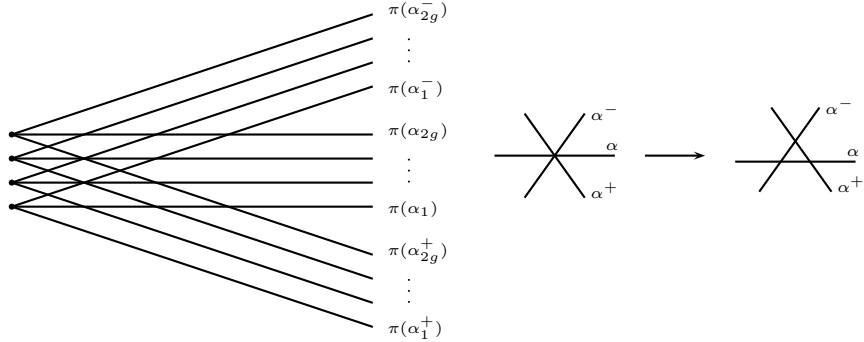
\begin{figure}[t]
\setlength{\unitlength}{8mm}
\begin{picture}(7,5.5)(0,-2.75)
\psset{unit=\unitlength}
\pscircle*(0,-0.6){0.05}
\pscircle*(0,-0.2){0.05}
\pscircle*(0,0.2){0.05}
\pscircle*(0,0.6){0.05}
\psline(0,-0.6)(6,-0.6)
\psline(0,-0.2)(6,-0.2)
\psline(0,0.2)(6,0.2)
\psline(0,0.6)(6,0.6)
\psline(0,-0.6)(6,1.4)
\psline(0,-0.2)(6,1.8)
\psline(0,0.2)(6,2.2)
\psline(0,0.6)(6,2.6)
\psline(0,-0.6)(6,-2.6)
\psline(0,-0.2)(6,-2.2)
\psline(0,0.2)(6,-1.8)
\psline(0,0.6)(6,-1.4)
\put(6.25,-0.7){\tiny $\pi(\alpha_1)$}
\put(6.25,0.6){\tiny $\pi(\alpha_{2g})$}
\psline[linestyle=dotted](6.6,0.2)(6.6,-0.2)
\put(6.25,-2.7){\tiny $\pi(\alpha^+_1)$}
\put(6.25,-1.4){\tiny $\pi(\alpha^+_{2g})$}
\psline[linestyle=dotted](6.6,-1.8)(6.6,-2.2)
\put(6.25,1.3){\tiny $\pi(\alpha^-_1)$}
\put(6.25,2.6){\tiny $\pi(\alpha^-_{2g})$}
\psline[linestyle=dotted](6.6,1.8)(6.6,2.2)
\end{picture}
\qquad
\begin{picture}(6,2)(-1,-3)
\psset{unit=\unitlength}
\psline(-1,0)(1,0)
\psline(-0.5,0.7)(0.5,-0.7)
\psline(-0.5,-0.7)(0.5,0.7)
\put(0.85,0.1){\tiny $\alpha$}
\put(0.6,0.6){\tiny $\alpha^-$}
\put(0.6,-0.7){\tiny $\alpha^+$}
\psline{->}(1.5,0)(2.5,0)
\psline(3,-0.1)(5,-0.1)
\psline(3.6,0.8)(4.6,-0.6)
\psline(3.4,-0.6)(4.4,0.8)
\put(4.85,0){\tiny $\alpha$}
\put(4.5,0.7){\tiny $\alpha^-$}
\put(4.7,-0.6){\tiny $\alpha^+$}
\end{picture}
\caption{The projection of the triple diagram
$(F',\alpha_i^-,\alpha_i,\alpha_i^+)$}\label{fig:triple}
\end{figure}

\proof 
As before, we use the fact that the thimbles $D_s^-$, $D_t$ and $D_u^+$ are
products of arcs in $F'$.
The image under $\pi$ of the triple diagram formed by these arcs
is depicted on Figure \ref{fig:triple} for convenience. This diagram has non-generic
triple intersections, which can be perturbed as in Figure \ref{fig:triple}
right.

Pick generators $z\in D_s^-\cap D_t$, $x\in D_t\cap D_u^+$, and $y\in
D_s^-\cap D_u^+$ (each viewed as $k$-tuples of intersections between arcs in
the diagram), and consider the homotopy class $\phi$ of a holomorphic
triangle contributing to the coefficient of $y$ in the product $z\cdot x$.
Projecting from the symmetric product to $F'$, we can think of
$\phi$ as a 2-chain in $F'$ with boundary on the arcs of the diagram,
staying within the bounded regions of the diagram.
Then the Maslov index $\mu(\phi)$ and the intersection number $i(\phi)$ of
$\phi$ with the diagonal divisor in $\Sym^k(F')$ are related to each
other by the following formula due to Sarkar \cite{Sarkar}:
\begin{equation}\label{eq:mutriangle}
\mu(\phi)=i(\phi)+2e(\phi)-k/2,
\end{equation}
where $e(\phi)$ is the {\it Euler measure} of the 2-chain $\phi$,
characterized by additivity and by the property that the Euler 
measure of an embedded $m$-gon with convex corners is $1-\frac{m}{4}$.
In our situation, we can draw the perturbed diagram in such a way that
all intersections occur at 60-degree and 120-degree angles as in Figure
\ref{fig:triple} right. The Euler measure of a convex polygonal region of the
diagram can then be computed by summing contributions from its vertices,
namely $+\frac{1}{12}$ for every vertex with a 60-degree angle, and
$-\frac{1}{12}$ for every vertex with a 120-degree angle; using additivity,
$e(\phi)$ can be expressed as a sum of local
contributions near the intersection points of the diagram covered by the
2-chain $\phi$.

View the 2-chain $\phi$ as the image of a holomorphic map $u$ from a Riemann
surface $S$ (with boundaries and strip-like ends) to $F'$ (as in
Lipshitz's approach to Heegaard-Floer theory), and fix an
intersection point $p$ in the triple diagram. If $u$ hits $p$ at an interior
point of $S$, then the local contributions to the multiplicities of $\phi$
in the four regions that meet at $p$ are all equal, hence the local
contribution to the Euler measure is zero. Likewise, if $u$ hits $p$ at a
point on the boundary of $S$, then (assuming $u$ is unbranched at $p$) 
locally the image of $u$ hits two of the
four regions that meet at $p$, one making a 60-degree angle and the other
making a 120-degree angle; in any case the local contributions to the Euler
measure cancel out. On the other hand, consider a strip-like end of
$S$ where $u$ converges to $p$ (i.e.\ an actual {\it corner} of the 
2-chain $\phi$), and recall the ordering condition on the boundaries of
$S$: going in the positive direction, the boundary of $S$ jumps from
some $\alpha_i^-$ to some $\alpha_i$, then to some
$\alpha_i^+$, then back to some $\alpha_i^-$ and so on. 
Looking at the local configurations of Figure \ref{fig:triple},
we see that locally $u$ maps into a region with a 60-degree angle 
at the vertex $p$ (unless there is a nearby boundary branch point, in which case $u$
locally maps into two 60-degree regions and one 120-degree region). Thus
each corner of $\phi$ contributes $+\frac{1}{12}$ to the Euler measure.
Summing over all $3k$ strip-like ends of $S$, we deduce that
$$e(\phi)=k/4.$$
The Floer product counts holomorphic discs such that $\mu(\phi)=0$; by
(\ref{eq:mutriangle}) these are precisely the discs for which $i(\phi)=0$,
i.e., using positivity of intersections, those which do not intersect 
the diagonal in $\Sym^k(F')$. Such holomorphic discs in 
$\Sym^k(F')$ can be viewed as $k$-tuples of holomorphic discs in
$F'$ (i.e., the domain $S$ is a disjoint union of $k$ discs), and
the Maslov index for such a product of discs is easily seen to be
the sum of the individual Maslov indices. Next, we recall that rigid
holomorphic discs on a Riemann surface are immersed polygonal regions
with convex corners; i.e., there are no branch points. (This conclusion
can also be reached by using equation (6) of \cite{LMW} which expresses the
Maslov index in terms of the Euler measure and the total number of branch points.)

Hence, the conclusion is the same as if our triple diagram
had been ``nice''
in the sense of \cite{LMW,Sarkar}: the Floer product counts $k$-tuples of 
immersed holomorphic triangles in $F'$ such that the 
corresponding map to $\Sym^k(F')$ does not hit the diagonal. 

Moreover, closer inspection of the triple diagram shows that immersed triangles are
actually embedded, and are contained either in a small neighborhood
$\mathcal{V}$ of
$V$ or in a small neighborhood $\mathcal{V}'$ of $V'$. (Recall that $V,V'$ are the two components
of \hbox{$\pi^{-1}(\{\Re z>0\})$}; in the limit where we consider the
unperturbed diagram of Figure~\ref{fig:triple} with triple intersections at the branch points $q_i$
the triangles cannot cross over the branch locus to jump from $V$ to $V'$, 
hence after perturbation they must lie in a small neighborhood of
either $V$ or $V'$.) 

Given a pair of triangles $T$ and $T'$ contained in $\mathcal{V}$, realized 
as the images of holomorphic maps $u,u'$ from the unit disc with three 
boundary marked points, the intersection number of the product map $(u,u')$ with the 
diagonal in $\Sym^2(F')$ can be evaluated by considering the rotation
number of the boundaries around each other: namely, embedding
$\mathcal{V}$ into $\R^2$, the restriction of $u'-u$ to the unit circle
defines a loop in $\R^2\setminus \{0\}$, whose degree is easily seen to
equal the intersection number of $(u,u')$ with the diagonal. One then checks
that configurations where $T$ and $T'$ are disjoint or intersect in a 
triangle (``head-to-tail overlap'') lead to an intersection number
of 0 and are hence allowed; however, all other configurations, e.g.\
when $T$ and $T'$ are contained inside each other or
intersect in a quadrilateral, lead to an intersection number of
1 and are hence forbidden. Similarly for triangles in $\mathcal{V}'$.

We conclude that the Floer product
counts $k$-tuples of embedded triangles in $F'$ which either are
disjoint or overlap head-to-tail (compare \cite[Lemma~2.6]{LMW}).

Recall that $\alpha_i^-$, $\alpha_j$ and $\alpha_l^+$ intersect pairwise
if and only if $i\le j\le l$. In that case, these curves bound exactly two 
embedded triangles $T_{ijl}$ and $T'_{ijl}$, the former contained in 
$\mathcal{V}$ and the latter contained in $\mathcal{V}'$, unless $i=j=l$
in which case there is a single triangle $T_{iii}=T'_{iii}$ 
obtained by deforming the triple
intersection at the branch point $p_i$ (see Figure \ref{fig:triple}).
Under the dictionary introduced in the proof of Proposition
\ref{prop:complex12}, the triangle $T_{ijl}$ corresponds to the
concatenation of strands connecting $a_i$ to $a_j$ and $a_j$ to $a_l$ to obtain
a strand connecting $a_i$ to $a_l$, while $T'_{ijl}$ corresponds to the
concatenation of strands connecting $a_{2g+i}$ to $a_{2g+j}$ and $a_{2g+l}$
to $a_{2g+l}$ to obtain a strand connecting $a_{2g+i}$ to $a_{2g+l}$;
the special case $i=j=l$ corresponds to the concatenation
of pairs of horizontal dotted lines.

Finally, consider two triangles $T_{ijl}$ and $T_{i'j'l'}$ where $i\le j\le
l$, $i'\le j'\le l'$, and $i<i'$: the concatenation of the strands 
connecting $a_i$ to $a_j$ and $a_j$ to $a_l$ intersects the concatenation
of the strands connecting $a_{i'}$ to $a_{j'}$ and $a_{j'}$ to $a_{l'}$
twice if and only if $j>j'$ and $l<l'$, i.e.\ the forbidden case is
$i<i'\le j'<j\le l<l'$. A tedious but straightforward enumeration of cases
shows that this is precisely the scenario in which the triangles $T_{ijl}$
and $T_{i'j'l'}$ overlap in a forbidden manner (other than head-to-tail).
Thus, the rules defining the product operations in $\A_{1/2}(F',k)$ and
$\mathcal{F}(f_{2g,k})$ agree with each other.
\endproof

The last ingredient is the following:

\begin{prop}\label{prop:higher12}
The higher compositions involving the thimbles $D_s$ $(s\in \mathcal{S}_k)$ in
$\F(f_{2g,k})$ are identically zero.
\end{prop}

\proof
The argument is similar to the first part of the proof of Proposition
\ref{prop:prod12}. Namely, the $\ell$-fold composition $m_\ell$ is determined by
picking $\ell+1$ different perturbations of the thimbles, and identifying 
them in the relevant portion of $\Sym^k(F')$ with products of arcs
obtained by perturbing the $\alpha_i$. The resulting diagram generalizes
in the obvious manner that of Figure \ref{fig:triple} (with $\ell+1$
sets of $2g$ arcs).

Consider the class $\phi$ of a holomorphic
$(\ell+1)$-pointed disc in $\Sym^k(F')$ that contributes to $m_\ell$:
then by Theorem 4.2 of \cite{Sarkar} we have
$$\mu(\phi)=i(\phi)+2e(\phi)-(\ell-1)k/2.$$
We can calculate the Euler measure as in the proof of
Proposition~\ref{prop:prod12} by setting up a perturbation of the diagram
in which all intersections occur at angles that are multiples of
$\pi/(\ell+1)$, and summing local contributions. (The local contribution
of a vertex with angle $r\pi$ to the Euler measure is $\frac14-\frac{r}{2}$).
The same argument as
before shows that each of the $(\ell+1)k$ corners contributes
$\frac14-\frac{1}{2(\ell+1)}=\frac{\ell-1}{4(\ell+1)}$ to the Euler measure, so that
$e(\phi)=(\ell-1)k/4$ and $\mu(\phi)=i(\phi)\ge 0$.

On the other hand, $m_\ell$ counts rigid holomorphic discs, i.e.\ discs of
Maslov index $2-\ell$. The above calculation shows that for $\ell\ge 3$ 
there are no such discs.
\endproof

Theorem \ref{thm:A12} follows from Propositions \ref{prop:complex12},
\ref{prop:prod12} and \ref{prop:higher12}.

\begin{rem}
Seidel's definition of the Fukaya category of a Lefschetz fibration
\cite{SeBook} is slightly more restrictive
than the version we gave in Section \ref{ss:fukayalf} above, in that
the only non-compact Lagrangians he allows are thimbles; the difference 
between the two versions is not expected to be significant when one passes 
to twisted complexes, but the cautious reader may wish to impose this
additional restriction. With this understood,
Theorem 18.24 of \cite{SeBook} implies that the Fukaya category of the
Lefschetz fibration $f_{2g,k}$ is generated by the exceptional collection
of thimbles $\{D_s,\ s\in \mathcal{S}_k\}$, in the sense that, after passing
to twisted complexes, the inclusion of the finite directed subcategory
$\A_{1/2}(F',k)$ into $\F(f_{2g,k})$ induces a quasi-equivalence $Tw
\A_{1/2}(F',k)\to Tw \F(f_{2g,k})$.
\end{rem}

\begin{rem}\label{rmk:A12inF}
In the next sections we will consider the slightly larger
surface $F$ and the Lefschetz fibration $f_{2g+1,k}:\Sym^k(F)\to\C$.
Assume that the points $p_j=i\theta_j$ have been chosen so that
$\theta_1<\dots<\theta_{2g}<0<\theta_{2g+1}$ and $|\theta_{2g+1}|\gg
|\theta_1|$: then the double covers $F\to\C$ and $F'\to\C$ can be
identified outside of a neighborhood of
the positive imaginary axis.
Passing to symmetric products, the
Lefschetz fibrations $f_{2g+1,k}$ and $f_{2g,k}$ agree over a large
convex open subset $\mathcal{U}$ which includes the $\binom{2g}{k}$ critical 
points of $f_{2g,k}$ and the corresponding thimbles. In this situation,
the Fukaya category $\F(f_{2g.k})$
embeds as a full $A_\infty$-subcategory 
of $\F(f_{2g+1,k})$, namely the subcategory generated by the
thimbles $D_s$, $s\in \mathcal{S}_k$ ($=\mathcal{S}^{2g}_k\subsetneq
\mathcal{S}^{2g+1}_k$). Indeed, the Lagrangian submanifolds and holomorphic
discs considered above all lie within $\mathcal{U}$ and do not see the
difference between $f_{2g,k}$ and $f_{2g+1,k}$. This alternative
description of $\A_{1/2}(F',k)$ as a subcategory of $\F(f_{2g+1,k})$ amounts
to viewing it as the strands algebra associated to a {\it
twice pointed} matched circle, rather than a pair of pointed circles.
\end{rem}

\section{Partially wrapped Fukaya categories and the algebra
$\A(F,k)$}\label{s:A}

\subsection{Partially wrapped Fukaya categories}

The Fukaya category of a Lefschetz fibration, as discussed in Section
\ref{ss:fukayalf}, is a particular instance of a more general construction,
which also encompasses the so-called wrapped Fukaya category 
(see \cite{AS}). In both cases, the idea is to allow noncompact
Lagrangian manifolds with appropriate behavior at infinity, and to
define their intersection theory by means of suitable Hamiltonian perturbations
which achieve a certain geometric behavior at infinity.

Let $(M,\omega)$ be an exact symplectic manifold with contact boundary.
Let $\hat{M}$ be the completion of $M$, i.e.\ the symplectic manifold
obtained by attaching 
to $M$ the positive part $([1,\infty)\times \partial M, d(r\alpha))$ of the 
symplectization of $\partial M$.
Let $H:\hat{M}\to \R$ be a Hamiltonian function such that $H\ge 0$ everywhere
and $H(r,y)=r$ on $[1,\infty)\times \partial M$.

The objects of the wrapped Fukaya category of $M$ (or $\hat{M}$) are exact 
Lagrangian submanifolds of $\hat{M}$ with cylindrical
ends modelled on Legendrian submanifolds of $\partial M$. The morphisms are defined
by $\hom(L_1,L_2)=\lim_{w\to +\infty} CF^*(\phi_{wH}(L_1),L_2)$, where 
$\phi_{wH}$ is the Hamiltonian diffeomorphism generated by $wH$; in the
symplectization, this Hamiltonian isotopy ``wraps'' $L_1$ by the time 
$w$ flow of the Reeb vector field.
The differential, composition, and higher products are defined in terms
of suitably perturbed versions of the holomorphic curve equation; i.e., they
can be understood in terms of holomorphic discs with boundary on 
increasingly perturbed versions of the Lagrangians. The reader is
referred to \S 3 of \cite{AS} for details.

We now consider ``partially wrapped'' Fukaya categories, tentatively
defined in the following manner:

\begin{qdefi}\label{def:fukrho}
Given a smooth function $\rho:\partial M\to [0,1]$, let
$H_\rho:\hat{M}\to\R$ be
a Hamiltonian function such that $H_\rho\ge 0$ everywhere and
$H_\rho(r,y)=\rho(y)\,r$ on the positive symplectization $[1,\infty)
\times \partial M$. The objects of the ``$\rho$-wrapped'' Fukaya
category $\F(M,\rho)$ are exact Lagrangian submanifolds of $\hat{M}$
with cylindrical ends modelled
on Legendrian submanifolds of $\partial M\setminus \rho^{-1}(0)$, and the
morphisms and compositions are defined by perturbing the Lagrangians 
by the long-time flow generated by $H_\rho$. Namely, $$\hom(L_1,L_2)=
\lim_{w\to+\infty} CF^*(\phi_{wH_\rho}(L_1),L_2),$$
and the differential, composition, and higher products are defined as
in \cite{AS} by
counting solutions of the Cauchy-Riemann equations perturbed by the
Hamiltonian flow of $H_\rho$.
\end{qdefi}

\noindent
At the boundary, the flow generated by $H_\rho$ 
can be viewed as the Reeb flow for the contact form $\rho^{-1}\alpha$ 
on the non-compact hypersurface 
$\{r=\rho^{-1}\}\simeq \partial M\setminus \rho^{-1}(0)$.
The effect of this modification is to slow down the wrapping
so that the long time flow never quite reaches $\rho^{-1}(0)$.

The direct limit in
Definition \ref{def:fukrho} relies on the existence of well-defined
continuation maps from $CF^*(\phi_{wH_\rho}(L_1),L_2)$ to
$CF^*(\phi_{w'H_\rho}(L_1),L_2)$ for $w'>w$. Even though exactness
prevents bubbling and the positivity of $H_\rho$ implies an a priori
energy bound on perturbed holomorphic discs, it is not entirely clear
that the construction is well-defined in full generality.%
\footnote{Ongoing work of Mohammed Abouzaid provides a treatment of
the important case where $\rho$ is lifted from an open book on $\partial M$.}
Here, we will only consider settings in
which $\phi_{wH_\rho}(L_1)$ and $L_2$ are transverse to each other
for all sufficiently large~$w$, and in particular no intersections
appear or disappear for $w\gg 0$. This simplifies things greatly,
as the complex stabilizes for large enough $w$. The continuation
maps can then be constructed by the ``homotopy method'' (see Appendix
\ref{appendix}), and turn out
to be the obvious ones for $w,w'$ large enough. The product maps can
also be defined similarly by counting ``cascades'' of (unperturbed)
holomorphic discs,
i.e.\ trees of rigid holomorphic discs with boundaries on the Lagrangian
submanifolds $\phi_{w_i H_\rho}(L_i)$ (where the parameter $w_i$
is sometimes fixed, and sometimes allowed to vary);
see Appendix \ref{appendix} for details.
However, in our case the upshot will be that the complexes, differentials
and products behave exactly as if one simply
considered sufficiently perturbed copies of the Lagrangians.

\begin{rem}\label{rmk:corner}
In many situations (exact Lefschetz fibrations over the disc with convex
fibers, symmetric products of Riemann surfaces with boundary, etc.), one
is naturally given an exact symplectic manifold with corners; one then needs
to ``round the corners'' to obtain a contact boundary. Concretely,
in the case of a product of Stein domains $M_1\times M_2$, we consider the
completed Stein manifolds $(\hat{M}_i,dd^c\varphi_i)$ and equip their product 
with the plurisubharmonic function $\pi_1^*\varphi_1+\pi_2^*\varphi_2$,
then restrict to a sublevel set to obtain a Stein domain again. More
importantly for our purposes, a similar procedure can be used to round
the corners of the symmetric product of a Riemann surface with boundary.
\end{rem}

The Fukaya category of a Lefschetz fibration over the disc
can now be understood as
a partially wrapped Fukaya category for a suitably chosen $\rho$, which
vanishes in the direction of the fiberwise boundary (recall that one
only considers Lagrangians on which the projection is proper) and also
in the fiber above one point of the boundary (or a subinterval of the
boundary).


\medskip

Another property that we expect of partially wrapped Fukaya categories
is the existence of ``acceleration'' $A_\infty$-functors $\F(M,\rho)\to
\F(M,\rho')$ whenever $\rho\le \rho'$ (i.e., from a ``less wrapped'' Fukaya
category to one that is ``more wrapped''). Specifically, because $H_\rho\le
H_{\rho'}$ one should have well-defined continuation maps from
$CF^*(\phi_{wH_\rho}(L_1),L_2)$ to $CF^*(\phi_{wH_{\rho'}}(L_1),L_2)$,
which (taking direct limits) define the linear term of the functor. 
However, the construction in the general case is well beyond the scope of this
paper. In our case, we will consider a very specific setting in which the
``less wrapped'' Floer complex turns out to be a {\em subcomplex}\/ of the 
``more wrapped'' one, and the acceleration functor is simply given by the
inclusion map.

\subsection{Partially wrapped categories for symmetric products}
\label{ss:Fz}

Let $S$ be a Riemann surface with boundary, equipped with an exact area form, 
and fix a point $z\in \partial S$.
Then $M=\Sym^k(S)$ is an exact symplectic
manifold with corners, and $V=\{z\}\times \Sym^{k-1}(S)\subset \partial M$.
As in Remark \ref{rmk:corner}, we can complete $M$ to
$\hat{M}=\Sym^k(\hat{S})$ where $\hat{S}$ is a punctured Riemann surface
obtained by attaching cylindrical ends to $S$, and use a plurisubharmonic
function on $\hat{M}$ to round the corners of $M$.

Consider a Lagrangian submanifold of $\hat{M}$ of the form
$\hat{L}=\hat\lambda_1\times \dots\times \hat\lambda_k$, where 
$\hat\lambda_i$ are disjoint properly embedded arcs in $\hat{S}$ obtained
by extending arcs $\lambda_i\subset S$ into the cylindrical ends. We assume
that the end points of $\lambda_i$ lie away from $z$, so that $\hat{L}$ is
tentatively an object of the partially wrapped Fukaya category.

Away from the diagonal strata, the exact symplectic structure on 
$\hat{M}$ is the product one, and the Hamiltonian $H$ that defines wrapped 
Floer homology in $\hat{M}$ is just a sum $H([z_1,\dots,z_k])=\sum_i
h(z_i)$, where $h$ is a Hamiltonian on $\hat{S}$. Thus, wrapping preserves
the product structure away from the diagonal: wrapping the product
Lagrangian $\hat{L}$ inside the
symmetric product $\hat{M}$ is equivalent to wrapping each factor
$\hat\lambda_i$ inside $\hat{S}$.

Due to the manner in which the smooth structure on the symmetric
product $\hat{M}=\Sym^k(\hat{S})$ is defined near the diagonal, it is
impossible for a nontrivial smooth Hamiltonian on 
$\hat{M}$ to preserve the product structure everywhere. Thus, if we wish
to preserve the interpretation of holomorphic discs in $\hat{M}$
in terms of holomorphic curves in $\hat{S}$, we cannot perturb the
holomorphic curve equations by an inhomogeneous Hamiltonian term. This
is one of the key reasons why we choose to set up wrapped Floer theory in the
language of cascades: then we consider genuine holomorphic discs (for the
product complex structure, or occasionally a small perturbation of the
product $J$ if needed to
achieve transversality) with boundary on product Lagrangian submanifolds
(recall that $H$ preserves the product structure
away from a small neighborhood of the diagonal, and in particular near
the Lagrangian submanifolds that we consider).

When we work relatively to $V=\{z\}\times \Sym^{k-1}(S)$, we are
``slowing down'' the wrapping whenever one of the $k$ components approaches
$z$ or, in the completion, the ray $\hat{Z}=\{z\}\times [1,\infty)$ generated 
by $z$ in the cylindrical end. Observe that $\{0\}\times \Sym^{k-1}(\C)\subset \Sym^k(\C)$ is the
(transverse) zero set of the $k$-th elementary symmetric function
$\sigma_k(x_1,\dots,x_k)=\prod x_i$. Hence, a natural way to associate
a partially wrapped Fukaya category to the pair $(M,V)$ is to
use a function $\rho$ which decomposes
as a product: $\rho([z_1,\dots,z_k])=\prod \rho_{S,z}(z_i)$, where
$\rho_{S,z}:S\to [0,1]$ is a smooth function that vanishes to order 2 at
$z\in\partial S$.

In this situation the wrapping flow no 
longer preserves the product structure as soon as one of the points
$z_i$ gets too close to $\hat{Z}$, even away from the diagonal.
So, if we consider two product Lagrangians $\hat{L}$,
$\hat{L}'$ which are disjoint from the support of $1-\rho$, the wrapping
perturbation applied to $\hat{L}$ only preserves the product structure
until $\phi_{wH_\rho}(\hat{L})$ enters the neighborhood of 
$\hat{Z}\times \Sym^{k-1}(\hat S)$ where $\rho\neq 1$. While it can be
checked that this is not an issue when it comes to the definition of the
Floer complexes and differentials, it is not entirely clear at this point that 
the product operations are well-defined and reduce to calculations
in the surface $S$. Thus, to avoid technical difficulties,
we will use a different choice of $\rho$ to
construct the $A_\infty$-category $\F_z$.
\medskip

Let us specialize right away to the case at hand, and consider again the
situation where $\hat{S}=\hat{F}$ is a punctured genus $g$ surface,
equipped with a double covering map
$\pi:\hat{F}\to\C$ with branch points $p_1,\dots,p_{2g+1}\in \C$ (with
$\mathrm{Im}\,p_1<\dots<\mathrm{Im}\,p_{2g+1}$), the subsurface 
$F\subset\hat{F}$
is the preimage of some large disc, say of radius $a$, and $z\in\partial F$ is one of
the two points in $\pi^{-1}(-a)$. 
\medskip

{\bf First version.}
We first equip $\hat{F}$ with a Hamiltonian constructed 
as follows. Let $a'>0$ be such that $\max |p_j|\ll a'\ll a$,
define $U=\pi^{-1}(D^2(a'))\subset \hat{F}$,
and let $\hat{Z}\subset \hat{F}$ be the component of 
$\pi^{-1}((-\infty,-a'])$ which passes through $z$.
We define $h_\rho(w)=\chi(w)\,|\pi(w)|^2$, where
$\chi:\hat{F}\to [0,1]$ is a smooth function which 
vanishes on $\hat{Z}\cup U$ and
equals 1 everywhere away from $\hat{Z}\cup U$.
Note that $h_\rho$
has no critical points outside of $\hat{Z}\cup U$, and it has the right
growth rate at infinity for the purposes of constructing a partially
wrapped Fukaya category for the pair $(F,\{z\})$.

The long-time flow of $X_{h_\rho}$ acts on 
properly embedded arcs in ${\hat{F}}$ in a straightforward manner:
the flow is identity inside the subset $U$, while in the cylindrical end
the flow wraps in the positive direction and
accumulates onto the ray $\hat{Z}$ (if $\chi$ is chosen suitably). 
To be more specific,
we identify $\hat{F}\setminus U$ with a cylinder, with radial
coordinate $|\pi(\cdot)|^2$ and angular coordinate $\vartheta=\frac12
\arg\pi(\cdot)$ (with, say, $\vartheta=\pi/2$ at $z$ to fix things).
The level sets of $h_\rho$ are asymptotic (from both sides) to the ray 
$\smash{\hat{Z}}$, where $\vartheta=\pi/2$; thus the wrapping by the
positive (resp.\ negative) time flow generated by $h_\rho$ moves any point
outside $\hat{Z}\cup U$ towards infinity, with $\vartheta$
increasing (resp.\ decreasing) towards~$\pi/2$.

\begin{figure}[t]
\setlength{\unitlength}{1cm}
\begin{picture}(5.5,3)(-0.5,-1.5)
\psset{unit=\unitlength}
\psellipticarc[linewidth=0.5pt,linestyle=dashed,dash=2pt 2pt](4.5,0)(0.2,1.5){90}{-90}
\psellipticarc(4.5,0)(0.2,1.5){-90}{90}
\psline[linearc=1.5](4.5,1.5)(-1,1.5)(-1,-1.5)(4.5,-1.5)
\pscircle*(4.58,1.3){0.07} \put(4.7,1.3){\tiny $z$}
\psellipse(0.35,0)(0.35,0.2)
\psellipse(1.65,0)(0.35,0.2)
\psellipticarc(3,0)(1.03,0.21){90}{180}
\psellipticarc(3,0)(1.71,0.51){90}{180}
\psellipticarc(3,0)(2.3,0.81){90}{180}
\psellipticarc(3,0)(3,1.11){90}{180}
\psline(3,0.2)(4.68,0.2)
\psline(3,0.5)(4.65,0.5)
\psline(3,0.8)(4.65,0.8)
\psline(3,1.1)(4.6,1.1)
\psellipticarc[linestyle=dashed,dash=2pt 2pt](3,0)(1.03,0.215){180}{270}
\psellipticarc[linestyle=dashed,dash=2pt 2pt](3,0)(1.71,0.515){180}{270}
\psellipticarc[linestyle=dashed,dash=2pt 2pt](3,0)(2.3,0.815){180}{270}
\psellipticarc[linestyle=dashed,dash=2pt 2pt](3,0)(3,1.115){180}{270}
\psline[linestyle=dashed,dash=2pt 2pt](3.05,-0.2)(4.32,-0.2)
\psline[linestyle=dashed,dash=2pt 2pt](3.05,-0.5)(4.32,-0.5)
\psline[linestyle=dashed,dash=2pt 2pt](3.05,-0.8)(4.35,-0.8)
\psline[linestyle=dashed,dash=2pt 2pt](3.05,-1.1)(4.38,-1.1)
\put(4.79,0.15){\tiny $\alpha_1$}
\put(4.75,1.05){\tiny $\alpha_{2g}$}
\psline[linestyle=dotted](4.9,0.85)(4.9,0.45)
\end{picture}
\qquad\quad
\begin{picture}(5.5,3)(-0.5,-1.5)
\psset{unit=\unitlength}
\psellipticarc[linewidth=0.5pt,linestyle=dashed,dash=2pt 2pt](4.5,0)(0.2,1.5){90}{-90}
\psellipticarc(4.5,0)(0.2,1.5){-90}{90}
\psline[linearc=1.5](4.5,1.5)(-1,1.5)(-1,-1.5)(4.5,-1.5)
\pscircle*(4.58,1.35){0.05} \put(4.7,1.3){\tiny $z$}
\psset{linecolor=red}
\psellipticarc[linestyle=dashed,dash=2pt 2pt](2.2,0)(0.23,0.215){180}{270}
\psellipticarc[linestyle=dashed,dash=2pt 2pt](2.2,0)(0.91,0.365){180}{270}
\psellipticarc[linestyle=dashed,dash=2pt 2pt](2.2,0)(1.5,0.515){180}{270}
\psellipticarc[linestyle=dashed,dash=2pt 2pt](2.2,0)(2.2,0.665){180}{270}
\pscurve[linestyle=dashed,dash=2pt 2pt](2.2,-0.2)(2.3,-0.15)(3.3,1.02)(4.35,1.17)
\pscurve[linestyle=dashed,dash=2pt 2pt](2.2,-0.35)(2.37,-0.3)(3.38,0.87)(4.35,1.02)
\pscurve[linestyle=dashed,dash=2pt 2pt](2.2,-0.5)(2.45,-0.45)(3.45,0.72)(4.3,0.87)
\pscurve[linestyle=dashed,dash=2pt 2pt](2.2,-0.65)(2.52,-0.6)(3.53,0.57)(4.3,0.72)
\pscurve[linestyle=dashed,dash=2pt 2pt](3.4,-1.5)(3.5,-1.45)(4.1,0.5)(4.3,0.57)
\pscurve[linestyle=dashed,dash=2pt 2pt](3.6,-1.5)(3.7,-1.45)(4.15,0.35)(4.3,0.42)
\pscurve[linestyle=dashed,dash=2pt 2pt](3.8,-1.5)(3.88,-1.45)(4.2,0.2)(4.3,0.27)
\pscurve[linestyle=dashed,dash=2pt 2pt](4,-1.5)(4.05,-1.45)(4.2,-0.1)(4.3,0.12)
\psset{linecolor=blue}
\psellipticarc(3,0)(1.03,0.81){90}{180}
\psellipticarc(3,0)(1.71,0.96){90}{180}
\psellipticarc(3,0)(2.3,1.11){90}{180}
\psellipticarc(3,0)(3,1.26){90}{180}
\psline(3,0.8)(4.65,0.8)
\psline(3,0.95)(4.63,0.95)
\psline(3,1.1)(4.6,1.1)
\psline(3,1.25)(4.58,1.25)
\psellipticarc[linestyle=dashed,dash=2pt 2pt](1.8,0)(1.8,1.5){180}{270}
\pscurve(1.8,-1.5)(2.3,-1.4)(3.6,0.5)(4.65,0.65)
\psellipticarc[linestyle=dashed,dash=2pt 2pt](2.1,0)(1.4,1.5){180}{270}
\pscurve(2.1,-1.5)(2.6,-1.4)(3.7,0.35)(4.68,0.5)
\psellipticarc[linestyle=dashed,dash=2pt 2pt](2.4,0)(1.1,1.5){180}{270}
\pscurve(2.4,-1.5)(2.9,-1.4)(3.8,0.2)(4.68,0.35)
\psellipticarc[linestyle=dashed,dash=2pt 2pt](2.7,0)(0.72,1.5){180}{270}
\pscurve(2.7,-1.5)(3.1,-1.4)(3.87,0.05)(4.7,0.2)
\psset{linecolor=red}
\psellipticarc(2.2,0)(0.23,0.215){90}{180}
\psellipticarc(2,0)(0.71,0.365){90}{180}
\psellipticarc(2.2,0)(1.5,0.515){90}{180}
\psellipticarc(2.2,0)(2.2,0.665){90}{180}
\pscurve(2.2,0.2)(2.3,0.18)(3.2,-1.45)(3.4,-1.5)
\pscurve(2,0.35)(2.4,0.3)(3.4,-1.45)(3.6,-1.5)
\pscurve(2.2,0.5)(2.6,0.42)(3.6,-1.45)(3.8,-1.5)
\pscurve(2.2,0.65)(2.7,0.55)(3.8,-1.45)(4,-1.5)
\psset{linecolor=black}
\psellipse(0.35,0)(0.35,0.2)
\psellipse(1.65,0)(0.35,0.2)
\pscircle*(0,0){0.04} \pscircle*(0.7,0){0.04} \pscircle*(1.3,0){0.04} \pscircle*(2,0){0.04}
\pscircle*(1.01,0.55){0.04} \pscircle*(1.75,0.63){0.04} \pscircle*(2.4,0.63){0.04} 
\pscircle*(1.5,0.44){0.04} \pscircle*(2.2,0.50){0.04} \pscircle*(2.09,0.35){0.04} 
\pscircle*(0.79,-0.5){0.04} \pscircle*(1.4,-0.61){0.04} \pscircle*(2.07,-0.65){0.04} 
\pscircle*(1.36,-0.41){0.04} \pscircle*(2.03,-0.5){0.04} \pscircle*(2.01,-0.35){0.04} 
\pscircle*(3.13,-0.07){0.04} \pscircle*(3.22,-0.33){0.04} 
\pscircle*(3.32,-0.66){0.04} \pscircle*(3.41,-0.89){0.04} 
\pscircle*(3.03,-0.26){0.04} \pscircle*(3.12,-0.59){0.04} 
\pscircle*(3.24,-0.92){0.04} \pscircle*(3.32,-1.11){0.04} 
\pscircle*(2.91,-0.56){0.04} \pscircle*(3.01,-0.87){0.04} 
\pscircle*(3.13,-1.13){0.04} \pscircle*(3.23,-1.27){0.04} 
\pscircle*(2.79,-0.78){0.04} \pscircle*(2.91,-1.07){0.04} 
\pscircle*(3.04,-1.28){0.04} \pscircle*(3.13,-1.38){0.04} 
\put(1.4,-1.8){\tiny $\tilde\alpha^-_{2g}\,\cdots\,\tilde\alpha^-_1$}
\put(3,-1.8){\tiny $\tilde\alpha^+_1\cdots\tilde\alpha^+_{2g}$}
\end{picture}

\caption{The arcs $\alpha_j$ and $\tilde\alpha_j^\pm$ on $\hat{F}$}
\label{fig:sigma'}
\end{figure}
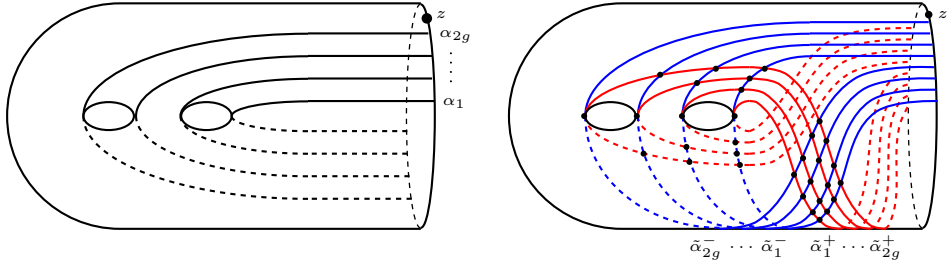

In particular, the positive (resp.\ negative) time flow of $h_\rho$ maps the arcs
$\alpha_j=\pi^{-1}(p_j+\R_{\ge 0})\subset \hat{F}$ to arcs which,
after a compactly supported isotopy, look like
the arcs $\tilde\alpha_j^-$ (resp.\ $\tilde\alpha_j^+$) pictured
in Figure \ref{fig:sigma'} below (the last arc
$\alpha_{2g+1}$ is not pictured, but behaves in a similar manner). Note however that, due to
the degeneracy of $h_\rho$ inside $U$, the arcs $\phi_{\pm
wh_\rho}(\alpha_i)$ are never transverse to each other: without further perturbation,
the flow of $h_\rho$ only yields an $A_\infty$-precategory, i.e.\ morphisms
and compositions are only defined for objects which are mutually transverse
within $U$ (and in particular, endomorphisms are not well-defined). In order
to construct an honest $A_\infty$-category one needs to choose further 
(compactly supported) Hamiltonian perturbations in a consistent manner; 
see below.

With $h_\rho$ at hand, we equip $\hat{M}=\Sym^k(\hat{F})$ with a Hamiltonian
$H_\rho$ such that, outside of a small neighborhood of the diagonal strata,
$H_\rho([z_1,\dots,z_k])=\sum_i h_\rho(z_i)$. In particular, the Hamiltonian
flow generated by $H_\rho$ preserves the product structure away from a small
neighborhood of the diagonal.
Thus, given $k$ disjoint embedded arcs
$\hat{\lambda}_1,\dots,\hat{\lambda}_k\subset \hat{F}$, for suitable
values of $w$ the flow maps
$\hat{L}=\hat{\lambda}_1\times\dots\times\hat{\lambda}_k$ to
\begin{equation}\label{eq:prodflow}
\phi_{wH_\rho}(\hat{L})=\phi_{wh_\rho}(\hat{\lambda}_1)\times\dots\times
\phi_{wh_\rho}(\hat{\lambda}_k).
\end{equation}

\begin{rem}\label{rmk:diagonal}
Due to the specifics of the construction, for large $w$ the image under
$\phi_{wH_\rho}$ of
a product of disjoint arcs does approach the diagonal, where the product structure
is not preserved by the flow; we will want to correct this and
ensure that (\ref{eq:prodflow}) holds for all $w$.
There are several ways to proceed. A first option
would be to modify the definition of $h_\rho$ appropriately in order to
control the manner in which things can accumulate towards the diagonal;
this comes at the expense of
making $h_\rho$ non-constant over $U$, which complicates the geometric behavior
of the flow. A second possibility, suggested by
the referee, is to let the Hamiltonian $H_\rho$ be singular along the diagonal. This
is valid because in our technical setup the Hamiltonian is never used 
to perturb the Cauchy-Riemann equation (see Appendix \ref{appendix}); instead, we consider honest
holomorphic curves with boundary on the images of the Lagrangians under
the flow, and these remain smooth for Lagrangians which do not
intersect the diagonal. One would also need to make the K\"ahler
form singular along the diagonal, which is actually not a problem in our
case. A third approach, strictly equivalent to the previous one and
which we will use instead, is to allow the
choice of $H_\rho$ near the diagonal to depend on the product Lagrangian $\hat{L}$ 
under consideration; it is then not hard to ensure that (\ref{eq:prodflow})
holds for all $w$.
\end{rem}

{\bf Hamiltonian perturbations.}
One way to address the degeneracy of $h_\rho$ would be to replace it by a
non-degenerate Hamiltonian; however, this affects the long-term
dynamics inside $U$ in a counter-intuitive manner. Another approach is to
keep using a degenerate Hamiltonian, but further add 
small compactly supported Hamiltonian perturbations in order to achieve
transversality. This is conceptually similar to the approach taken by Seidel
in \cite{SeBook}, except we again consider cascades of honest holomorphic
curves with boundaries on perturbed Lagrangian submanifolds, rather than
perturbing the holomorphic curve equation.


Concretely, for each pair of Lagrangians $(L_1,L_2)$, we choose a family
of Hamiltonians $\{H'_{L_1,L_2,\tau}\}_{\tau\ge 0}$, uniformly bounded,
depending smoothly on $\tau$, and with $H'_{L_1,L_2,0}=0$, with
the property that $\phi_{wH_\rho+H'_{L_1,L_2,w}}(L_1)$ is transverse
to $L_2$ for all sufficiently large $w$.
We then define $$\hom(L_1,L_2)=\lim_{w\to +\infty}
CF^*(\phi_{wH_\rho+H'_{L_1,L_2,w}}(L_1),L_2).$$
The definition of product structures requires additional transversality
properties, and the choice of suitable
homotopies between the Hamiltonian perturbations; these are incorporated
into the definition of the $A_\infty$-operations via cascades. The details
can be found in \S \ref{ss:cascadeham} where, for
simplicity, we only describe the construction in the case where
the perturbation $H'_{L_1,L_2,w}=H'_{L_1,w}$ is chosen to depend only
on $L_1$ and $w$, not on $L_2$.
This assumption makes the construction much simpler, but prevents us
from achieving transversality for arbitrary pairs of Lagrangians.

In our case, we will essentially be able to use small multiples of the same 
Hamiltonian perturbation $H'$ for all the thimbles $D_s$. Namely,
we pick a Hamiltonian $h':\hat{F}\to\R$ with the following properties:
\begin{itemize}
\item the branch points $q_1,\dots,q_{2g+1}$ of the projection $\pi$
are nondegenerate critical points of $h'$;
\item $h'$ is bounded, and constant on the level sets of $h_\rho$ in the cylindrical
end of $\hat{F}$;
\item $h'_{|\alpha_j}$ is a Morse function
with a single minimum at $q_j$.
\end{itemize}
The second property ensures that the flow of $h'$ commutes with that of 
$h_\rho$ (which makes perturbed cascades more intuitive) and does not affect
the behavior at infinity; the third
one ensures that the images of the arcs $\alpha_j$ under the flow generated
by $wh_\rho+\epsilon h'$ (for $\epsilon>0$) behave exactly like the arcs
$\tilde{\alpha}_j^-$ pictured in Figure
\ref{fig:sigma'}.

As above, we define a Hamiltonian $H'$ on $\hat{M}=\Sym^k(\hat{F})$ such
that, outside of a small neighborhood of the diagonal strata,
$H'([z_1,\dots,z_k])=\sum_i h'(z_i)$. We can in particular arrange for its
Hamiltonian flow to preserve the product structure away from the diagonal
and commute with that of $H_\rho$.
Thus, given $k$ sufficiently disjoint embedded arcs
$\hat{\lambda}_1,\dots,\hat{\lambda}_k\subset \hat{F}$, the flow generated
by $wH_\rho+\epsilon H'$ maps
the product $\hat{L}=\hat{\lambda}_1\times\dots\times\hat{\lambda}_k$ to
$\phi_{wH_\rho+\epsilon H'}(\hat{L})=\phi_{wh_\rho+\epsilon h'}(\hat{\lambda}_1)\times\dots\times
\phi_{wh_\rho+\epsilon h'}(\hat{\lambda}_k)$, at least away from the
diagonal. As explained in Remark \ref{rmk:diagonal}, we can ensure that this
identity remains true for all large $w$ and small $\epsilon$ by letting the choices of $H_\rho$ and
$H'$ near the diagonal depend on the Lagrangian $\hat{L}$; we denote these
choices by $H_{\rho,\hat{L}}$ and $\smash{H'_{\hat{L}}}$, though we will often drop
the subscript from the notation. (Here again, another option would have been to let
$H'$ be singular along the diagonal).

For $s\in \mathcal{S}^{2g+1}_k$ and $\tau\ge 0$, we set
$H'_{D_s,\tau}=\epsilon(\tau)H'_{D_s}$, 
where $\epsilon$ is a monotonically
increasing smooth function with $\epsilon(0)=0$ and bounded by a small
positive constant.
By construction, the image of $D_s$ under the flow generated by
$wH_{\rho,D_s}+H'_{D_s,w}$ is transverse to $D_t$ for all large enough $w$,
without any intersections being created or cancelled;
moreover, the construction of $H'$ is flexible enough to ensure
that the appropriate moduli spaces of holomorphic discs are generically
regular (see below).
Thus, the necessary technical conditions (Definition \ref{ass:ass}, 
as modified in \S \ref{ss:cascadeham} to include the
perturbations) are satisfied.

\begin{defn}\label{def:Fz} 
We denote by $\F_z$ the $A_\infty$-(pre)category whose objects 
are \begin{enumerate}
\item closed exact Lagrangian submanifolds contained in $\Sym^k(U)\subset\Sym^k(\hat{F})$, and
\item exact Lagrangian submanifolds of the form $\hat\lambda_1\times
\dots\times \hat\lambda_k$, where the $\hat\lambda_i$ are disjoint
properly embedded arcs in $\hat{F}$ such that
$\hat\lambda_i\cap (\hat{F}\setminus U)$ consists of two components which
project via $\pi$ to straight lines contained in the right half-plane
$\mathrm{Re}\,\pi>0$,
\end{enumerate}
with morphisms and compositions defined by partially wrapped Floer
theory (in the sense of Appendix \ref{appendix}) with respect to a small
perturbation of the
product complex structure $J$, the Hamiltonian $H_\rho$, and suitably
chosen small bounded Hamiltonian perturbations.
\end{defn}

\begin{rem}
As pointed out by the referee, in general it is not known whether
transversality can be achieved for product Lagrangian submanifolds using
the product complex structure; thus, we  need to work with a small
generic perturbation of the product $J$ within a suitable class of
almost-complex structures (as in Heegaard-Floer theory). However,
where the thimbles $D_s$ are concerned (e.g.\ in the proof of Theorem
\ref{thm:A}) this makes no difference:
indeed, for these specific objects one easily checks that the moduli 
spaces are transversely cut out for the product complex structure, and
hence the moduli spaces remain exactly the same after a slight perturbation
of the complex structure.
\end{rem}

\noindent
We leave the Hamiltonian perturbations $H'_{L,\tau}$ unspecified except for the
thimbles $D_s$. Indeed, the actual choice is
immaterial, and the Fukaya categories constructed for different
choices of perturbations are quasi-equivalent (the argument is essentially
the same as in \cite{SeBook}). The only key requirement is that we need the
perturbations to be small and bounded so as to not significantly affect 
the behavior at infinity of the long-time flow (for
non-compact objects as in Definition \ref{def:Fz}(2), the properness of
$h_\rho$ away from the ray $\hat{Z}$ ensures that a small bounded Hamiltonian
perturbation pulled back from $\hat{F}$ does not modify the large-scale behavior). 

We also note that, since the compact objects in Definition \ref{def:Fz}(1)
are required to lie in $\Sym^k(U)$, over which $H_\rho$ vanishes, they
are not affected by the wrapping.

In general, due to our simplifying assumption on the Hamiltonian
perturbations we cannot expect transversality in the sense of 
\S \ref{ss:cascadeham} to hold for arbitrary Lagrangian submanifolds,
so that $\F_z$ is only an $A_\infty$-precategory, i.e.\ morphisms
and compositions are only defined for objects which satisfy the 
transversality conditions. The issue is fairly mild, and can be
ignored for all practical purposes, since any ordered sequence
of thimbles $D_s$ is transverse. Nonetheless, the cautious reader may
wish to restrict the set of objects of $\F_z$ to some fixed countable 
collection of Lagrangians (such that
every isotopy class is represented, and including the thimbles $D_s$)
for which transversality can be achieved.

\begin{rem}\label{rmk:lesswrapped}
If we modify the construction of $h_\rho$ to make the cut-off function $\chi$ vanish
on {\it both} components of $\pi^{-1}((-\infty,-a'])$, then we obtain a ``less
wrapped'' category which is fairly closely related to the Fukaya category
of the Lefschetz fibration $f_{2g+1,k}$, at least as far as the thimbles
$D_s=\prod_{i\in s} \alpha_i$ are concerned. Indeed, the flow still
preserves the product structure, but since the Hamiltonian now vanishes
over the entire preimage of an arc connecting $p_{2g+1}$ to $-\infty$,
the wrapping now accumulates on the {\it two} 
infinite rays $\vartheta=\pm \pi/2$ in the cylindrical
end and never crosses the preimage of the negative real axis. Thus
the flow now maps the arcs $\alpha_i$ to a configuration which, for all
practical purposes, behaves interchangeably with the arcs $\alpha_i^-$
previously introduced. It is an exercise left to the reader to adapt
the argument below and show that, 
in this ``less wrapped'' Fukaya category, the
$A_\infty$-algebra associated to the thimbles $D_s$, $s\in
\mathcal{S}^{2g}_k$ is again $\A_{1/2}(F',k)$, just as in 
$\F(f_{2g+1,k})$ (cf.\ Remark \ref{rmk:A12inF}).
\end{rem}

In the rest of this section, we will be considering the thimbles
$D_s=\prod_{j\in s}\alpha_j$, where 
$s\in \mathcal{S}^{2g}_k$ ranges over all $k$-element subsets 
of $\{1,\dots,2g\}$, viewed as
objects of the partially wrapped Fukaya category $\F_z$. The following
lemma says that we can ignore the technicalities of the construction
of the partially wrapped Fukaya category, and simply perturb
$D_s$ to $\tilde{D}_s^\pm = \prod_{j\in s} \tilde\alpha_j^\pm $, 
where the $\tilde\alpha_j^\pm$ are the arcs pictured in
Figure~\ref{fig:sigma'}.

\begin{lem}\label{l:simplify}
The full subcategory of $\F_z$ with objects $D_s,\ s\in \mathcal{S}^{2g}_k$
is quasi-isomorphic to the $A_\infty$-category with the same objects,
$\hom(D_s,D_t)=CF^*(\tilde{D}_s^-,\tilde{D}_t^+)$, and product operations
given by counting holomorphic discs bounded by suitably perturbed versions
of the $D_s$ (using the long-time flow of $H_\rho$ and the Hamiltonian
perturbation $H'$).
\end{lem}

\proof Lemma \ref{l:trivcontham} gives
a criterion under which the infinitely generated complex used to define
$\hom(D_s,D_t)$ in the partially wrapped Fukaya category $\F_z$ can be 
replaced by the ordinary Floer complex $CF^*(\phi_{wH_\rho+\epsilon(w)H'}(D_s),D_t)$
(which is naturally isomorphic to $CF^*(\tilde{D}_s^-,\tilde{D}_t^+)$),
and the cascades used to define $A_\infty$-operations are
simply rigid holomorphic discs with boundaries on the images of the 
given Lagrangians under $\phi_{\tau H_\rho+\epsilon(\tau)H'}$ (for
sufficiently different values of $\tau$). 

The first 
assumption of the lemma, i.e.\ the transversality of
$\phi_{(\tau+w)H_\rho+\epsilon(\tau+w)H'}(D_s)$ to $\phi_{\tau H_\rho+
\epsilon(\tau)H'}(D_t)$ for all $s,t\in \mathcal{S}^{2g}_k$, 
$\tau\ge 0$ and large enough $w$, follows from the
construction of $H'$ (using the fact that
the function $\epsilon$ is monotonically increasing).
Thus we only need to check that, for $s_0,\dots,s_\ell\in
\mathcal{S}^{2g}_k$ and $\tau_0\gg \tau_1\gg \dots\gg \tau_\ell\ge 0$,
the Lagrangian submanifolds 
$\phi_{\tau_i H_\rho+\epsilon(\tau_i)H'}(D_{s_i})$ 
never bound any holomorphic discs of Maslov
index less than $2-\ell$. 

We claim that the diagram formed by the arcs $\alpha_j$ and their images under
the flow generated by $\tau_i h_\rho+\epsilon(\tau_i)h'$ has the same nice properties as the diagram
considered in Section \ref{s:A12}. Namely, one can draw the bounded 
regions of the diagram formed by $\ell+1$ different increasingly wrapped
perturbations of the arcs $\alpha_j$ ($1\le j\le 2g$) in such a way that
all intersections occur at angles that are multiples of $\pi/(\ell+1)$,
and find as in the proof of Proposition \ref{prop:higher12} that the
Maslov index of any holomorphic disc is equal to its intersection number
with the diagonal strata, $\mu(\phi)=i(\phi)\ge 0$. (See the argument
below and Figures \ref{fig:rectangles} and \ref{fig:triple'}.)
This immediately implies the absence of discs of index less than
$2-\ell$ except in the case $\ell=1$.

Next, we observe that a Maslov index 0 holomorphic strip would have to be disjoint
from the diagonal strata in $\Sym^k(\hat{F})$ (since $\mu(\phi)=i(\phi)=0$).
Thus, such a strip can be viewed as a $k$-tuple of holomorphic strips in
$\hat{F}$; however, $\phi_{wh_\rho+\epsilon(w)h'}(\alpha_i)$ and $\alpha_j$ (or
equivalently, $\tilde\alpha_i^-$ and $\tilde\alpha_j^+$)
do not bound any non-trivial discs in $\hat{F}$.
Hence there are no nonconstant Maslov index 0 holomorphic strips, which 
completes the verification of the assumptions of Lemma \ref{l:trivcontham}.
The result follows.
\endproof

\subsection{Proof of Theorem \ref{thm:A}}\label{ss:pfthmA}
The proof of Theorem \ref{thm:A} goes along the same lines as 
that of Theorem \ref{thm:A12}, but using the arcs $\tilde\alpha_j^\pm$
instead of $\alpha_j^\pm$. The theorem follows from Lemma \ref{l:simplify}
and the following three propositions.

\begin{prop}\label{prop:complex}
The chain complexes $\hom_{\F_z}(D_s,D_t)$
and $\hom_{\A(F,k)}(s,t)$ are isomorphic for all $s,t\in \mathcal{S}^{2g}_k$.
\end{prop}

\proof The intersections of $\tilde{D}_s^-$ with $\tilde{D}_t^+$ consist of
$k$-tuples of intersections between the arcs $\tilde\alpha_i^-$, $i\in s$
and $\tilde\alpha_j^+$, $j\in t$. These can be determined by looking at
Figure~\ref{fig:sigma'}. Namely, the ``left half'' of $\hat{F}$ looks
similar to the configuration of Section \ref{s:A12}, while in the ``right
half'' the wrapping creates one new intersection between each
$\tilde\alpha_i^-$ and each $\tilde\alpha_j^+$. Intersections of the
first type are again interpreted as strands which do not cross the
interval $[2g,2g+1]$ on the pointed matched circle, while the new
intersection point between $\tilde\alpha_i^-$ and $\tilde\alpha_j^+$ is
interpreted as a strand connecting $a_i$ to $a_{2g+j}$. The dictionary
between intersection points and strands is now as follows:
\begin{itemize}
\item For $i<j$, $\tilde\alpha^-_i\cap \tilde\alpha^+_j$ consists of three
points; the point at the upper-left on the front part of Figure
\ref{fig:sigma'} is interpreted as 
$\left[\begin{smallmatrix}2g+i\\2g+j\end{smallmatrix}\right]$, while the
point at the lower-left on the back part of the figure corresponds to
$\left[\begin{smallmatrix}i\\j\end{smallmatrix}\right]$, and the
point in the lower-right part of the figure corresponds to
$\left[\begin{smallmatrix}i\\2g+j\end{smallmatrix}\right]$;
\item For $i=j$, $\tilde\alpha^-_i\cap \tilde\alpha^+_j$ consists of two
points; the branch point of $\pi$ in the left part of the figure corresponds
to the double dotted line $\left[\begin{smallmatrix}i\\\vphantom{i}\end{smallmatrix}\right]$, while 
the point in the lower-right part of the figure corresponds to
$\left[\begin{smallmatrix}i\\2g+i\end{smallmatrix}\right]$;
\item For $i>j$, $\tilde\alpha^-_i\cap \tilde\alpha^+_j$ consists of a
single point, interpreted as
$\left[\begin{smallmatrix}i\\2g+j\end{smallmatrix}\right]$.
\end{itemize}
As before, by considering the set of $k$-tuples for which the labels in $s$
and $t$ each appear exactly once we obtain a bijection between the generators 
of $\hom_{\F_z}(D_s,D_t)$ and $\hom_{\A(F,k)}(s,t)$.

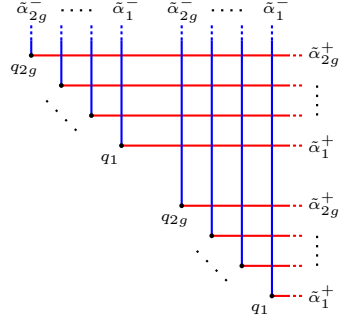
\begin{figure}[t]
\setlength{\unitlength}{4mm}
\begin{picture}(10,10)(-1.2,-1.2)
\psset{unit=\unitlength}
\psset{linecolor=red}
\psline(-1,7)(7.5,7) \psline(0,6)(7.5,6) \psline(1,5)(7.5,5) \psline(2,4)(7.5,4)   
\psline(4,2)(7.5,2) \psline(5,1)(7.5,1) \psline(6,0)(7.5,0) \psline(7,-1)(7.5,-1)
\multips(7.5,-1)(0,1){4}{\psline[linestyle=dashed, dash=1pt 1pt](0.5,0)(0,0)}
\multips(7.5,4)(0,1){4}{\psline[linestyle=dashed, dash=1pt 1pt](0.5,0)(0,0)}
\psset{linecolor=blue}
\psline(7,-1)(7,7.5) \psline(6,0)(6,7.5) \psline(5,1)(5,7.5) \psline(4,2)(4,7.5)   
\psline(2,4)(2,7.5) \psline(1,5)(1,7.5) \psline(0,6)(0,7.5) \psline(-1,7)(-1,7.5)
\multips(-1,7.5)(1,0){4}{\psline[linestyle=dashed, dash=1pt 1pt](0,0.5)(0,0)}
\multips(4,7.5)(1,0){4}{\psline[linestyle=dashed, dash=1pt 1pt](0,0.5)(0,0)}
\psset{linecolor=black}
\multips(-1,7)(1,-1){4}{\pscircle*(0,0){0.08}}
\multips(4,2)(1,-1){4}{\pscircle*(0,0){0.08}}
\put(-1.7,6.5){\tiny $q_{2g}$}
\put(1.3,3.5){\tiny $q_1$}
\psline[linestyle=dashed, dash=1pt 5pt](-0.5,5.5)(0.5,4.5)
\put(3.3,1.5){\tiny $q_{2g}$}
\put(6.3,-1.5){\tiny $q_1$}
\psline[linestyle=dashed, dash=1pt 5pt](4.5,0.5)(5.5,-0.5)
\put(-1.5,8.4){\tiny $\tilde\alpha^-_{2g}$}
\put(1.7,8.4){\tiny $\tilde\alpha^-_{1}$}
\psline[linestyle=dashed, dash=1pt 3pt](0,8.5)(1,8.5)
\put(3.5,8.4){\tiny $\tilde\alpha^-_{2g}$}
\put(6.7,8.4){\tiny $\tilde\alpha^-_{1}$}
\psline[linestyle=dashed, dash=1pt 3pt](5,8.5)(6,8.5)
\put(8.2,6.9){\tiny $\tilde\alpha^+_{2g}$}
\put(8.2,3.8){\tiny $\tilde\alpha^+_{1}$}
\psline[linestyle=dashed, dash=1pt 3pt](8.5,6)(8.5,5)
\put(8.2,1.9){\tiny $\tilde\alpha^+_{2g}$}
\put(8.2,-1.2){\tiny $\tilde\alpha^+_{1}$}
\psline[linestyle=dashed, dash=1pt 3pt](8.5,1)(8.5,0)
\end{picture}
\caption{The bounded regions of the diagram $(\hat{F},\{\tilde\alpha_i^-\},
\{\tilde\alpha_i^+\})$}\label{fig:rectangles}
\end{figure}

Next we consider the Floer differential. One easily checks that the
bounded regions of $\hat{F}$ delimited by the arcs
$\tilde{\alpha}_1^-,\dots,\tilde{\alpha}_{2g}^-$ and $\tilde\alpha_1^+,
\dots,\tilde\alpha_{2g}^+$ are all rectangles; see 
Figure~\ref{fig:rectangles} for a picture of the relevant portion of the
diagram (Figure~\ref{fig:rectangles} is obtained from Figure~\ref{fig:sigma'} 
by cutting open $\hat{F}$
at the back in a manner that splits each arc $\tilde{\alpha}^\pm_i$ at the 
branch point $q_i$; thus, pairs of rectangles which touch 
by a corner at $q_i$ are now separated).

Let us mention in passing that our dictionary between
intersections and strands is easy to understand in terms of Figure
\ref{fig:rectangles}: the columns of the
diagram, from right to left, can be viewed as the $4g$ starting
positions for strands, while the rows, from bottom to top, correspond to
the ending positions. The intersection at column $i$ and row $j$ is then
the strand $\left[\begin{smallmatrix}i\\j\end{smallmatrix}\right]$; however
the intersection at the branch point $q_i$ appears in two places in the
diagram, namely at $(i,i)$ and at $(2g+i,2g+i)$.

Since the diagram $(\hat{F},\{\tilde\alpha_i^-\},\{\tilde\alpha_i^+\})$
is nice, the Floer differential on $CF^*(\tilde{D}_s^-,\tilde{D}_t^+)$ counts
empty embedded rectangles. As in the proof of Proposition
\ref{prop:complex12}, rectangles correspond to resolutions of
crossings in the strand diagram, and the emptiness condition amounts to
the requirement that the resolution does not create any double crossing.
Thus the differentials agree.
\endproof

Next we compare the products in $\F_z$ and $\A(F,k)$. Given $s,t,u\in
\mathcal{S}^{2g}_k$, the composition $\hom(D_s,D_t)\otimes \hom(D_t,D_u)\to
\hom(D_s,D_u)$ in $\F_z$ can be computed by wrapping the thimbles in such a way
that each pair lies in the correct relative position at infinity.
Concretely, we can consider $\tilde{D}_s^-$, $D_t$, and $\tilde{D}_u^+$,
which are products of arcs as in Figure~\ref{fig:sigma'} (with the
understanding that the end points of the $\alpha_i$ all lie on the portion of
$\partial F$ in between the end points of the $\tilde\alpha_i^+$ and
those of the $\tilde\alpha_i^-$).

\begin{prop}\label{prop:prod}
The isomorphism of Proposition \ref{prop:complex} intertwines the product
structures of $\F_z$ and $\A(F,k)$.
\end{prop}

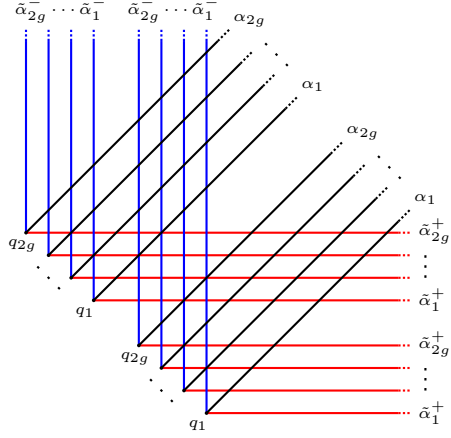
\begin{figure}[t]
\setlength{\unitlength}{3mm}
\begin{picture}(18,18.5)(-1.2,-1.2)
\psset{unit=\unitlength}
\psset{linecolor=red}
\psline(-1,7)(15.5,7) \psline(0,6)(15.5,6) \psline(1,5)(15.5,5) \psline(2,4)(15.5,4)   
\psline(4,2)(15.5,2) \psline(5,1)(15.5,1) \psline(6,0)(15.5,0) \psline(7,-1)(15.5,-1)
\multips(15.5,-1)(0,1){4}{\psline[linestyle=dashed, dash=1pt 1pt](0.5,0)(0,0)}
\multips(15.5,4)(0,1){4}{\psline[linestyle=dashed, dash=1pt 1pt](0.5,0)(0,0)}
\psset{linecolor=blue}
\psline(7,-1)(7,15.5) \psline(6,0)(6,15.5) \psline(5,1)(5,15.5) \psline(4,2)(4,15.5)   
\psline(2,4)(2,15.5) \psline(1,5)(1,15.5) \psline(0,6)(0,15.5) \psline(-1,7)(-1,15.5)
\multips(-1,15.5)(1,0){4}{\psline[linestyle=dashed, dash=1pt 1pt](0,0.5)(0,0)}
\multips(4,15.5)(1,0){4}{\psline[linestyle=dashed, dash=1pt 1pt](0,0.5)(0,0)}
\psset{linecolor=black}
\multips(0,0)(-1,1){4}{\psline(7,-1)(15.5,7.5)}
\multips(0,0)(-1,1){4}{\psline[linestyle=dashed, dash=1pt 1pt](15.5,7.5)(16,8)}
\multips(-5,5)(-1,1){4}{\psline(7,-1)(15.5,7.5)}
\multips(-5,5)(-1,1){4}{\psline[linestyle=dashed, dash=1pt 1pt](15.5,7.5)(16,8)}
\multips(-1,7)(1,-1){4}{\pscircle*(0,0){0.08}}
\multips(4,2)(1,-1){4}{\pscircle*(0,0){0.08}}
\put(-1.8,6.4){\tiny $q_{2g}$}
\put(1.2,3.4){\tiny $q_1$}
\psline[linestyle=dashed, dash=1pt 5pt](-0.5,5.5)(0.5,4.5)
\put(3.2,1.4){\tiny $q_{2g}$}
\put(6.2,-1.6){\tiny $q_1$}
\psline[linestyle=dashed, dash=1pt 5pt](4.5,0.5)(5.5,-0.5)
\put(-1.5,16.7){\tiny $\tilde\alpha^-_{2g}\cdots\tilde\alpha^-_1$}
\put(3.5,16.7){\tiny $\tilde\alpha^-_{2g}\cdots\tilde\alpha^-_1$}
\put(16.4,6.9){\tiny $\tilde\alpha^+_{2g}$}
\put(16.4,3.8){\tiny $\tilde\alpha^+_{1}$}
\psline[linestyle=dashed, dash=1pt 3pt](16.7,6)(16.7,5)
\put(16.4,1.9){\tiny $\tilde\alpha^+_{2g}$}
\put(16.4,-1.2){\tiny $\tilde\alpha^+_{1}$}
\psline[linestyle=dashed, dash=1pt 3pt](16.7,1)(16.7,0)
\put(16.2,8.4){\tiny $\alpha_1$}
\put(13.2,11.4){\tiny $\alpha_{2g}$}
\put(11.2,13.4){\tiny $\alpha_1$}
\put(8.2,16.4){\tiny $\alpha_{2g}$}
\psline[linestyle=dashed, dash=1pt 5pt](14.6,10.4)(15.6,9.4)
\psline[linestyle=dashed, dash=1pt 5pt](9.6,15.4)(10.6,14.4)
\end{picture}
\caption{The bounded regions of the diagram $(\hat{F},\{\tilde\alpha_i^-\},
\{\alpha_i\},\{\tilde\alpha_i^+\})$}\label{fig:triple'}
\end{figure}

\proof
The argument is similar to the proof of Proposition \ref{prop:prod12}.
Namely, the arcs $\tilde\alpha_i^-$, $\alpha_i$ and $\tilde\alpha_i^+$ can
be drawn on $\hat{F}$ so as to form a diagram with non-generic triple 
intersections; after cutting $\hat{F}$ open at the $q_i$, the relevant
portion of the diagram is shown on Figure \ref{fig:triple'}. The triple
intersections can be perturbed as in Figure \ref{fig:triple} right.

By the same argument as in the proof of Proposition \ref{prop:prod12}, the Euler measure of any
2-chain $\phi$ that contributes to the Floer product is equal to $k/4$, 
and the condition $\mu(\phi)=0$ then implies that $\phi$ is disjoint from
the diagonal strata in $\Sym^k(\hat{F})$. Hence
the product counts $k$-tuples of embedded triangles in $\hat{F}$ which
either are disjoint or overlap head-to-tail. Finally, the same argument as
before shows that embedded triangles correspond to strand concatenations,
and that the forbidden overlaps correspond to concatenations that create
double crossings.
\endproof

\begin{prop}\label{prop:higher}
The higher compositions involving the thimbles $D_s$ $(s\in \mathcal{S}^{2g}_k)$
in $\F_z$ are identically zero.
\end{prop}

The proof is identical to that of Proposition \ref{prop:higher12} and simply
relies on a Maslov index calculation to show that there are no rigid discs.

\subsection{Other matchings}\label{ss:moreA}

In \cite{LOT}, Lipshitz, Ozsv\'ath and Thurston construct the algebra
$\A(F,k)$ for an arbitrary pointed matched circle, i.e.\ the $2g$ pairs
of labels assigned to the $4g$ points on the circle need not be in the
configuration $1,\dots,2g,1,\dots,2g$ that we have used throughout.
The only requirement is that the surface obtained by attaching bands
connecting the pairs of identically labelled points and filling in a disc
should have genus $g$ and a single boundary component.

We claim that Theorem \ref{thm:A} admits a natural extension to this more
general setting. Namely, take the configuration of arcs depicted in
Figure \ref{fig:triple'} and view it as lying in a disc $D$, with the $4g$
end points (previously labelled $q_1,\dots,q_{2g},q_1,\dots,q_{2g}$) lying on the boundary. (So there are now $4g$ 
marked points on the boundary of $D$, and $4g$ $\alpha$-arcs 
emanating from them). Next, attach $2g$ bands
to the disc, in such a way that the two ends of each band are attached to
small arcs in $\partial D$ containing end points which carry the same label;
and push the end points into the bands until they come together in pairs.
In this manner one obtains a configuration of $2g$ properly embedded arcs
$\eta_1,\dots,\eta_{2g}$
in a genus $g$ surface with boundary $\mathbb{S}$, as well as their perturbed versions
$\tilde \eta^\pm_i$ which enter in the construction of the partially wrapped Fukaya category.

The $\binom{2g}{k}$ objects of the partially wrapped Fukaya category of the 
$k$-fold symmetric product which correspond to the primitive idempotents of 
$\A(F,k)$ are again products $\Delta_s=\prod_{j\in s}\eta_j$; morphisms,
differentials and products can be understood by cutting $\mathbb{S}$ open
in each band, to obtain diagrams identical to those of Figures
\ref{fig:rectangles} and \ref{fig:triple'} except for a change in labels.
The proof of Theorem \ref{thm:A} then extends without modification.

\section{Generating the partially wrapped category $\F_z$}\label{s:generate}

The goal of this section is to sketch a proof of Theorem
\ref{thm:generate}. The argument is based on a careful analysis of
the relation between the Fukaya category $\F(f_{2g+1,k})$ of the Lefschetz 
fibration $f_{2g+1,k}$ and the partially wrapped category $\F_z$.

In the definition of $\F_z$, we
restricted ourselves to a specific set of noncompact objects with two
useful properties. First, the restriction of $\Re\,f_{2g+1,k}$ to these 
objects is proper and bounded below, and the imaginary part is bounded
by a multiple of the real part. This allows us to view them as objects
of $\F(f_{2g+1,k})$ (after generalizing the notion of admissible
Lagrangian to allow objects to project to a convex angular sector rather
than just to a straight line; this does not significantly affect the 
construction). Second, we only consider products of disjoint properly
embedded arcs, for which the behavior of the flow of $H_\rho$ near infinity
is easy
to understand: namely, in the cylindrical end the flow preserves the
product structure and rotates each arc towards the ray $\vartheta=\pi/2$.
\medskip

{\bf Step 1: The acceleration functor.}
The first ingredient is the existence of a natural $A_\infty$-functor from
$\F(f_{2g+1,k})$ (or rather from a full subcategory whose objects are also
objects of $\F_z$) to $\F_z$; this is a special case of more general ``acceleration''
functors between partially wrapped Fukaya categories, from a less wrapped
category to a more wrapped one. This functor is identity on objects,
and in the simplest cases (e.g.\ for the thimbles $D_s$) it is simply
given by an inclusion of morphism spaces.

Closed exact Lagrangians contained in $\Sym^k(U)$ (as in Definition
\ref{def:Fz}~(1)) are not affected by the flow of $X_{H_\rho}$, and neither
are their intersections with other Lagrangians. Hence, assuming the two
categories $\F(f_{2g+1,k})$ and $\F_z$ are built using the same auxiliary
Hamiltonian perturbations, as far as morphisms to/from compact objects are concerned the 
acceleration functor is simply given by the identity map on Floer complexes. 
Thus we can restrict our attention to noncompact objects.

Let $L=\lambda_1\times\dots\times\lambda_k$ and $L'=\lambda'_1
\times\dots\times\lambda'_k$, 
where $\lambda_1,\dots,\lambda_k$ and $\lambda'_1,\dots,\lambda'_k$ are 
mutually transverse $k$-tuples of disjoint properly embedded arcs as in
Definition~\ref{def:Fz}~(2).
When we view $L$ and $L'$ as objects of $\F(f_{2g+1,k})$, 
morphisms from $L$ to $L'$  are defined
by perturbing $L$ near infinity (in the complement of $U$)
until its slope becomes larger than that of $L'$, i.e.\
by perturbing each $\lambda_i$ in the positive direction to obtain a
new arc $\lambda_i^-$ whose
image under $\pi$ lies closer to the positive imaginary axis than the
images of $\lambda'_1,\dots,\lambda'_k$. (If needed we also choose a
small auxiliary Hamiltonian perturbation to achieve transversality inside $U$).
In other terms, we wrap
the arcs $\lambda_1,\dots,\lambda_k$ by a flow that accumulates
on the two infinite rays $\vartheta=\pi/4$ and $\vartheta=5\pi/4$ (recall
$\vartheta=\frac12 \arg\,\pi(\cdot)$). The complex
$\hom_{\F(f_{2g+1,k})}(L,L')$ is then generated by the intersections of
$L^-=\lambda^-_1\times\dots\times \lambda^-_k$ with $L'$.
(One could also keep perturbing the arcs $\lambda_i$ until they
approach the rays $\vartheta=\pm\pi/2$; this further perturbation does
not affect things in any significant manner,
see Remark \ref{rmk:lesswrapped}.)

The construction of $\hom_{\F_z}(L,L')$ involves the complexes
$CF^*(\phi_{wH_\rho+H'_{L,w}}(L),L')$ for $w\gg 1$. The long-time flow generated
by $H_\rho$ wraps each arc $\lambda_i$ in the positive direction until
it approaches the ray $\vartheta=\pi/2$. Assuming the auxiliary Hamiltonian
perturbations are chosen in the same manner in both categories, the resulting arc 
$\tilde{\lambda}_i^-$ can be viewed as
a perturbation of $\lambda_i^-$ in the cylindrical end $\hat{F}\setminus U$,
further wrapping the arc in the 
positive direction to approach $\vartheta=\pi/2$. Specifically, for each
arc $\lambda_i$ we fix such an isotopy between $\tilde{\lambda}_i^-$ and $\lambda_i^-$, to
be used consistently throughout. With this understood, we set 
$\tilde{L}^-=\tilde\lambda_1^-\times\dots\times\tilde\lambda_k^-$, and
consider the induced isotopy from $L^-$ to $\tilde{L}^-$.
The key point is that the isotopy
from $\lambda_i^-$ to $\tilde{\lambda}_i^-$ only {\it creates}
intersections with the arcs $\lambda'_1,\dots,\lambda'_k$. 
Hence, we can keep track of the intersection points under the isotopy,
which allows us to identify $L^-\cap L'$ with a subset of $\tilde{L}^-
\cap L'$.

\begin{lem}\label{l:subcomplex}
No intersection point created in the isotopy
from $L^-$ to $\tilde{L}^-$ can be the outgoing end of a $J$-holomorphic
strip in $\Sym^k(\hat{F})$ with boundary in $\tilde{L}^-\cup L'$ whose incoming end is a
previously existing intersection point (i.e., one that arises by deforming
a point of $L^-\cap L'$).
\end{lem}

\proof By contradiction, assume such a $J$-holomorphic strip exists. 
Lifting to a branched cover and projecting to $\hat{F}$, we can view it
as a holomorphic map from a bordered Riemann surface
to $\hat{F}$ (with the boundary mapping to the arcs $\tilde\lambda_i^-$ and
$\lambda'_j$). The argument is then purely combinatorial, but is best
understood in terms of the maximum principle applied to the radial
coordinate $r=|\pi|^2$. Namely, after a compactly supported isotopy 
that does not affect intersections, we can assume that, among the points
of $\tilde{\lambda}_i^-\cap \lambda'_j$, those which come from 
$\lambda_i^-\cap \lambda'_j$ have smaller $r$ than the others
(i.e., they lie less far in the cylindrical end; see e.g.\ Figure
\ref{fig:sigma'} right). Moreover, in the cylindrical end the various arcs 
at hand are all graphs (i.e., the angular coordinate $\vartheta$ can be
expressed as a function of the radial coordinate $r$), with the property
that at a point of $\tilde{\lambda}_i^-\cap \lambda'_j$
the slope of $\tilde{\lambda}_i^-$
is always greater than that of $\lambda'_j$. Thus, if an
outgoing strip-like end converges to such an intersection point
(i.e., the boundary of the
holomorphic curve jumps from $\lambda'_j$ to $\tilde{\lambda}_i^-$), then
the radial coordinate $r$ does not have a local maximum. The maximum of $r$
is then necessarily achieved at an incoming strip-like end converging to
an intersection point that lies further in the cylindrical end of $\hat{F}$,
i.e.\ one of the intersections created by the isotopy from $L^-$ to
$\tilde{L}^-$. This contradicts the assumption about the incoming end of
the strip.
\endproof

In other terms, the portion of $CF^*(\tilde{L}^-,L')$
generated by the intersection points that come from $L^-\cap L'$ is a
subcomplex. However the naive map from $CF^*(L^-,L')$ to 
$CF^*(\tilde{L}^-,L')$ obtained by ``following'' the existing generators
through the isotopy is not necessarily a chain map; rather, one should
construct a continuation map using linear cascades as in 
Appendix \ref{appendix}.

More generally, the same argument applies to the $J$-holomorphic discs bounded
by $\ell+1$ Lagrangians obtained by partial wrapping of (mutually
transverse) products of disjoint properly embedded arcs. Namely, using 
appropriate isotopies, we can again view the intersection points
which define morphisms in $\F(f_{2g+1,k})$ as a subset of those which define
morphisms in $\F_z$; the maximum principle applied to the radial
coordinate then implies that a $J$-holomorphic disc whose incoming 
ends all map to previously existing intersection points must have its
outgoing end also mapping to a previously existing intersection point.
In other terms, the wrapping isotopy from $L^-$ to $\tilde{L}^-$
satisfies a property similar to condition (2) in Definition \ref{ass:ass}.
For collections of product Lagrangians which satisfy appropriate
transversality properties, this allows us to use cascades of $J$-holomorphic
discs to build an $A_\infty$-functor whose linear term is given by the
above-mentioned continuation maps.

The behavior of the acceleration functor is significantly simpler if we 
consider the thimbles $D_s$, $s\in \mathcal{S}^{2g}_k$: namely,
in that case an argument similar to that of Lemma~\ref{l:simplify}
implies that the wrapping isotopy does not produce any exceptional
holomorphic discs (of Maslov index less than $2-\ell$), and hence there
are no non-trivial cascades. The acceleration functor is then simply
given by the naive embedding of one Floer complex into the other, obtained
by following the intersection points through the isotopy. Or, to state
things more explicitly via Theorems \ref{thm:A12} and \ref{thm:A}, the
acceleration functor simply corresponds to the obvious embedding of
$\A_{1/2}(F',k)$ as a subalgebra of $\A(F,k)$.

One last property we need to know about the acceleration functor is that
it is cohomologically unital (i.e., the induced functor on cohomology is
unital). When the auxiliary Hamiltonian perturbations are chosen suitably
and identically in both theories, this essentially follows from the fact that
the cohomological unit is given by the ``same'' generator of the 
Floer complex in $\F(f_{2g+1,k})$ and $\F_z$. (The general case is not much
harder). For compact objects contained
in $\Sym^k(U)$ this is clear. For products of properly embedded arcs, the
small-time flow of $H_\rho$ pushes each 
arc slightly off itself in the positive direction at infinity, and choosing
the perturbation suitably we can arrange for
each arc to intersect its pushoff exactly once; the Floer complex
then has a single generator, whose image under the relevant continuation
maps (or, in the case at hand, inclusion of the Floer complex) is a
cohomological unit. (For instance, in the case of the thimbles $D_s$,
this singles out the generator of $\hom(D_s,D_s)$ which consists only of
branch points of $\pi$; that generator turns out to be a strict unit.)
The behavior of the continuation maps which
make up the acceleration functor then ensures unitality of the induced
functor on cohomology.
\medskip

{\bf Step 2: Generation by thimbles.} The next ingredient is Seidel's
result which states that the Fukaya category $\F(f_{2g+1,k})$ is generated
by a collection of Lefschetz thimbles, e.g.\ the $\binom{2g+1}{k}$
product thimbles $D_s$, $s\in \mathcal{S}^{2g+1}_k$ (Theorem 18.24 of
\cite{SeBook}). To be more precise, the only non-compact Lagrangians 
allowed by Seidel are Lefschetz thimbles, so while his result
implies that any compact exact Lagrangian is quasi-isomorphic to a
twisted complex built out of the thimbles $D_s$, the argument in
\cite{SeBook} does not apply to the products of disjoint properly 
embedded arcs that we also wish to allow as objects. On the other hand,
those objects can be shown ``by hand'' to be generated by the $D_s$,
by interpreting arc slides as mapping cones.

Consider $k+1$ disjoint properly embedded arcs
$\lambda_1,\dots,\lambda_k,\lambda'_1$ in $\hat{F}$, all satisfying the 
conditions in Definition \ref{def:Fz}~(2), and such that one extremity
of $\lambda'_1$ lies immediately next to one extremity of $\lambda_1$ in the
cylindrical end $\hat{F}\setminus U$, say in the positive direction from
it. Let $\lambda''_1$ be the arc obtained by sliding $\lambda_1$ along
$\lambda'_1$. Finally,
denote by $\lambda_1^-,\dots,\lambda_k^-$ a collection of arcs obtained
by slightly perturbing $\lambda_1,\dots,\lambda_k$ in the positive 
direction in the cylindrical end, with each $\lambda_i^-$ intersecting 
$\lambda_i$ in a single point $x_i\in U$, and $\lambda_1^-$ intersecting 
$\lambda'_1$ in a single point $x'_1$ which lies near the cylindrical end;
see Figure \ref{fig:arcslide}. Let $L=\lambda_1\times\dots\times \lambda_k$,
$L'=\lambda'_1\times\lambda_2\times\dots\times \lambda_k$, and
$L''=\lambda''_1\times\lambda_2\times\dots\times \lambda_k$.
Then the point $(x'_1,x_2,\dots,x_k)\in
(\lambda_1^-\times\dots\times \lambda_k^-)\cap (\lambda'_1\times\lambda_2
\times\dots\times \lambda_k)$ determines (via the appropriate continuation
map between Floer complexes, to account for the need to further perturb
$L$) an element of $\hom(L,L')$, which we call $u$. We claim:

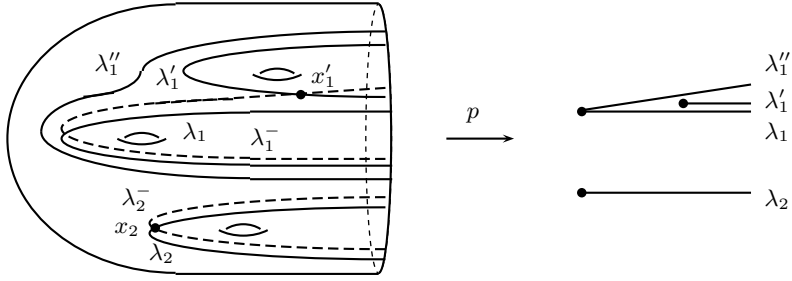
\begin{figure}[t]
\setlength{\unitlength}{9mm}
\begin{picture}(12,4)(0,0)
\psset{unit=\unitlength}
\psellipse[linewidth=0.5pt,linestyle=dashed,dash=2pt 2pt](5.5,2)(0.2,2)
\psellipticarc(5.5,2)(0.2,2){-90}{90}
\psellipticarc(4,2.75)(0.4,0.3){40}{140}
\psellipticarc(4,3.25)(0.6,0.4){-140}{-40}
\psellipticarc(2,1.8)(0.4,0.3){40}{140}
\psellipticarc(2,2.3)(0.6,0.4){-140}{-40}
\psellipticarc(3.5,0.45)(0.4,0.3){40}{140}
\psellipticarc(3.5,0.95)(0.6,0.4){-140}{-40}
\psline(5.5,4)(2.5,4) \psline(5.5,0)(2.5,0)
\psellipticarc(2.5,2)(2.5,2){90}{-90}
\psellipticarc(5.6,3)(3,0.4){90}{-90}
\psellipticarc(5.6,0.6)(3.5,0.4){90}{-90}
\psellipticarc[linestyle=dashed,dash=4pt 2pt](5.6,0.75)(3.5,0.4){90}{-90}
\psellipticarc(3.6,2)(2.8,0.4){90}{-90}
\psline(3.6,2.4)(5.7,2.4)
\psline(3.6,1.6)(5.7,1.6)
\psellipticarc[linestyle=dashed,dash=4pt 2pt](3.6,2.15)(2.8,0.45){120}{-90}
\psline[linestyle=dashed,dash=4pt 2pt](2.2,2.52)(5.6,2.75)
\psline[linestyle=dashed,dash=4pt 2pt](3.6,1.7)(5.6,1.7)
\psellipticarc(5.6,3)(3.6,0.6){90}{180}
\psellipticarc(1,3)(1,0.4){-70}{0}
\psellipticarc(1.8,2.1)(1.3,0.6){110}{180}
\psellipticarc(3.6,2.1)(3.1,0.7){180}{270}
\psline(3.6,1.4)(5.7,1.4)
\put(1.3,3){\small $\lambda''_1$}
\put(2.2,2.8){\small $\lambda'_1$}
\put(2.6,2){\small $\lambda_1$}
\put(3.6,1.9){\small $\lambda_1^-$}
\put(1.7,1.02){\small $\lambda_2^-$}
\put(2.1,0.2){\small $\lambda_2$}
\pscircle*(2.2,0.68){0.07}
\put(1.6,0.6){\small $x_2$}
\pscircle*(4.35,2.65){0.07}
\put(4.5,2.82){\small $x'_1$}
\psline{->}(6.5,2)(7.5,2)
\put(6.8,2.3){\small $p$}
\pscircle*(8.5,1.2){0.07}
\pscircle*(8.5,2.4){0.07}
\pscircle*(10,2.52){0.07}
\psline(8.5,1.2)(11,1.2)
\psline(8.5,2.4)(11,2.4)
\psline(8.5,2.42)(11,2.8)
\psline(10,2.52)(11,2.52)
\put(11.2,1){\small $\lambda_2$}
\put(11.2,2){\small $\lambda_1$}
\put(11.2,2.45){\small $\lambda'_1$}
\put(11.2,3){\small $\lambda''_1$}
\end{picture}
\caption{Sliding $\lambda_1$ along $\lambda'_1$, and the covering $p$} \label{fig:arcslide}
\end{figure}

\begin{lem}\label{l:arcslide}
In $Tw\,\F(f_{2g+1,k})$, $L''$ is quasi-isomorphic to the mapping cone of
$u$.
\end{lem}

\proof
The surface $\hat{F}$ admits a simple branched covering map $p:\hat{F}\to\C$
(i.e., a Lefschetz fibration) with the following properties: (1) the arcs
$\lambda_1,\lambda'_1,\lambda_2,\dots,\lambda_k$ are thimbles for $k+1$ 
of the critical points of $p$ (i.e., lifts of half-lines parallel to the
real axis and connecting critical values $y_1,y'_1,\dots,y_k$ to infinity), with the critical
value for $\lambda'_1$ lying immediately above and very close to the
vanishing path for $\lambda_1$; (2) the monodromies around the critical
points of $p$ corresponding to $\lambda_1$ and $\lambda'_1$ are
two transpositions with one common index, and
sliding the vanishing arc that lifts to $\lambda_1$ around that which lifts
to $\lambda'_1$ yields a new vanishing arc, whose Lefschetz thimble is isotopic to
$\lambda''_1$. See Figure \ref{fig:arcslide} right.
(The covering $p$, whose degree may be very large, can be
built by first projecting a neighborhood of $\lambda_1\cup\lambda'_1$
to $\C$ by a 3:1 map with two branch points, and a neighborhood
of every other $\lambda_i$ by a 2:1 map with a single branch point, and
then extending the map over the rest of~$\hat{F}$). Note that $p$ is not
holomorphic with respect to the given complex structure on $\hat{F}$, but
we can arrange for it to be holomorphic near the branch points.

As in Section \ref{s:lf}, we use $p$ to build a symplectic Lefschetz fibration
$P:\Sym^k(\hat{F})\to\C$, defined by $P([z_1,\dots,z_k])=\sum p(z_i)$
(at least away from the diagonal strata; smoothness requires a
slight modification of $P$ near the diagonal, which is irrelevant for our
purposes). As before, the critical points of $P$ are tuples of distinct
critical points of $p$, and the thimbles associated to straight line
vanishing arcs are just products of the corresponding thimbles for $p$. In particular,
the thimbles associated to the two critical points $[y_1,\dots,y_k]$ and
$[y'_1,y_2,\dots,y_k]$ of $F$ are respectively $L$ and $L'$, and sliding
one vanishing arc over the other one turns $L$ into a product Lagrangian
isotopic to $L''$. (The thimble obtained is not strictly speaking $L''$,
because the sliding operation forces us to consider
vanishing arcs with a small positive slope, so the factors
$\lambda_2,\dots,\lambda_k$ need to be perturbed accordingly.)

It is then a result of Seidel \cite[Proposition 18.23]{SeBook} that
$L''$ is quasi-isomorphic to the mapping cone of the unique generator
of $\hom(L,L')$ in the Fukaya
category of the Lefschetz fibration $P$. (Or, in other terms, the objects
$L$, $L'$ and $L''$ sit in an exact triangle). In order to return to
the Fukaya category of $f_{2g+1,k}$, we observe that the construction
of homomorphisms in the Fukaya category of the Lefschetz fibration $P$ 
requires less wrapping in the positive direction than when we consider
$f_{2g+1,k}$. (In fact, the perturbation needed to bring admissible
Lagrangians into positive position with respect to $P$ can be made
arbitrarily small by choosing $p$ of sufficiently high degree).
Thus, there is again an ``acceleration'' $A_\infty$-functor from the
Fukaya category of $P$ to that of $f_{2g+1,k}$. Taking the image of the
exact triangle involving $L,L',L''$ by this functor (and recalling that
$A_\infty$-functors are exact) yields the result.
\endproof

The other useful fact is that sliding one factor of $L$ over
another factor of $L$ only affects $L$ by a Hamiltonian isotopy.
For instance, if we denote by $\tilde\lambda_1$ the arc obtained by
sliding $\lambda_1$ along $\lambda_2$, then $\tilde{L}=\tilde\lambda_1
\times\lambda_2\times\dots\times \lambda_k$ is Hamiltonian isotopic to $L$.
This follows immediately from the main result in \cite{PerHH}. (More
precisely, the result in \cite{PerHH} is for product tori in symmetric
products of closed surfaces; one can reduce to that case by doubling
$F$ along its boundary to obtain a closed surface and reflecting the arcs
$\lambda_1,\dots,\lambda_k$ to obtain disjoint closed
curves; the arc slide operation then becomes a handle slide and the
result of \cite{PerHH} applies.)

With these two results about arc slides in hand, it is fairly easy to
show that any product of disjoint properly embedded arcs in $\hat{F}$
(satisfying the conditions in Definition~\ref{def:Fz}~(2)) is
quasi-isomorphic in $Tw\,\F(f_{2g+1,k})$ to a complex built out of 
copies of the thimbles $D_s$, $s\in \mathcal{S}^{2g+1}_k$.
Using the exactness of the acceleration $A_\infty$-functor constructed in
Step 1, it now follows that every object of $\F_z$ is quasi-isomorphic in
$Tw\,\F_z$ to a complex built out of the thimbles $D_s$, $s\in
\mathcal{S}^{2g+1}_k$.
\medskip

{\bf Step 3: Eliminating $\alpha_{2g+1}$.}
We now show that, even though all $\binom{2g+1}{k}$ thimbles are needed
to generate $\F(f_{2g+1,k})$, in the case of $\F_z$ it is enough to consider the
$\binom{2g}{k}$ thimbles $D_s$ for $s\subseteq \{1,\dots,2g\}$. For
simplicity, let us assume as in Remark \ref{rmk:A12inF}
that, of the $2g+1$ critical
values $p_j=i\theta_j$ of $\pi:\hat{F}\to\C$, $p_1,\dots,p_{2g}$ lie close 
to the origin along the negative imaginary axis, while $p_{2g+1}$ lies 
further away along the positive imaginary axis; for instance,
let's say that $|\theta_j|<\frac{1}{k}$ for $j\le 2g$, whereas $\theta_{2g+1}>1$.

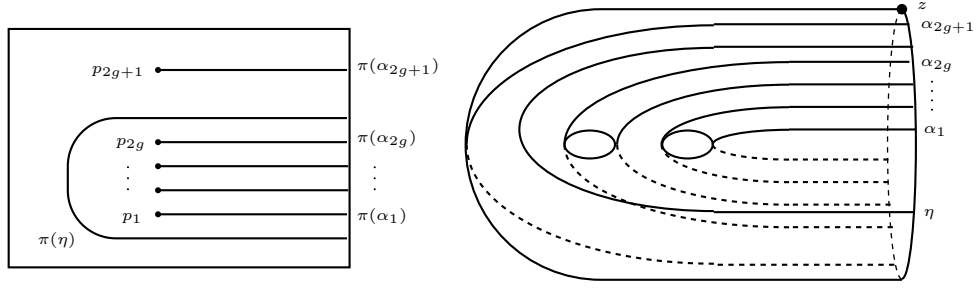
\begin{figure}[t]
\setlength{\unitlength}{8mm}
\begin{picture}(7.5,4.2)(-3,-1.5)
\psset{unit=\unitlength}
\psframe(-2.5,-1.5)(3.2,2.5)
\pscircle*(0,-0.6){0.05}
\pscircle*(0,-0.2){0.05}
\pscircle*(0,0.2){0.05}
\pscircle*(0,0.6){0.05}
\pscircle*(0,1.8){0.05}
\psline(0,-0.6)(3.14,-0.6)
\psline(0,-0.2)(3.16,-0.2)
\psline(0,0.2)(3.16,0.2)
\psline(0,0.6)(3.14,0.6)
\psline(0,1.8)(3.14,1.8)
\put(-0.6,-0.7){\tiny $p_1$}
\put(-0.7,0.55){\tiny $p_{2g}$}
\put(-1.1,1.75){\tiny $p_{2g+1}$}
\put(3.3,-0.7){\tiny $\pi(\alpha_1)$}
\put(3.3,0.6){\tiny $\pi(\alpha_{2g})$}
\put(3.3,1.8){\tiny $\pi(\alpha_{2g+1})$}
\psline[linestyle=dotted](3.6,0.2)(3.6,-0.2)
\psline[linestyle=dotted](-0.5,-0.2)(-0.5,0.2)
\psline[linearc=0.8](3.14,1)(-1.5,1)(-1.5,-1)(3.14,-1)
\put(-2,-1.1){\tiny $\pi(\eta)$}
\end{picture}
\qquad \quad
\setlength{\unitlength}{1cm}
\begin{picture}(5.5,3.5)(-0.5,-1.65)
\psset{unit=\unitlength}
\psellipticarc[linewidth=0.5pt,linestyle=dashed,dash=2pt 2pt](4.5,0)(0.2,1.8){90}{-90}
\psellipticarc(4.5,0)(0.2,1.8){-90}{90}
\psline[linearc=1.8](4.5,1.8)(-1.3,1.8)(-1.3,-1.8)(4.5,-1.8)
\pscircle*(4.5,1.8){0.07} \put(4.7,1.8){\tiny $z$}
\psellipse(0.35,0)(0.35,0.2)
\psellipse(1.65,0)(0.35,0.2)
\psellipticarc(3,0)(1.03,0.21){90}{180}
\psellipticarc(3,0)(1.71,0.51){90}{180}
\psellipticarc(3,0)(2.3,0.81){90}{180}
\psellipticarc(3,0)(3,1.11){90}{180}
\psellipticarc(2,0)(3.3,1.6){90}{180}
\psellipticarc(2,0.2)(2.6,1.1){90}{270}
\psline(2,1.3)(4.65,1.3)
\psline(2,-0.9)(4.65,-0.9)
\psline(3,0.2)(4.68,0.2)
\psline(3,0.5)(4.65,0.5)
\psline(3,0.8)(4.65,0.8)
\psline(3,1.1)(4.6,1.1)
\psline(2,1.6)(4.6,1.6)
\psellipticarc[linestyle=dashed,dash=2pt 2pt](3,0)(1.03,0.215){180}{270}
\psellipticarc[linestyle=dashed,dash=2pt 2pt](3,0)(1.71,0.515){180}{270}
\psellipticarc[linestyle=dashed,dash=2pt 2pt](3,0)(2.3,0.815){180}{270}
\psellipticarc[linestyle=dashed,dash=2pt 2pt](3,0)(3,1.115){180}{270}
\psellipticarc[linestyle=dashed,dash=2pt 2pt](2,0)(3.3,1.6){180}{270}
\psline[linestyle=dashed,dash=2pt 2pt](3.05,-0.2)(4.32,-0.2)
\psline[linestyle=dashed,dash=2pt 2pt](3.05,-0.5)(4.32,-0.5)
\psline[linestyle=dashed,dash=2pt 2pt](3.05,-0.8)(4.35,-0.8)
\psline[linestyle=dashed,dash=2pt 2pt](3.05,-1.1)(4.38,-1.1)
\psline[linestyle=dashed,dash=2pt 2pt](2.05,-1.6)(4.4,-1.6)
\put(4.79,0.15){\tiny $\alpha_1$}
\put(4.75,1.05){\tiny $\alpha_{2g}$}
\put(4.75,1.55){\tiny $\alpha_{2g+1}$}
\put(4.79,-0.95){\tiny $\eta$}
\psline[linestyle=dotted](4.9,0.85)(4.9,0.45)
\end{picture}
\caption{Isotoping $\alpha_{2g+1}$ into $\pi^{-1}(\{\Im w<0\})$}
\label{fig:move2g+1}
\end{figure}

The key observation is that $\alpha_{2g+1}$ can be isotoped, without
crossing the ray $\vartheta=\pi/2$ nor any of the arcs $\alpha_1,\dots,\alpha_{2g}$, to a properly embedded arc $\eta$
contained within the open subset $\pi^{-1}(\{\Im w<0\})$ (which can
be identified with the subsurface $F'$ considered in Section \ref{s:A12});
see Figure \ref{fig:move2g+1}.
Hence, for $s=\{i_1,\dots,i_{k-1},2g+1\}\in \mathcal{S}^{2g+1}_k$ the thimble
$D_s=\alpha_{i_1}\times\dots\times \alpha_{i_{k-1}}\times \alpha_{2g+1}$
can be isotoped without crossing the diagonal nor $\hat{Z}\times
\Sym^{k-1}(\hat{F})$ to the product $\Delta=\alpha_{i_1}\times\dots\times
\alpha_{i_{k-1}}\times\eta$. By construction, $\Delta$ lies within 
the preimage by $f_{2g+1,k}$ of the lower half-plane $\{\Im w<0\}$, 
which can be identified with an open subset of the Lefschetz fibration
$f_{2g,k}$ (see Remark \ref{rmk:A12inF}). The generation argument we have
outlined in Step 2 above then implies that, in $Tw\,\F(f_{2g,k})$, $\Delta$
is quasi-isomorphic to a cone built out of the $\binom{2g}{k}$ thimbles
corresponding to the elements of $\mathcal{S}^{2g}_k$ (we will make this
more explicit below). Recalling that $\F(f_{2g,k})$ embeds as a full
subcategory into $\F(f_{2g+1,k})$, the same result holds in
$Tw\,\F(f_{2g+1,k})$; and hence in $Tw\,\F_z$ as well, via the acceleration
functor of Step 1.

However, the isotopy from $\alpha_{2g+1}$ to $\eta$ does not cross the
ray $\vartheta=\pi/2$. Hence $D_s$ is isotopic to $\Delta$ among
product Lagrangians for which partially wrapped Floer theory
(with respect to the Hamiltonian $H_\rho$) is well-defined, and
the continuation map induced by the isotopy 
(defined using cascades as in Appendix \ref{appendix})
yields a quasi-isomorphism between these two objects
in $\F_z$. (Note here that one could have allowed more general objects in
the category $\F_z$, since the construction of partially wrapped Floer theory
does not require the arcs to project to the right half-plane, as long as
they stay away from the ray $\vartheta=\pi/2$.)
Hence $D_s$ is quasi-isomorphic in $Tw\,\F_z$ to a complex
built out of the thimbles $D_t$, $t\subseteq \{1,\dots,2g\}$. This completes
the proof.

It is not hard to write down explicitly a complex to which
$D_s$ is quasi-isomorphic. Observe that $\eta$ can be obtained by first
sliding $\alpha_1$ along $\alpha_2$ (at the end which lies at the back on 
Figure~\ref{fig:move2g+1} right), then sliding the resulting arc successively 
along $\alpha_3,\dots,\alpha_{2g}$ (at the front of the picture when
sliding over odd $\alpha_i$'s, and at the back when sliding over even
$\alpha_i$'s). For instance, in the case $k=1$, this sequence of arc slides tells us that
$\alpha_{2g+1}$ is quasi-isomorphic to the complex
$$\alpha_1 \stackrel{\left[\begin{smallmatrix}1\\2\end{smallmatrix}\right]}{\longrightarrow}
\alpha_2 \stackrel{\left[\begin{smallmatrix}2g+2\\2g+3\end{smallmatrix}\right]}{\longrightarrow} 
\alpha_3 \stackrel{\left[\begin{smallmatrix}3\\4\end{smallmatrix}\right]}{\longrightarrow} 
\alpha_4 \stackrel{\left[\begin{smallmatrix}2g+4\\2g+5\end{smallmatrix}\right]}{\longrightarrow} 
\dots
\stackrel{\left[\begin{smallmatrix}2g-1\\2g\end{smallmatrix}\right]}{\longrightarrow}\alpha_{2g}$$
(using the notations from $\A(F,k=1)$ to describe the morphisms).
For $k>1$, we can similarly express $\alpha_{i_1}\times\dots\times
\alpha_{i_{k-1}}\times\alpha_{2g+1}$ in terms of the generators by using the same
sequence of arc slides; however, some of the moves now amount to Hamiltonian
isotopies while the others are mapping cones. 

\section{$\widehat{CFA}$ and the pairing theorem}\label{s:CFApairing}

\subsection{Lagrangian correspondences and partially wrapped Fukaya
categories}

As explained in \S \ref{ss:1.1}, work in progress of Lekili and Perutz
\cite{LP} shows that Heegaard-Floer homology can be understood
in terms of quilted Floer homology (cf.\ \cite{WW,WW2}) for Lagrangian 
correspondences between symmetric products. The relevant correspondences
were introduced by Perutz in his thesis \cite{Perutz}; the construction
requires a non-exact perturbation of the K\"ahler form the symmetric
product. 

Given a Riemann surface $\Sigma$, Perutz equips 
$\Sym^k(\Sigma)$ with a K\"ahler form in a class of the form $s\eta_\Sigma+
t\theta_\Sigma$, where $s,t\in \R_+$, $\eta_\Sigma$ is Poincar\'e dual to $\{pt\}\times
\Sym^{k-1}(\Sigma)$, and $\theta_\Sigma-g\eta_\Sigma$ is Poincar\'e dual to
$\sum_1^g a_i\times b_i\times \Sym^{k-2}(\Sigma)$ where $\{a_i,b_i\}$ is a
symplectic basis of $H_1(\Sigma)$ (see \cite{Perutz}). In our case $\Sigma$ is a punctured
Riemann surface, so $\eta_\Sigma$ is trivial, and we choose $[\omega]$ to be a positive 
multiple of $\theta_\Sigma$, or equivalently, a negative multiple of the
first Chern class $c_1(T\Sym^k(\Sigma))=(n+1-g)\eta_\Sigma-\theta_\Sigma$.
Moreover, we arrange for $\omega$ to coincide with the product K\"ahler
form on $\Sigma^k$ away from the diagonal; this ensures that
the Hamiltonian flow used in the construction of the partially wrapped
Fukaya category still preserves the product structure away from the diagonal.


With this understood, let $\gamma$ be a non-separating simple closed 
curve on $\Sigma$, and $\Sigma_\gamma$ the surface obtained from $\Sigma$
by deleting a tubular neighborhood of $\gamma$ and gluing in two discs.
Equip $\Sigma_\gamma$ with a complex structure which agrees with that of $\Sigma$
away from $\gamma$, and equip $\Sym^k(\Sigma)$ and
$\Sym^{k-1}(\Sigma_\gamma)$ with K\"ahler forms $\omega$ and $\omega_\gamma$ chosen as above.

\begin{thm}[Perutz \cite{Perutz}]
The simple closed curve $\gamma$ determines a Lagrangian correspondence
$T_\gamma$ in the product $(\Sym^{k-1}(\Sigma_\gamma)\times \Sym^k(\Sigma),
-\omega_\gamma\oplus \omega)$, canonically up to Hamiltonian isotopy.
\end{thm}

Given $r$ disjoint simple closed curves $\gamma_1,\dots,\gamma_r$, linearly
independent in $H_1(\Sigma)$, we can consider the sequence of
correspondences that arise from successive surgeries along
$\gamma_1,\dots,\gamma_r$. The main properties of these correspondences
(see Theorem A in \cite{Perutz}) imply immediately that their composition
defines an embedded Lagrangian correspondence $T_{\gamma_1,\dots,\gamma_r}$ 
in $\Sym^{k-r}(\Sigma_{\gamma_1,\dots,\gamma_r})\times \Sym^k(\Sigma)$.

When $r=k=g(\Sigma)$, this construction yields
a Lagrangian torus in $\Sym^k(\Sigma)$, which by \cite[Lemma 3.20]{Perutz}
is smoothly isotopic to the product torus $\gamma_1\times\dots\times\gamma_k$; 
Lekili and Perutz show that these two tori are in fact Hamiltonian isotopic \cite{LP}.

Now, consider as in the introduction a 3-manifold $Y$ with connected
boundary $\partial Y\simeq F\cup_{S^1} D^2$ of genus $g$. Viewing $Y$ as a succession
of elementary cobordisms from $D^2$ to $F$ (e.g.\ by considering a Morse
function $f:Y\to [0,1]$ with index 1 and 2 critical points only, with
$f^{-1}(1)=D^2$ and $f^{-1}(0)=F$), $Y$ can be described
by a Heegaard diagram consisting of a once punctured surface $\Sigma$ of 
genus $\bar{g}$ carrying $\bar{g}$ simple closed curves $\beta_1,\dots,
\beta_{\bar{g}}$ (corresponding to the index 2 critical points) and
$\bar{g}-g$ simple closed curves $\alpha_1^c,\dots,\alpha_{\bar{g}-g}^c$
(determined by the index 1 critical points).  These determine respectively
the product torus $T_\beta=\beta_1\times\dots\times \beta_{\bar{g}}\subset
\Sym^{\bar{g}}(\Sigma)$ and a correspondence $\bar{T}_\alpha$ from 
$\Sym^{\bar{g}}(\Sigma)$ to $\Sym^g(F)$. The formal composition
of $T_\beta$ and $\bar{T}_\alpha$ then defines an object $\mathbb{T}_Y$
of the extended Fukaya category $\F^\sharp(\Sym^g(F))$ (in the sense of
Ma'u, Wehrheim and Woodward~\cite{MWW}).

\begin{thm}[Lekili-Perutz \cite{LP}]\label{th:LP}
Up to quasi-isomorphism the object $\mathbb{T}_Y$ 
is independent of the choice of Heegaard diagram for $Y$.
\end{thm}

Even though we are no longer in the exact setting, technical difficulties
in the definition of Floer homology can be avoided by considering {\em balanced}
(also known as {\it Bohr-Sommerfeld monotone}) Lagrangian submanifolds. Namely, equip the
anticanonical bundle $K^{-1}=\det TM^{1,0}$ of $M=\Sym^g(F)$ (resp.\
$\Sym^{\bar{g}}(\Sigma)$) with a connection $\nabla$ whose curvature is 
a constant multiple of the K\"ahler form. 
We say that an orientable Lagrangian submanifold $L$ is balanced with respect to
$\nabla$ if the restriction
of $\nabla$ to $L$ (which is automatically flat) has trivial holonomy, and if moreover the trivialization
of $K^{-1}_{|L}$ induced by a $\nabla$-parallel section is homotopic to the natural
trivialization given by projecting a basis of $TL$ to~$TM^{1,0}$. 

In the context of Heegaard-Floer theory, 
the balancing condition is closely related to admissibility of the Heegaard 
diagram, and can be similarly ensured by a suitable perturbation of the
diagram. Its usefulness is due to the following observation:
if $L$ and $L'$ are balanced, then the symplectic area of a pseudo-holomorphic
strip with boundary on $L,L'$ connecting two given intersection points 
is determined {\it a priori} by its Maslov index (cf.\ \cite[Lemma
4.1.5]{WW2}).  Moreover, the Lagrangians
that we consider do not bound any holomorphic discs, because the inclusion of $L$ into $M$
is injective on fundamental groups and hence $\pi_2(M,L)=\pi_2(M)=0$
(recall that we are considering symmetric products of punctured surfaces);
this prevents bubbling and makes Floer homology well-defined.

These properties allow us to extend the construction of the partially 
wrapped Fukaya category $\F_z$ to this setting, essentially without
modification (considering balanced Lagrangians with $\pi_2(M,L)=0$ instead
of exact ones). Moreover, we can enlarge $\F_z$ to allow sufficiently
well-behaved generalized Lagrangians. Namely, denote by $\F^\sharp_z$ the
$A_\infty$-(pre)category whose objects are

\begin{enumerate}
\item closed balanced Lagrangian tori constructed as products of disjoint,
homologically linearly independent simple closed curves, and generalized 
Lagrangians obtained as
images of such balanced product tori under balanced correspondences
between symmetric products arising from Perutz's construction;\smallskip
\item products of disjoint properly embedded arcs as in Definition
\ref{def:Fz}(2);
\end{enumerate}
with morphisms and compositions defined by partially wrapped Floer theory
using the Hamiltonian $H_\rho$ on $\Sym^g(\hat{F})$ and suitably chosen
small Hamiltonian perturbations.
As in \S \ref{ss:Fz}, we require the closed objects to be contained inside the
bounded subset $\Sym^g(U)$, where $H_\rho$ vanishes; thus these objects
and their intersections with other Lagrangians are not affected by the wrapping.

\begin{prop}\label{prop:generate}
The statement of Theorem \ref{thm:generate} remains valid if $\F_z$ is
replaced by~$\F^\sharp_z$.
\end{prop}

\proof[Sketch of proof]
The general strategy of proof is the same as in \S \ref{s:generate}.
However, we now associate to the Lefschetz fibration
$f_{2g+1,k}$ an extended Fukaya category $\F^\sharp(f_{2g+1,k})$,
whose compact closed objects are
the same balanced generalized Lagrangian submanifolds as in (1)
above (whereas the non-compact objects remain the same as in
$\F(f_{2g+1,k})$). The key point is that Seidel's generation result
still holds in this setting, namely 
$\F^\sharp(f_{2g+1,k})$ is generated by the thimbles
$D_s$, $s\in \mathcal{S}^{2g+1}_k$.

Seidel's argument relies on viewing the Fukaya category of a Lefschetz
fibration as a piece of the $\Z/2$-equivariant Fukaya category of a branched
double cover ramified along a smooth reference fiber (i.e., the pullback by a 
2:1 base change). In our case, we choose the reference fiber to be 
disjoint from $\Sym^k(U)$, e.g.\ we take $f_{2g+1,k}^{-1}(c)$ for $c\in\R_+$
sufficiently large. The thimbles $D_s$, viewed as Lagrangian
discs with boundary in the reference fiber, lift to Lagrangian spheres
$\tilde{D}_s$ in the double cover $\tilde{M}$, while a compact object $L$
lifts to the disjoint union of its two preimages $\tilde{L}=\tilde{L}_+\cup
\tilde{L}_-$. (All these lifts have to be equipped with suitable
$\Z/2$-equivariant structures.) 

Compact generalized Lagrangian submanifolds contained in $\Sym^k(U)$ also lift naturally to the
disjoint union of two compact generalized Lagrangians in $\tilde{M}$. These
behave in the same manner as ordinary Lagrangians. In
particular, the product of the Dehn twists about the Lagrangian spheres
$\tilde{D}_s$ interchanges the two preimages $\tilde{L}_\pm$ of a compact 
object $L$ of $\F^\sharp(f_{2g+1,k})$ (cf.\ \S 18 of \cite{SeBook}). 
Moreover, Seidel's long exact sequence for
Dehn twists generalizes to the quilted setting: namely,
Wehrheim and Woodward show that the graph of the Dehn twist about
$\tilde{D}_s$ fits into an exact triangle in the extended Fukaya category 
of $\tilde{M}\times \tilde{M}$, from which
the long exact sequence follows
(see \S 7 of \cite{WWsequence}). This in turn implies
by the same argument as in \cite[Lemma~18.15 and
Proposition~18.17]{SeBook} that, in $Tw\,\F^\sharp(f_{2g+1,k})$,
compact objects of $\F^\sharp(f_{2g+1,k})$ are quasi-isomorphic to 
twisted complexes built out of the thimbles $D_s$.


With this understood, the rest of the argument works as in \S
\ref{s:generate}. Namely, using the acceleration 
$A_\infty$-functor from $\F^\sharp(f_{2g+1,k})$ to $\F^\sharp_z$
we conclude that $\F^\sharp_z$ is also generated by the thimbles
$D_s$, and the final step (reducing
from $\mathcal{S}^{2g+1}_k$ to $\mathcal{S}^{2g}_k$) is unchanged.
\endproof

\subsection{$\widehat{CFA}$ via Lagrangian correspondences}

We now give a brief outline of the proof of Theorem \ref{th:CFA}. As before,
we represent a 3-manifold $Y$ with parameterized boundary
$\partial Y\simeq F\cup_{S^1} D^2$ by a Heegaard diagram consisting
of a surface $\Sigma$ of genus $\bar{g}\ge g$ with one boundary component,
carrying a base point $z\in \partial \Sigma$ and:
\begin{itemize}
\item $\bar{g}-g$ simple closed curves
$\alpha_1^c,\dots,\alpha_{\bar{g}-g}^c$, which determine a Lagrangian
correspondence $T_\alpha$ from $\Sym^g(F)$ to $\Sym^{\bar{g}}(\Sigma)$ and
the opposite correspondence $\bar{T}_\alpha$ from $\Sym^{\bar{g}}(\Sigma)$
to $\Sym^g(F)$;
\item $2g$ arcs $\alpha_1^a,\dots,\alpha_{2g}^a$, which after surgery
along $\alpha_1^c,\dots,\alpha_{\bar{g}-g}^c$ are assumed to correspond exactly to the arcs
$\alpha_1,\dots,\alpha_{2g}\subset F$ considered in previous sections;
\item $\bar{g}$ simple closed
curves $\beta_1,\dots,\beta_{\bar g}$, which determine a product torus
$T_\beta$ in $\Sym^{\bar{g}}(\Sigma)$.
\end{itemize}
As in the case of $F$, we complete $\Sigma$ to a punctured Riemann
surface $\hat{\Sigma}$, whose cylindrical end can be identified naturally
with that of $\hat{F}$, and consider partially wrapped Floer theory for
balanced product Lagrangians in the symmetric product
$\Sym^{\bar{g}}(\hat\Sigma)$.

Namely, we associate to $\Sym^{\bar{g}}(\hat\Sigma)$ a
partially wrapped category $\bar{\F}^\sharp=\F_z^\sharp(\Sym^{\bar g}
(\hat\Sigma))$, defined similarly to $\F^\sharp_z$ 
except we allow noncompact objects which are balanced products of
mutually disjoint simple closed curves and properly embedded arcs in
$\hat\Sigma$. As before, the simple closed curves are constrained to lie within a bounded
subset $U'$ (corresponding to $U\subset \hat{F}$ after
surgery along the curves $\alpha_i^c$, and assumed to contain all the
closed curves of the Heegaard diagram), while the properly embedded arcs
are constrained to go to infinity in the same manner as in Definition
\ref{def:Fz}(2). 

The Hamiltonian $\bar{H}_\rho$ used to define wrapped Floer homology 
is constructed 
exactly as in \S \ref{ss:Fz}. Namely, away from the diagonal strata it is pulled back 
from a Hamiltonian $\bar{h}_\rho:\hat\Sigma\to \R$ which vanishes over
$U'$, so that the flow of $\bar{H}_\rho$ preserves 
the product structure away from the diagonal and is trivial inside
$\Sym^{\bar{g}}(U')$. Moreover, we pick $\bar{h}_\rho$ to agree with
$h_\rho$ over $\hat{\Sigma}\setminus U'\simeq \hat{F}\setminus U$,
so that the wrapping flow acts similarly on a noncompact object of $\F^\sharp_z$
and on its image under the Lagrangian correspondence~$T_\alpha$.

For $s\in \mathcal{S}^{2g}_g$, we consider the object
$\Delta_{\alpha,s}=\prod\limits_{i\in s} \alpha_i^a\times
\prod\limits_{j=1}^{\bar{g}-g}\alpha_j^c$ of $\bar{\F}^\sharp$.

\begin{lem}
$\Delta_{\alpha,s}$ is Hamiltonian isotopic to the image $T_\alpha(D_s)$
of $D_s\subset \Sym^g(\hat{F})$ under the correspondence $T_\alpha$.
\end{lem}

\noindent
This follows directly from the results in \cite{LP} (since after doubling
$F$ and $\Sigma$ along their boundaries we can reduce to the case of product
tori).

As in \S \ref{ss:Fz}, we choose Hamiltonian perturbations 
for $\Delta_{\alpha,s}$ in such a way that they preserve the product 
structure and commute with the flow of $\bar{H}_\rho$. More specifically,
we choose a Hamiltonian $\bar{h}':\hat\Sigma\to\R$ which agrees with
$h':\hat{F}\to\R$ away from the $\alpha^c_i$, and whose restriction to
each $\alpha^c_i$ is a Morse function with only two critical points,
and we use it to construct a Hamiltonian $\bar{H}'$ on 
$\Sym^{\bar{g}}(\hat{\Sigma})$. This choice of perturbation ensures that
$\hom_{\bar{\F}^\sharp}(\Delta_{\alpha,s},\Delta_{\alpha,t})\simeq
\hom_{\F^\sharp_z}(D_s,D_t)\otimes H^*(T^{\bar{g}-g},\Z_2)$.

By the work of Ma'u-Wehrheim-Woodward \cite{MWW}, the Lagrangian 
correspondences $T_\alpha$ and $\bar{T}_\alpha$ induce $A_\infty$-functors
$\Phi_\alpha:\F^\sharp_z\to \bar{\F}^\sharp$ and $\bar\Phi_\alpha:
\bar{\F}^\sharp\to \F^\sharp_z$. (More precisely, we only have
$A_\infty$-functors between suitable full subcategories, due to the
slightly different restrictions we placed on objects in $\F^\sharp_z$
and $\bar{\F}^\sharp$.) The presence of wrapping Hamiltonians
does not create any significant technical difficulties, since $H_\rho$
and $\bar{H}_\rho$ were chosen compatibly near infinity, and the
$\alpha^c_i$ are contained inside $U'$ where $\bar{h}_\rho$ vanishes
identically.

The functor $\Phi_\alpha$ induces an $A_\infty$-homomorphism from
$\mathcal{A}(F,g)=\bigoplus_{s,t} \hom_{\F^\sharp_z}(D_s,D_t)$ to
$\bar{\mathcal{A}}=\bigoplus_{s,t} \hom_{\bar{\F}^\sharp}(\Delta_{\alpha,s},
\Delta_{\alpha,t})$. In fact, with the choices of perturbations given
above, this map is simply the embedding of $\mathcal{A}(F,g)$
into $\bar{\mathcal{A}}\simeq \mathcal{A}(F,g)\otimes H^*(T^{\bar{g}-g},
\Z_2)$ given by $x\mapsto x\otimes 1$. This makes any
$A_\infty$-module over $\bar{\mathcal{A}}$ into a module over $\A(F,g)$.
With this understood, we have:

\begin{prop}\label{prop:CFA}
$\widehat{CFA}(Y)$ is quasi-isomorphic to $\bigoplus_{s\in
\mathcal{S}^{2g}_g} \hom_{\bar{\F}^\sharp}(T_\beta,\Delta_{\alpha,s})$.
\end{prop}

\proof[Sketch of proof]
Recall from \cite{LOT} that $\widehat{CFA}(Y)$ is generated as a
$\Z_2$-vector spaces by $\bar{g}$-tuples of intersections between
the closed loops $\beta_i$ and the loops and arcs $\alpha^c_i,\alpha^a_i$ 
such that each of $\beta_1,\dots,\beta_{\bar{g}}$ is used exactly once, 
each $\alpha^c_i$ is used exactly once, and each $\alpha^a_i$ is used at
most once. Denoting by $s$ the set of $\alpha^a_i$ which are involved in
the intersection, these tuples correspond exactly to points of
$T_\beta\cap \Delta_{\alpha,s}$. Thus the two sides can be identified
as $\Z_2$-vector spaces.

The $A_\infty$-module structure on $\widehat{CFA}(Y)$ comes from considering
holomorphic curves in $[0,1]\times \R\times \hat{\Sigma}$ with additional
strip-like ends mapping to Reeb chords between the $\alpha^a_i$. Meanwhile, the 
$A_\infty$-module structure on $\bigoplus_s \hom(T_\beta,\Delta_{\alpha,s})$
comes from perturbing the arcs $\alpha^a_i$ by the flow of $\bar{h}_\rho$,
which turns all Reeb chords avoiding the base point $z$ into intersection
points (as seen in \S \ref{ss:pfthmA}), and counting holomorphic discs in
$\Sym^{\bar{g}}(\hat\Sigma)$. There are two steps involved in relating these
two holomorphic curve counts. 

The first step is to view holomorphic discs in $\Sym^{\bar{g}}(\hat\Sigma)$
as curves in $[0,1]\times \R \times \hat{\Sigma}$. This is essentially
identical to Lipshitz's ``cylindrical'' reformulation of Heegaard-Floer
homology. Namely, consider a
holomorphic map $u$ from the disc to $\Sym^{\bar{g}}(\hat\Sigma)$, with
boundary mapping to $T_\beta$ and to suitably wrapped 
copies of objects $\Delta_{\alpha,s_i}$ ($i=1,\dots,k$) in that order
(where the wrapping times $\tau_i$ satisfy $\tau_1\gg \tau_2\gg
\dots \gg \tau_k$). Up to translation there exists a unique biholomorphism $\varphi:D^2\to
(0,1)\times\R$ such 
that the boundary marked points corresponding to the intersections
involving $T_\beta$ are sent to the strip-like ends at $\pm \infty$ while
the boundary marked points corresponding to the intersections between the
perturbed $\Delta_{\alpha,s_i}$'s are sent to points $t_1,\dots,t_{k-1}$
of $\{1\}\times \R$.
Denoting by $\pi:S\to D^2$ a suitable ramified $\bar{g}$:1 covering,
we can turn $u$ into a holomorphic map $\hat{u}:S\to [0,1]\times \R\times
\hat{\Sigma}$, whose first component is $\varphi\circ \pi$ and whose
second component maps the $\bar{g}$ preimages of a point $x\in D^2$ 
to the $\bar{g}$ elements of $u(x)$. The boundary components of $S$
lying above $\{0\}\times\R$ map to the closed curves $\beta_i$, while
the boundary components lying above $\{1\}\times\R$ map to perturbed
copies of the $\alpha$ arcs and curves (switching from one to another
above each $t_i\in \{1\}\times\R$).

(The above discussion assumes that $\Sym^g(\hat\Sigma)$ is equipped 
with the product complex structure; otherwise, the argument
proceeds via perturbation into the class of
``quasi-nearly-symmetric'' almost-complex structures, see \S 13 of
\cite{lipshitz}.)

The second step is to get rid of Hamiltonian perturbations and replace 
the intersection points occurring at the punctures above each $t_i$
by Reeb chords. The main idea is to ``stretch the neck'' near 
$\partial \bar{U}$, i.e.\ deform the complex structure on $\hat{\Sigma}$ 
so that the compact subsurface $\Sigma$ is separated from the region where
the wrapping Hamiltonian $\bar{h}_\rho$ is nonzero by a cylinder of
arbitrarily large modulus. (Equivalently, we do not modify $\hat{\Sigma}$
but change the choice of $\bar{h}_\rho$ so that its support lies further
and further out at infinity.) 
Simultaneously, we turn off the auxiliary
perturbation $\bar{h}'$, so that the $k$ different versions 
of the $\alpha$-arcs and curves converge towards each other in arbitrarily
large subsets of $\hat{\Sigma}$. Under this deformation, holomorphic curves
in $[0,1]\times\R\times\hat{\Sigma}$ converge to multi-stage curves (in
the sense of symplectic field theory). 

Naturally, this homotopy (deforming the complex structure and isotoping
the Lagrangians) needs to be carried out simultaneously and
consistently for all moduli spaces determining the $A_\infty$-module
structure. When deforming a given holomorphic disc, since the homotopy 
affects the values of $t_1,\dots,t_{k-1}$ at which the boundary jumps
from one set of $\alpha$-arcs to the next one, it may happen that two 
of the $t_i$ become equal at some point along the homotopy. When such an
exceptional configuration occurs, the structure maps of the
$A_\infty$-module change; however, since the quasi-isomorphism
class of the $A_\infty$-module $\bigoplus_s \hom(T_\beta,\Delta_{\alpha,s})$
is independent of the choice of the complex structure or admissible Hamiltonian
perturbations, each such modification amounts to composition with a suitable
quasi-isomorphism. Thus, without loss of generality we can assume that the
initial choice of complex structure was already sufficiently stretched to
avoid encountering any such exceptional configurations along the homotopy.

On the other hand, because we are considering rigid holomorphic curves,
and because replacing intersections by Reeb chords preserves
the index of curves in a given homotopy class (under the natural
identification between the respective homotopy groups), a generic choice
of almost-complex structures ensures that the limiting bottom-stage
curve cannot have two Reeb chords lying over the same coordinate in
$[0,1]\times\R$ unless those corresponded to a same input of the
$A_\infty$-module map (i.e., for index reasons, generically
the times $t_i$ cannot become equal in the limit). 

With this understood,
the ``bottom'' stage of the limit curve is again a holomorphic curve in
$[0,1]\times\R\times\hat{\Sigma}$, but the $k$ portions of the boundary over
$\{1\}\times\R$ now all map to unperturbed $\alpha$-arcs and curves.  The
strip-like ends which used to converge to intersection points lying outside
of $\bar{U}$ in the wrapped setting now map to ``Reeb chords'', i.e.\
unbounded strips with boundary on $\alpha$-arcs in the cylindrical end of
$\hat{\Sigma}$, as expected in bordered Heegaard-Floer theory.  Meanwhile,
wherever in the wrapped setting one had a strip-like end converging to an
intersection point lying inside $\bar{U}$ (hence, an intersection between
copies of a {\it same}\/ $\alpha$-arc or curve), the limit curve has a
smooth boundary point, together with a gradient flow trajectory for the
restriction of $\bar{h}'$ to the appropriate arc or loop.  Since we are
considering rigid curves, the limit curve has no intermediate stages, and
the top stage is constant in the $[0,1]\times \R$ factor and consists of
strips in the infinite cylinder $\R\times S^1$ each connecting a Reeb chord
(at the negative end of the cylinder) to the corresponding intersection
point between the wrapped $\alpha$-arcs. 

(In principle, when several Reeb
chords occur at a same $t_i$, the upper stage could also be a more
complicated, non-immersed curve. However, comparing the index formulas,
the presence of branch points in the upper stage causes the index of
the lower stage curve to strictly decrease compared to the case where the
upper stage is ``trivial''. Therefore, since we are considering rigid curves
the upper stage is generically as claimed.)

We claim that the two-stage limit configurations we have just described are in
one-to-one correspondence with the curves used to define the module
structure on $\widehat{CFA}(Y)$. This follows from two observations.

First, the upper stage of the limit curve is uniquely determined by the
bottom stage, since each Reeb chord between two of the arcs $\alpha^a_i$
(not passing over the base point) is connected to the corresponding intersection point between
appropriately wrapped versions of the arcs by a unique rigid holomorphic
strip  in the cylinder
$\R\times S^1$.

Second, whenever an intersection point between an arc
or loop $\eta\in \{\alpha^c_i,\alpha^a_i\}$ and its image under the
perturbation $\bar{h}'$ lies inside $\bar{U}$ and occurs in a generator of
$\Phi_\alpha(\A(F,g))\subset \bar{\A}$, it is necessarily the {\it minimum}
of the restriction of $\bar{h}'$ to $\eta$. Indeed, in the case of
$\alpha^a_i$ the only intersection inside $\bar{U}$ (corresponding to the
pair of horizontal dotted lines $\left[\begin{smallmatrix}i\\
\vphantom{j}\end{smallmatrix}\right]$ in the notation of \cite{LOT}) is
by construction the minimum of $\bar{h}'$ on $\alpha^a_i$; and in the case
of the closed loop $\alpha^c_i$, the claim follows from the description
of the embedding of $\A(F,g)$ into $\bar{\A}$ given just before the
statement of the proposition.

Thus, at each boundary marked point which does not degenerate to
a Reeb chord, the limit curve instead has a smooth boundary on some arc
$\eta$, together with a Morse gradient flow line of $\bar{h}'_{|\eta}$
from the marked point on the boundary of the limit curve to the minimum.
Since every generic point of $\eta$ is connected to the minimum by a
{\it unique} gradient flow line of $\bar{h}'_{|\eta}$, we conclude that
turning the Hamiltonian perturbation $\bar{h}'$ on or off does not
affect the count of holomorphic curves.
\endproof

\begin{prop}\label{prop:adjoint}
The $\A(F,g)$-modules $\bigoplus_s
\hom_{\bar{\F}^\sharp}(T_\beta,\Delta_{\alpha,s})$ and
$\bigoplus_s \hom_{\F^\sharp_z}(\mathbb{T}_Y,D_s)$ are
quasi-isomorphic.
\end{prop}

\begin{rem}
Recalling that $\Delta_{\alpha,s}\simeq \Phi_\alpha(D_s)$ and $\mathbb{T}_Y=
\bar{\Phi}_\alpha(T_\beta)$, this proposition is a special case of a more
general statement, namely that the $A_\infty$-functors
$\Phi_\alpha:\F^\sharp_z\to \bar{\F}^\sharp$ and
$\bar\Phi_\alpha:\bar{\F}^\sharp\to \F^\sharp_z$ induced by the Lagrangian
correspondence $T_\alpha$ are mutually {\em adjoint}. As evident from the
proof, this is a general feature of functors induced by Lagrangian
correspondences and in no way specific to the specific example at hand.
\end{rem}

\proof[Sketch of proof]
The fact that
$\hom_{\bar{\F}^\sharp}(T_\beta,\Phi_\alpha(D_s))$ and
$\hom_{\F^\sharp_z}(\bar\Phi_\alpha(T_\beta),D_s)$ are isomorphic as vector
spaces follows directly
from the definition of extended Fukaya categories, since both are
given by the quilted Floer complex $CF^*(T_\beta,T_\alpha,D_s)$.

\begin{figure}[b]
\setlength{\unitlength}{1.45cm}
\begin{picture}(3,5.4)(0,-0.2)
\psset{unit=\unitlength}
\newrgbcolor{lt1}{0.94 0.94 0.94}
\newrgbcolor{lt2}{0.88 0.88 0.88}
\psframe[linestyle=none,fillstyle=solid,fillcolor=lt2](0,0)(0.9,5)
\pscustom[linestyle=none,fillstyle=solid,fillcolor=lt1]{%
\psline(0.9,0)(0.9,5)
\psellipticarc(3,5)(1,0.5){180}{270}   
\psellipticarc(3,3.5)(1,0.5){90}{180}    
\psellipticarc(3,1.5)(1,0.5){180}{270}
\psellipticarc(3,0)(1,0.5){90}{180}    
}
\psline(0,0)(0,5)
\psline(0.9,0)(0.9,5)
\psellipticarc(3,0)(1,0.5){90}{180}
\psellipticarc(3,1.5)(1,0.5){180}{270}
\psellipticarc(3,3.5)(1,0.5){90}{180}
\psellipticarc(3,5)(1,0.5){180}{270}
\psline[linestyle=dotted](2,2.2)(2,2.8)
\put(0.45,2.5){\makebox(0,0)[ct]{\small $\Sym^{\bar{g}}\hat\Sigma$}}
\put(1.3,4){\small $\Sym^g\hat{F}$}
\put(-0.05,0.3){\makebox(0,0)[rt]{\small $T_\beta$}}
\put(0.95,0.3){\makebox(0,0)[lt]{\small $T_\alpha$}}
\put(-0.05,4.95){\makebox(0,0)[rt]{\small $T_\beta$}}
\put(0.95,4.95){\makebox(0,0)[lt]{\small $T_\alpha$}}
\put(2.2,0.3){\makebox(0,0)[lt]{\small $D_{s_1}$}}
\put(2.2,1.45){\makebox(0,0)[lt]{\small $D_{s_2}$}}
\put(2.2,3.8){\makebox(0,0)[lt]{\small $D_{s_{k-1}}$}}
\put(2.2,4.95){\makebox(0,0)[lt]{\small $D_{s_k}$}}
\end{picture}
\hspace{1.6cm}
\begin{picture}(4.5,5.4)(0,-0.2)
\psset{unit=\unitlength}
\newrgbcolor{lt1}{0.94 0.94 0.94}
\newrgbcolor{lt2}{0.88 0.88 0.88}
\pscustom[linestyle=none,fillstyle=solid,fillcolor=lt2]{%
\psline(0,0)(0,5)
\psellipticarc(2.5,5)(1.5,0.7){180}{270}
\psellipticarc(2.5,3)(1.5,0.7){90}{180}
\psellipticarc(2.5,2)(1.5,0.7){180}{270}
\psellipticarc(2.5,0)(1.5,0.7){90}{180}
}
\pscustom[linestyle=none,fillstyle=solid,fillcolor=lt1]{%
\psellipticarcn(2.5,0)(1,0.4){180}{90}
\psellipticarc(2.5,0)(1.5,0.7){90}{180}
}
\pscustom[linestyle=none,fillstyle=solid,fillcolor=lt1]{%
\psellipticarcn(2.5,2)(1,0.4){270}{180}
\psellipticarcn(2.5,3)(1,0.4){180}{90}
\psellipticarc(2.5,3)(1.5,0.7){90}{180}
\psellipticarc(2.5,2)(1.5,0.7){180}{270}
}
\pscustom[linestyle=none,fillstyle=solid,fillcolor=lt1]{%
\psellipticarcn(2.5,5)(1,0.4){270}{180}
\psellipticarc(2.5,5)(1.5,0.7){180}{270}
}
\psline(0,0)(0,5)
\psellipticarc(2.5,0)(1,0.4){90}{180}
\psellipticarc(2.5,2)(1,0.4){180}{270}
\psellipticarc(2.5,3)(1,0.4){90}{180}
\psellipticarc(2.5,5)(1,0.4){180}{270}
\psellipticarc(2.5,0)(1.5,0.7){90}{180}
\psellipticarc(2.5,2)(1.5,0.7){180}{270}
\psellipticarc(2.5,3)(1.5,0.7){90}{180}
\psellipticarc(2.5,5)(1.5,0.7){180}{270}
\psline[linestyle=dotted](1,2.2)(1,2.8)
\psline[linestyle=dotted](1.5,2.2)(1.5,2.8)
\put(1,1){\makebox(0,0)[cc]{\small $\Sym^{\bar{g}}\hat\Sigma$}}
\put(1,4){\makebox(0,0)[cc]{\small $\Sym^{\bar{g}}\hat\Sigma$}}
\put(1.4,4.5){\small $\Sym^g\hat{F}$}
\put(1.4,0.4){\small $\Sym^g\hat{F}$}
\put(1.4,1.5){\small $\Sym^g\hat{F}$}
\put(1.4,3.4){\small $\Sym^g\hat{F}$}
\put(-0.05,0.3){\makebox(0,0)[rt]{\small $T_\beta$}}
\put(0.7,0.3){\makebox(0,0)[lt]{\small $T_\alpha$}}
\put(0.7,3.3){\makebox(0,0)[lt]{\small $T_\alpha$}}
\put(-0.05,4.95){\makebox(0,0)[rt]{\small $T_\beta$}}
\put(0.75,4.95){\makebox(0,0)[lt]{\small $T_\alpha$}}
\put(0.75,1.95){\makebox(0,0)[lt]{\small $T_\alpha$}}
\put(2.2,0.3){\makebox(0,0)[lt]{\small $D_{s_1}$}}
\put(2.2,1.95){\makebox(0,0)[lt]{\small $D_{s_{j_1}}$}}
\put(2.2,4.95){\makebox(0,0)[lt]{\small $D_{s_k}$}}
\pscustom[linestyle=none,fillstyle=solid,fillcolor=lt1]{%
\pscurve(2.7,0.4)(3.2,0.36)(4,0.04)(4.5,0)
\psellipticarcn(4.5,0.7)(0.5,0.2){270}{180}
\psellipticarcn(4.5,1.3)(0.5,0.2){180}{90}
\pscurve(4.5,2)(4.5,2)(4,1.96)(3.2,1.64)(2.7,1.6)
}
\pscurve(2.7,0.4)(3.2,0.36)(4,0.04)(4.5,0)
\pscurve(2.7,1.6)(3.2,1.64)(4,1.96)(4.5,2)
\psellipticarc[fillstyle=solid,fillcolor=lt2](2.7,1)(0.7,0.3){-90}{90}
\psellipticarc(4.5,1.3)(0.5,0.2){90}{180}
\psellipticarc(4.5,0.7)(0.5,0.2){180}{270}
\psline[linestyle=dotted](4,0.8)(4,1.2)
\put(3.2,1.4){\small $\Sym^g\hat{F}$}
\put(2.5,1){\makebox(0,0)[lc]{\small $\Sym^{\bar g}\hat{\Sigma}$}}
\put(3.2,0.6){\small $T_\alpha$}
\put(4.2,-0.06){\makebox(0,0)[lt]{\small $D_{s_1}$}}
\put(4.2,0.56){\makebox(0,0)[lb]{\small $D_{s_2}$}}
\put(4.2,2.06){\makebox(0,0)[lb]{\small $D_{s_{j_1}}$}}
\pscustom[linestyle=none,fillstyle=solid,fillcolor=lt1]{%
\pscurve(2.7,3.4)(3.2,3.36)(4,3.04)(4.5,3)
\psellipticarcn(4.5,3.7)(0.5,0.2){270}{180}
\psellipticarcn(4.5,4.3)(0.5,0.2){180}{90}
\pscurve(4.5,5)(4.5,5)(4,4.96)(3.2,4.64)(2.7,4.6)
}
\pscurve(2.7,3.4)(3.2,3.36)(4,3.04)(4.5,3)
\pscurve(2.7,4.6)(3.2,4.64)(4,4.96)(4.5,5)
\psellipticarc[fillstyle=solid,fillcolor=lt2](2.7,4)(0.7,0.3){-90}{90}
\psellipticarc(4.5,4.3)(0.5,0.2){90}{180}
\psellipticarc(4.5,3.7)(0.5,0.2){180}{270}
\psline[linestyle=dotted](4,3.8)(4,4.2)
\put(3.2,4.4){\small $\Sym^g\hat{F}$}
\put(2.5,4){\makebox(0,0)[lc]{\small $\Sym^{\bar g}\hat{\Sigma}$}}
\put(3.2,3.6){\small $T_\alpha$}
\put(4.2,4.4){\makebox(0,0)[lt]{\small $D_{s_{k-1}}$}}
\put(4.2,5.06){\makebox(0,0)[lb]{\small $D_{s_k}$}}
\psline[linestyle=dotted](3.5,2.2)(3.5,2.8)
\end{picture}

\caption{The $\A(F,g)$-module structure on $\bigoplus
CF^*(\bar\Phi_\alpha(T_\beta),D_s)$ (left) and the
$\Phi_\alpha(\A(F,g))$-module structure on $\bigoplus
CF^*(T_\beta,\Phi_\alpha(D_s))$ (right)}
\label{fig:adjoint}
\end{figure}
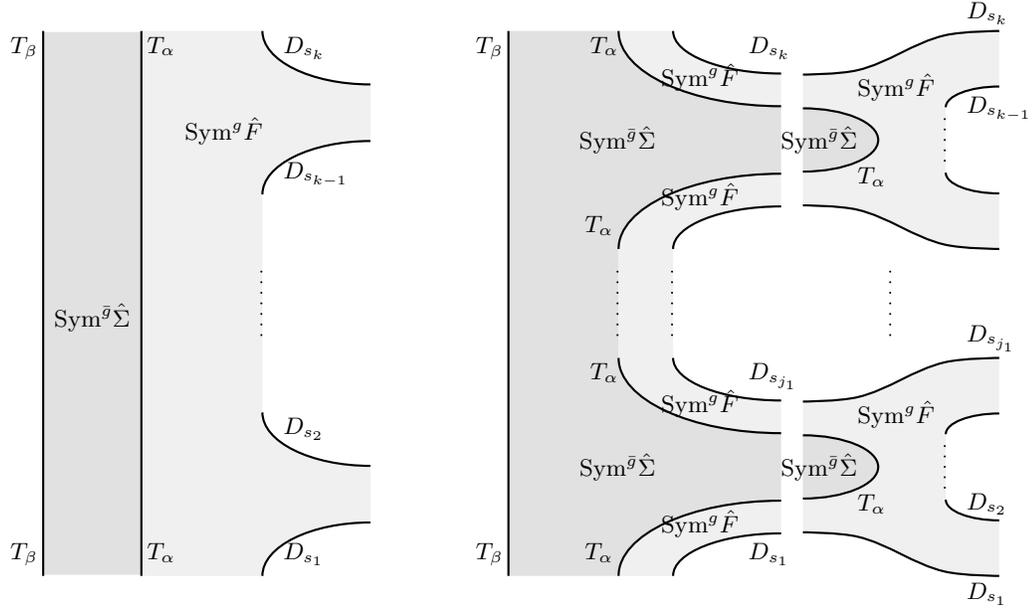

In order to compare the module structures, we describe the relevant
operations graphically in terms of quilted holomorphic curves. In
$\F^\sharp_z$, the $k$-fold product
$$\hom(\bar\Phi_\alpha(T_\beta),D_{s_1})
\otimes\hom(D_{s_1},D_{s_2})\otimes\dots\otimes \hom(D_{s_{k-1}},
D_{s_k})\to \hom(\bar\Phi_\alpha(T_\beta),D_{s_k})$$
is given by a count of quilted holomorphic discs with boundaries on
$T_\beta$ and $D_{s_1},\dots,D_{s_k}$, with a seam mapping to the
correspondence $T_\alpha$, as depicted in the left half of Figure \ref{fig:adjoint}.
On the other hand, the right half of Figure \ref{fig:adjoint} represents 
the quilted discs which contribute to the product operation
$$\hom(T_\beta),\Phi_\alpha(D_{s_1}))
\otimes\hom(\Phi_\alpha(D_{s_1}),\Phi_\alpha(D_{s_{j_1}}))
\otimes\dots \to \hom(T_\beta,\Phi_\alpha(D_{s_k}))$$
in $\bar{\F}^\sharp$, together with the quilted discs which govern
the $A_\infty$-homomorphism from $\A(F,k)$ to $\bar{\A}$ induced by
$\Phi_\alpha$. (Actually, in our case the higher order terms of this
$A_\infty$-homo\-morphism vanish, so the latter quilted discs 
have only one input and look like those in \cite[Figure 10]{WW}.)
The right-hand side picture can be deformed to that on the left-hand side
by gluing the
various components together and moving the seam across; thus the
two module maps agree up to a chain homotopy.

In order to ensure that the various chain homotopies which
arise in this manner for different module maps are consistent, 
the deformation of the seam from the configuration of
Figure \ref{fig:adjoint} left to that of Figure \ref{fig:adjoint} right
needs to carried out along an isotopy chosen in a systematic manner -- 
namely, the isotopies of the seam should depend continuously on the domain,
and compatibly with respect to degenerations.
We will not give an explicit procedure here; however, thinking of the domain of Figure \ref{fig:adjoint}
left as a strip $[0,1]\times\R$ with additional strip-like ends at various 
points $t_i\in \{1\}\times\R$, it is not hard to come up with a construction 
of consistent isotopies of the seam in families over the compactified 
moduli~spaces.
\endproof

\noindent
Finally, Theorem \ref{th:CFA} is a direct corollary of Propositions \ref{prop:CFA}
and \ref{prop:adjoint}.

\subsection{The pairing theorem}
We now sketch the proof of Theorem \ref{th:pairing}. 
Consider a closed 3-manifold
$Y$ which decomposes as the union of two
3-manifolds $Y_1$ and $Y_2$ with $\partial Y_1=-\partial Y_2=F\cup_{S^1} D^2$. 
As in the previous section, Heegaard diagrams for $Y_1$ and for $-Y_2$
allow us to associate to these manifolds two objects
$\mathbb{T}_{Y_1}$ and $\mathbb{T}_{-Y_2}$ of $\F^\sharp_z$. These
generalized Lagrangian submanifolds of $\Sym^g(F)$ can also be constructed
by viewing $Y_1$ and $-Y_2$ as successions of elementary cobordisms between
Riemann surfaces, starting from $D^2$ and ending with $F$. From this
perspective, $Y_2$ is obtained by considering the same sequence of elementary
cobordisms as for $-Y_2$ but in reverse order, starting from $F$ and ending
with $D^2$; thus $Y_2$ defines the opposite correspondence
$\mathbb{T}_{Y_2}=\bar{\mathbb{T}}_{-Y_2}$ from $\Sym^g(F)$ to
$\Sym^0(D^2)=\mathrm{pt}$.

By the work of Lekili and Perutz \cite{LP}, these Lagrangian correspondences allow
us to compute the Heegaard-Floer homology of $Y$, namely $\widehat{CF}(Y)$
is quasi-isomorphic to the quilted Floer complex of the sequence
of correspondences $(\mathbb{T}_{Y_1},\mathbb{T}_{Y_2})$. (Indeed, this
sequence arises from a particular way of
representing the complement of a ball in $Y$ as a cobordism from $D^2$ to
$D^2$; the claim then follows from Theorem \ref{th:LP}, which we now apply in the context
of the manifold $Y\setminus B^3$ with boundary $S^2=D^2\cup_{S^1} D^2$.) Thus, we have
$$\widehat{CF}(Y)\simeq
\hom_{\F^\sharp_z}(\mathbb{T}_{Y_1},\mathbb{T}_{-Y_2}).$$

Next, recall that we have a
contravariant Yoneda functor $\mathcal{Y}:\F^\sharp_z\to \A(F,g)\text{-mod}$,
given on objects by
$$\textstyle \mathbb{L}\mapsto \mathcal{Y}(\mathbb{L})=\bigoplus_s
\hom_{\F_z^\sharp}(\mathbb{L},D_s),$$
and that by Theorem \ref{th:CFA} we have $\widehat{CFA}(Y_1)\simeq
\mathcal{Y}(\mathbb{T}_{Y_1})$ and $\widehat{CFA}(-Y_2)\simeq
\mathcal{Y}(\mathbb{T}_{-Y_2})$.

\begin{prop}
$\mathcal{Y}$ is a cohomologically full and faithful (contravariant) embedding.
\end{prop}

\proof The usual Yoneda embedding of $\F^\sharp_z$ into
$\F^\sharp_z\text{-mod}$ is cohomologically full and faithful (cf.\ e.g.\
\cite[Corollary 2.13]{SeBook}). Moreover, by Proposition \ref{prop:generate}
(the analogue of Theorem \ref{thm:generate} for the extended category
$\F^\sharp_z$), the natural functor from $\F^\sharp_z$-mod to
$\A(F,g)$-mod is an equivalence. The result follows. \endproof

\noindent
Theorem \ref{th:pairing} follows, since we now have
\begin{eqnarray*}
\hom_{\A(F,g)\text{-mod}}(\widehat{CFA}(-Y_2),\widehat{CFA}(Y_1))&\simeq&
\hom_{\A(F,g)\text{-mod}}(\mathcal{Y}(\mathbb{T}_{-Y_2}),
\mathcal{Y}(\mathbb{T}_{Y_1}))\\
&\simeq& \hom_{\F^\sharp_z}(\mathbb{T}_{Y_1},\mathbb{T}_{-Y_2})\\
&\simeq& {\widehat{CF}(Y)}.
\end{eqnarray*}

\appendix
\section{Cascades and partially wrapped Floer theory}\label{appendix}

In this appendix, we outline the construction of the partially 
wrapped Floer complexes and their $A_\infty$-operations. Generally speaking,
things are very similar to the wrapped case defined by Abouzaid and Seidel
in \cite{AS}. However, instead of considering solutions of 
inhomogeneous Cauchy-Riemann equations with Hamiltonian perturbations,
we study trees of genuine $J$-holomorphic curves with boundaries on
perturbed Lagrangian submanifolds. This construction, which was pointed
out to us by Mohammed Abouzaid and is similar to that in Section 10e of
\cite{SeBook},
 allows us both to avoid compactness
issues, and to relate the outcome more directly to
Heegaard-Floer theory. On the other hand, we need to make some assumptions
about the behavior of Lagrangian intersections upon wrapping.

\subsection{Linear cascades and the partially wrapped Floer complex}

Let $(M,\omega)$ be an exact symplectic manifold with convex
contact boundary $(\partial M,\alpha)$, and let $\hat{M}$ be the
exact symplectic manifold obtained by attaching the positive
symplectization $([1,\infty)\times \partial M,d(r\alpha))$ along the
boundary of $M$. We consider a
Hamiltonian function $H_\rho:\hat{M}\to \R$ such that $H_\rho\ge 0$
everywhere and $H_\rho(r,y)=\rho(y)\,r$ on the positive symplectization,
where $\rho:\partial M\to [0,1]$ is a smooth function on the contact
boundary. To a pair of exact Lagrangians $L_1,L_2\subset \hat{M}$ with
cylindrical ends modelled on Legendrian submanifolds of $\partial M\setminus
\rho^{-1}(0)$, we wish to associate a chain complex $\hom(L_1,L_2)$
which behaves as the direct limit for $w\to \infty$ of the Floer
complexes $CF^*(\phi_{wH_\rho}(L_1),L_2)$. Following Abouzaid and Seidel
\cite{AS}, we actually define $\hom(L_1,L_2)$ to be the infinitely
generated complex $\bigoplus_{w=1}^\infty CF^*(\phi_{wH_\rho}(L_1),L_2)[q]$,
or rather the quasi-isomorphic truncation $\bigoplus_{w=m}^\infty
CF^*(\phi_{wH_\rho}(L_1),L_2)[q]$ for some $m\ge 1$ (see
Definition \ref{ass:ass} below),
where the formal variable $q$ has degree $-1$ and satisfies $q^2=0$, equipped with a
differential of the form
\vskip5mm

\begin{equation}\label{eq:CW}
\begin{psmatrix}
CF^*(\phi_{H_\rho}(L_1),L_2) & CF^*(\phi_{2H_\rho}(L_1),L_2) &
CF^*(\phi_{3H_\rho}(L_1),L_2) \\
q\,CF^*(\phi_{H_\rho}(L_1),L_2) &
q\,CF^*(\phi_{2H_\rho}(L_1),L_2) & \dots
\psset{arrows=->,nodesep=4pt}
\everypsbox{\scriptstyle}
\ncline{2,1}{1,1}\tlput{id}
\ncline{2,2}{1,2}\tlput{id}
\nccircle{1,1}{.35cm}^{\delta}
\nccircle{1,2}{.35cm}^{\delta}
\nccircle{1,3}{.35cm}^{\delta}
\nccircle[angleA=180]{2,1}{.35cm}_{\delta}
\nccircle[angleA=180]{2,2}{.35cm}_{\delta}
\ncline{2,1}{1,2}^{\kappa}
\ncline{2,2}{1,3}^{\kappa}
\end{psmatrix}
\end{equation}
\vskip1cm

\noindent
Here $\delta$ is the usual Floer differential, counting index 1
$J$-holomorphic strips with boundary on
$\phi_{wH_\rho}(L_1)$ and $L_2$, while $\kappa$ is a continuation map.
Before we give its definition, let us list the technical assumptions
that will enable our construction to be well-defined.

\begin{defn}\label{ass:ass}
We say that a collection $\{L_i,\ i\in I\}$ of exact Lagrangian submanifolds 
of $\hat{M}$
is transverse with respect to the Hamiltonian $H_\rho$ and the
almost-complex structure $J$ if the following properties hold.
\begin{enumerate}
\item $\phi_{wH_\rho}(L_i)$ is transverse to $L_j$ for all $i,j\in I$
and for all integer values of $w$ greater or equal to some lower
bound $m=m_{i,j}$.\medskip
\item For $w\ge m$, each point of $\phi_{wH_\rho}(L_i)\cap L_j$ lies on a unique 
maximal smooth arc $t\mapsto \gamma(t)$ parametrized by either the whole
interval $[m,\infty)$ or a subinterval of the form $(t_0,\infty)$, such that
$\gamma(t)$ is a transverse intersection of $\phi_{tH_\rho}(L_i)$ 
and $L_j$ for all $t$. In the second case $(t_0>m)$, $\gamma(t)$ 
tends to infinity as $t\to t_0$, and there exists $\epsilon>0$ such that
for $t\in (t_0,t_0+\epsilon)$ no $J$-holomorphic disc can have an outgoing
strip-like end converging to $\gamma(t)\in\phi_{tH_\rho}(L_i)\cap L_j$.
\medskip
\item Given any $i_0,\dots,i_\ell\in I$, and any integers 
$m_{i_{j-1},i_j}\le w_j^-\le
w_j^+$, $j=1,\dots,\ell$ and $0=w_{\ell+1}^-\le w_{\ell+1}^+$, consider
all $J$-holomorphic discs 
in $\hat{M}$ with boundary on the Lagrangian submanifolds
$\phi_{\tau_j H_\rho}(L_{i_j})$ $(0\le j\le\ell)$, where
$\tau_j=\sum_{k=j+1}^{\ell+1}w_k$
and $w_j\in [w_j^-,w_j^+]$, with $\ell+1$ marked points mapping to
given intersections (in the sense of condition (2) above) and 
representing a given relative class $\varphi$. Then the moduli space of
such discs is smooth and of the expected dimension $\mu(\varphi)+\ell-2+\#\{j\,|\,
w_j^-<w_j^+\}$, and all these discs are
regular (as elements of the parametrized moduli space). Moreover,
all nontrivial projections to the real parameters $\tau_j$ and $w_j$ are generic
and transverse to each other with respect to gluing operations
(whenever the outgoing marked point in one moduli space matches with
an incoming marked point in another moduli space).
\end{enumerate}
\end{defn}

Condition (2) can be stated more informally as follows: as $w$ increases
continuously from $m$ to $\infty$, existing intersections between 
$\phi_{wH_\rho}(L_i)$ and $L_j$ persist and remain transverse, whereas
new intersections may be created ``at infinity'' but only provided that,
each time this happens, the newly created intersection cannot be
the outgoing end of any $J$-holomorphic disc. In particular, given $p\in
\phi_{wH_\rho}(L_i)\cap L_j$ and $w'\ge w$, we can associate to $p$
a unique point of $\phi_{w'H_\rho}(L_i)\cap L_j$, which we denote by
$\vartheta_w^{w'}(p)$.

Finally, condition (3)
states that all the moduli spaces of holomorphic discs we will consider
are regular, and behave in the expected manner with respect to gluing;
the precise meaning of the transversality requirement will become clear
in the subsequent discussion.
As usual, we only need this property to hold for 
0- and 1-dimensional moduli spaces in order for the construction to 
be well-defined (while invariance properties also involve 2-dimensional
moduli spaces).

\begin{rem}\label{rmk:transv}
One should keep in mind the following subtlety: when defining higher products,
one sometimes needs to consider cascades in which two of the components
are given by the same data, in which case it is impossible to make the
projections to the time and width parameters transverse, so that condition
(3) does not hold. When it arises, this issue can
be addressed by picking perturbations which depend on the full boundary
data of the cascade (see Definition \ref{def:data}), and not just on the
component under consideration; see Sections 4.7 and 4.9 of 
\cite{AS} for details.
\end{rem}

It is worth mentioning that condition (2) is the key limiting technical
assumption in the approach we adopt. Conditions (1) and (3) can often
be achieved by introducing suitable perturbations 
into the story below (see Remark \ref{rmk:transv}; see also \S \ref{ss:cascadeham}).
On the other hand it is not clear as of this writing how
to construct continuation maps via cascades if (2) does not hold.
When constructing (partially) wrapped Fukaya categories, 
condition (2) usually follows from a finiteness property or from
an appropriate version of the maximum principle.

To keep the notations under control, in the discussion below we will
ignore ignore perturbations; we will also
assume that it is always possible to choose $m=1$ and define
$$\hom(L_i,L_j)=\bigoplus_{w=1}^\infty CF^*(\phi_{wH_\rho}(L_i),L_j)[q].$$
In the general case, we will leave it up to the reader to replace these
complexes by their quasi-isomorphic truncations (restricting to $w\ge m$,
or replacing $H_\rho$ by a multiple).
\medskip

Given two transverse exact Lagrangians $L_1,L_2$ and a positive integer $w$,
we can now define the continuation map $\kappa:CF^*(\phi_{wH_\rho}(L_1),L_2)
\to CF^*(\phi_{(w+1)H_\rho}(L_1),L_2)$ as follows. Given $p\in
\phi_{wH_\rho}(L_1)\cap L_2$ and $q\in \phi_{(w+1)H_\rho}(L_1)\cap L_2$,
a {\em $k$-step linear cascade} from $p$ to $q$ is a sequence of $k$ finite
energy
$J$-holomorphic strips $u_1,\dots,u_k:\R\times [0,1]\to \smash{\hat{M}}$ 
such that:
\begin{itemize}
\item $u_i(\R\times 0)\subset \phi_{w_iH_\rho}(L_1)$ and $u_i(\R\times 1)
\subset L_2$, for some $w_1\le \dots\le w_k$ in the interval
$[w,w+1]$;\medskip

\item denoting by $p_i^\pm\in \phi_{w_iH_\rho}(L_1)\cap L_2$ the intersection
points to which the strips $u_i$ converge at $\pm\infty$, and setting $p_0^+=p$
and $p_{k+1}^-=q$, the points $p_i^+$ and $p_{i+1}^-$ match up
in the sense of property \ref{ass:ass}(2), i.e.\ $p_{i+1}^-=
\vartheta_{w_i}^{w_{i+1}}(p_i^+)$
$\forall\, 0\le i\le k$.
\end{itemize}

\noindent
As a special case we allow $k=0$, i.e.\ the empty sequence of strips,
provided that $q=\vartheta_w^{w+1}(p)$.

We denote by $\mathcal{M}_1^{\{1\}}(L_1,L_2;w;p,q;\varphi)$ the moduli space of
all linear cascades from $p$ to $q$ which represent a given total relative
homotopy class $\varphi$ (the precise definition of the homotopy class
involves completing the broken trajectory to a continuous arc in
the path space using the Hamiltonian isotopy; the details are left to the
reader). The coefficient of $q$ in $\kappa(p)$ is then defined as a
count of rigid linear cascades from $p$ to $q$, i.e.\ those which represent
classes $\varphi$ for which the Maslov index $\mu(\varphi)$ is zero.
By the regularity assumption, these are cascades in which each component
is a Maslov index 0 holomorphic strip at which the linearized $\bar\partial$
operator has a one-dimensional cokernel (``exceptional'' holomorphic
strips).

Linear cascades are a special case of the more general cascades
we will introduce below. Informally, these objects can be understood by
considering
the perturbed holomorphic strips normally used to define Floer 
continuation maps, with a Hamiltonian perturbation term of the form
$\beta(t) X_{H_\rho}$ where the smooth function $\beta:\R\to\R$ tends to $w$ 
at $+\infty$ and $w+1$ at $-\infty$, and taking the limit where the
derivative of $\beta$ tends to zero; it is then reasonable to expect 
that perturbed holomorphic strips converge (in the sense of Gromov 
compactness) to linear cascades.

The algebraic properties of $\kappa$ are determined by the behavior
of one-dimensional moduli spaces of linear cascades.
These moduli spaces are obtained by gluing together various pieces, 
corresponding to different numbers of steps and/or individual homotopy 
classes of the components. Namely,
the part of the boundary of the moduli space of $k$-step cascades
where one of the $k$ components breaks into two $J$-holomorphic
strips is glued with the part of the boundary of the moduli space of 
$k+1$-step cascades where two values $w_i$ and $w_{i+1}$ become equal.
The only remaining boundaries correspond to the cases $w_1=w$ and
$w_k=w+1$, which amounts to breaking off of a $J$-holomorphic strip
contributing to the usual Floer differential $\delta$. Thus
$\kappa\delta=\delta\kappa$ (up to sign), i.e.\ 
$\kappa$ is a chain map between the Floer complexes
$CF^*(\phi_{wH_\rho}(L_1),L_2)$ and $CF^*(\phi_{(w+1)H_\rho}(L_1),L_2)$, 
and the differential on the complex (\ref{eq:CW}) squares to zero.

\subsection{Cascades and $A_\infty$ operations}
The construction of the partially wrapped Fukaya $A_\infty$-category
$\F(M,\rho)$ relies on that of a series of maps\medskip
$$
m_\ell^F:CF^*(\phi_{w_\ell H_\rho}(L_{\ell-1}),L_\ell)\otimes
\dots\otimes CF^*(\phi_{w_1H_\rho}(L_0),L_1)\to 
CF^*(\phi_{w_{out}H_\rho}(L_0),L_\ell),
$$
where $L_0,\dots,L_\ell$ are a transverse collection of exact Lagrangians
($\ell\ge 1$), $F$ is a subset of $\{1,\dots,\ell\}$, $w_1,\dots,w_\ell$ are
positive integers, and $w_{out}=w_1+\dots+w_\ell+|F|$. 

The maps $m_\ell^F$ generalize both the usual Floer-theoretic product operations,
which correspond to $F=\emptyset$, and the continuation map $\kappa$ defined
above, which corresponds to $\ell=1$ and $F=\{1\}$. Up to sign, $m_\ell^F$ is
precisely the part of the $\ell$-fold product operation which maps 
$q^{\epsilon_\ell} CF^*(\phi_{w_\ell H_\rho}(L_{\ell-1}),L_\ell)\otimes
\dots\otimes q^{\epsilon_1} CF^*(\phi_{w_1H_\rho}(L_0),L_1)$ to
$CF^*(\phi_{w_{out}H_\rho}(L_0),L_\ell)$, where $\epsilon_i=1$ if $i\in F$
and $0$ otherwise; see Section 3.8 of \cite{AS}.

We will define the map $m_\ell^F$ differently from the construction in
Section 3 of \cite{AS}, which involves counts of perturbed holomorphic
curves called ``popsicles''. We will instead use cascades of (unperturbed)
holomorphic discs.

\begin{defn}\label{def:data}
We call\/ {\em boundary data} a tuple $(\underline{L};\underline{w},F;
\underline{\smash p},q)$ where:
\begin{itemize}
\item $\underline{L}=(L_0,\dots,L_\ell)$ is a transverse collection of exact
Lagrangian submanifolds;
\item $\underline{w}=(w_1,\dots,w_\ell)\in \R_+^\ell$ are positive real numbers;
\item $F$ is a (possibly empty) subset of $\{1,\dots,\ell\}$; set
$w'_i=w_i+1$ if $i\in F$ and $w'_i=w_i$ otherwise, and
$w_{out}=\sum_{i=1}^\ell w_i+|F|=\sum_{i=1}^\ell w'_i$;
\item $\underline{\smash p}=(p_1,\dots,p_\ell)$, $p_i\in
\phi_{w_iH_\rho}(L_{i-1})\cap L_i$, and $q\in \phi_{w_{out}H_\rho}(L_0)
\cap L_\ell$ are transverse intersection points.
\end{itemize}
A {\em labelled planar tree} for the boundary data 
$(\underline{L};\underline{w},F;\underline{\smash p},q)$
consists of:
\begin{enumerate}
\item a planar tree $\Gamma$ with $\ell+1$ leaves (properly embedded in
$D^2$, with the leaves mapping to the boundary and the other vertices
mapping to the interior), together with a
distinguished leaf called {\em output}; all the edges of\/ $\Gamma$
are oriented so they point towards the output, and the 
components of\/ $D^2\setminus \Gamma$ are numbered by integers $0,\dots,\ell$
and labelled by the Lagrangians $L_0,\dots,L_\ell$ in counterclockwise order
starting  from the output;
\medskip
\item for each vertex $v$ of\/ $\Gamma$ and for each region 
$i$ adjacent to $v$, a ``time'' $\tau_{i,v}\in\R$. These are required to
satisfy the following conditions:
\begin{enumerate}
\item at the output leaf $v_{out}$, 
$\tau_{0,v_{out}}=w_{out}$ and $\tau_{\ell,v_{out}}=0$;
\item at the $i$-th input leaf $v_{in,i}$,
$\tau_{i-1,v_{in,i}}=w_i$ and $\tau_{i,v_{in,i}}=0$;
\item for every directed edge $e=(v^-,v^+)$ separating regions
$i$ and $j$ $(i<j)$, $\tau_{j,v^-}\le \tau_{j,v^+}$, and 
$w_{e,v^-}:=\tau_{i,v^-}-\tau_{j,v^-}\le w_{e,v^+}=\tau_{i,v^+}-\tau_{j,v^+}
\le \sum\limits_{i<k\le j} w'_k$;
\end{enumerate}
\medskip
\item for each vertex $v$ of\/ $\Gamma$ and each edge $e$ adjacent to $v$,
separating two regions $i$ and $j$, a point $p_{e,v}\in
\phi_{\tau_{i,v}H_\rho}(L_i)\cap \phi_{\tau_{j,v}H_\rho}(L_j)$. These
are required to satisfy the following conditions:
\begin{enumerate}
\item at the output leaf, $p_{e,v_{out}}=q$;
\item at the input leaves, $p_{e,v_{in,i}}=p_i$;
\item for every directed edge $e=(v^-,v^+)$ separating regions $i$ and
$j$ $(i<j)$, the points $p_{e,v^-}$ and $p_{e,v^+}$ match up in the
sense of property \ref{ass:ass}(2), i.e.\ $p_{e,v^+}=\phi_{\tau_{j,v^+}
H_\rho}\circ \vartheta_{w_{e,v^-}}^{w_{e,v^+}}\circ \phi_{\tau_{j,v^-}
H_\rho}^{-1}(p_{e,v^-})$.
\end{enumerate}
\end{enumerate}
We denote by $\ell_v$, $\underline{L}_v$, $\underline{\tau}_v$, $\underline{\smash p}_v$
and $q_v$ the number of inputs, Lagrangian submanifolds, times,
incoming and outgoing intersection points
associated to the vertex $v$.
\end{defn}

The elementary building blocks of cascades are $J$-holomorphic
discs with boundaries on the images of given Lagrangian submanifolds
by the Hamiltonian flow generated by $H_\rho$.
Given a transverse collection $\underline{L}=(L_0,\dots,L_\ell)$ of
exact Lagrangians, a tuple of real numbers
$\underline\tau=(\tau_0,\dots,\tau_\ell)\in \R^{\ell+1}$,
a tuple of intersection points $\underline{\smash p}=(p_1,\dots,p_\ell)$, $p_i
\in \phi_{\tau_{i-1}H_\rho}(L_{i-1})\cap \phi_{\tau_iH_\rho}(L_i)$,
$q\in \phi_{\tau_0H_\rho}(L_0)\cap \phi_{\tau_\ell H_\rho}(L_\ell)$, and
a relative homotopy class $\varphi$,
we denote by 
$\mathcal{M}_\ell^\mathrm{hol}(\underline{L};\underline{\tau};\underline{\smash p},q;\varphi)$
the moduli space of $J$-holomorphic maps from the disc with $\ell+1$ (ordered)
boundary marked points to $\hat{M}$, with the boundary arcs mapping to the Lagrangian 
submanifolds $\phi_{\tau_i H_\rho}(L_i)$ and the marked points mapping 
to $p_1,\dots,p_\ell,q$, representing the class $\varphi$.

The Floer product operation 
$$m_\ell=m_\ell^\emptyset:CF^*(\phi_{w_\ell H_\rho}(L_{\ell-1}),L_\ell)\otimes
\dots\otimes CF^*(\phi_{w_1H_\rho}(L_0),L_1)\to 
CF^*(\phi_{w_{out}H_\rho}(L_0),L_\ell)$$
(where $w_{out}=\sum w_i$)
corresponding to the case $F=\emptyset$ differs from a
mere count of $J$-holomorphic discs in that one needs to apply to all inputs the
$A_\infty$-functors intertwining Lagrangian intersection theory for 
$(\phi_{w_iH_\rho}(L_{i-1}),L_i)$ and for $(\phi_{\tau_{i-1}H_\rho}(L_{i-1}),
\phi_{\tau_iH_\rho}(L_i))$, where $\tau_i=\sum_{j>i} w_j$. The standard
way of doing this relies on a Hamiltonian perturbation of the holomorphic
curve equation; instead, the homotopy method leads us to consider cascades
of holomorphic discs. To distinguish the cascades for $F=\emptyset$ from
the more general case (for arbitrary $F$), we will sometimes call them
``plain cascades''.

\begin{defn} \label{def:cascade0}
Given boundary data
$(\underline{L};\underline{w},\emptyset;\underline{\smash p},q)$,
a {\em (plain) cascade} of $J$-holomorphic discs  consists of:
\begin{itemize}
\item a labelled planar tree $(\Gamma,\{\tau_{i,v}\},\{p_{e,v}\})$ for the
boundary data (in the sense of Definition \ref{def:data}); \medskip
\item for each interior vertex $v$ of\/ $\Gamma$, a holomorphic disc
$u_v\in \mathcal{M}_{\ell_v}^\mathrm{hol}(\underline{L}_v;\underline{\tau}_v;
\underline{\smash p}_v,q_v;\varphi_v)$ representing some homotopy class
$[u_v]=\varphi_v$.
\end{itemize}
We denote by
$\mathcal{M}_\ell^\emptyset(\underline{L};\underline{w};\underline{\smash p},
q;\varphi)$
the moduli space of such cascades representing a total homotopy class
$\sum [u_v]=\varphi$.
\end{defn}

\noindent
Note that, since $w_{out}=\sum w_i$, it must be the case that
in condition \ref{def:data}(2)(c)
the equality $w_{e,v^-}=w_{e,v^+}=\sum_{i<k\le j} w_k$ holds for
every directed edge $e=(v^-,v^+)$ of $\Gamma$ separating regions 
$i$ and $j$.

The transversality condition \ref{ass:ass}(3) implies that, when the $w_i$
are positive integers, the moduli space 
$\mathcal{M}^\emptyset_\ell(\underline{L};\underline{w};\underline{\smash
p},q;\varphi)$ is smooth and of the expected dimension, i.e.\
$\mu(\varphi)+\ell-2$. The coefficient of $q$ in $m^\emptyset_\ell(p_\ell,\dots,p_1)$ 
is then defined as a count of cascades in the moduli spaces
$\mathcal{M}^\emptyset_\ell(\underline{L};\underline{w};\underline{\smash
p},q;\varphi)$ for which $\mu(\varphi)=2-\ell$. 

The simplest case is when
the graph $\Gamma$ has a single interior vertex, and the cascade consists
of a single holomorphic disc in
$\mathcal{M}_\ell^{\mathrm{hol}}(\underline{L};\underline{\tau};
\underline{\smash p}^+,q;\varphi)$, where $\tau_i=\sum_{j>i} w_j$ and
$p^+_i=\phi_{\tau_i H_\rho}(p_i)$. More
generally,
the cascades which contribute to $m_\ell^\emptyset$ consist of a ``root component'' which is a rigid
holomorphic disc carrying the output marked point, and other components
which are exceptional holomorphic discs of index $1-\ell_v$ for a component
with $\ell_v$ inputs (indeed, the time parameters $\tau_{i,v}$ for the 
root component are completely fixed, while for the other components
they are only determined up to a simultaneous translation).

\begin{figure}[t]
\setlength{\unitlength}{9mm}
\begin{picture}(14.5,5.2)(-2.8,-2.7)
\psset{unit=\unitlength}
\pscircle(-1.5,1){0.7}
\pscircle(0.9,1){1.1}
\pscircle(3.3,1){0.7}
\pscircle(5.7,1){1.1}
\pscircle(0.9,-1.4){0.7}
\pscircle*(-2.2,1){0.07} \pscircle*(-0.8,1){0.07}
\pscircle*(-0.2,1){0.07} \pscircle*(2,1){0.07} \pscircle*(0.9,-0.1){0.07}
\pscircle*(0.9,-0.7){0.07} \pscircle*(0.9,-2.1){0.07}
\pscircle*(2.6,1){0.07} \pscircle*(4,1){0.07}
\pscircle*(4.6,1){0.07} \pscircle*(6.8,1){0.07} \pscircle*(5.7,-0.1){0.07}
\pscircle*(-2.8,1){0.07}
\pscircle*(0.9,-2.7){0.07} 
\put(6.1,2.2){\makebox(0,0)[cb]{\small $L_0,w_1+w_2+w_3$}}
\put(3.3,1.8){\makebox(0,0)[cb]{\small $L_0,\tau_{0,\mathrm{II}}$}}
\put(0.9,2.2){\makebox(0,0)[cb]{\small $L_0,\tau_{0,\mathrm{III}}$}}
\put(-1.5,1.8){\makebox(0,0)[cb]{\small $L_0,\tau_{0,\mathrm{IV}}$}}
\put(-1.5,0.25){\makebox(0,0)[ct]{\small $L_1,\tau_{1,\mathrm{IV}}$}}
\put(0,0.2){\makebox(0,0)[ct]{\small $L_1$}}
\put(0,-0.2){\makebox(0,0)[ct]{\small $\tau_{1,\mathrm{III}}$}}
\put(0.1,-1.4){\makebox(0,0)[rc]{\small $L_1,\tau_{1,\mathrm{V}}$}}
\put(1.7,-1.4){\makebox(0,0)[lc]{\small $L_2,\tau_{2,\mathrm{V}}$}}
\put(1.95,0.2){\makebox(0,0)[ct]{\small $L_2$}}
\put(1.95,-0.2){\makebox(0,0)[ct]{\small $\tau_{2,\mathrm{III}}$}}
\put(3.3,0.3){\makebox(0,0)[ct]{\small $L_2,\tau_{2,\mathrm{II}}$}}
\put(4.65,0.1){\makebox(0,0)[ct]{\small $L_2,w_3$}}
\put(6.8,0.2){\makebox(0,0)[ct]{\small $L_3,0$}}
\put(-1.5,1.2){\makebox(0,0)[cc]{\small IV}}
\put(-1.5,0.8){\makebox(0,0)[cc]{\small $\mu=0$}}
\put(0.9,1.2){\makebox(0,0)[cc]{\small III}}
\put(0.9,0.8){\makebox(0,0)[cc]{\small $\mu=-1$}}
\put(0.9,-1.2){\makebox(0,0)[cc]{\small V}}
\put(0.9,-1.6){\makebox(0,0)[cc]{\small $\mu=0$}}
\put(3.3,1.2){\makebox(0,0)[cc]{\small II}}
\put(3.3,0.8){\makebox(0,0)[cc]{\small $\mu=0$}}
\put(5.7,1.2){\makebox(0,0)[cc]{\small I}}
\put(5.7,0.8){\makebox(0,0)[cc]{\small $\mu=0$}}
\psline{->}(-2.8,1)(-2.2,1)
\psline{->}(-0.8,1)(-0.2,1)
\psline{->}(2,1)(2.6,1)
\psline{->}(4,1)(4.6,1)
\psline{->}(0.9,-2.7)(0.9,-2.1)
\psline{->}(0.9,-0.7)(0.9,-0.1)
\put(-2.8,1.2){\makebox(0,0)[cb]{\small $p_1$}}
\put(1.1,-2.7){\makebox(0,0)[lc]{\small $p_2$}}
\put(5.7,-0.3){\makebox(0,0)[ct]{\small $p_3$}}
\put(7,1.2){\makebox(0,0)[cb]{\small $q$}}
\put(8,0.8){$\tau_{0,\mathrm{II}}-\tau_{2,\mathrm{II}}=w_1+w_2$,}
\put(8,0.35){$\tau_{0,\mathrm{III}}-\tau_{1,\mathrm{III}}=w_1$,}
\put(8,-0.1){$\tau_{1,\mathrm{III}}-\tau_{2,\mathrm{III}}=w_2$,}
\put(8,-0.55){$\tau_{0,\mathrm{IV}}-\tau_{1,\mathrm{IV}}=w_1$,}
\put(8,-1){$\tau_{1,\mathrm{V}}-\tau_{2,\mathrm{V}}=w_2$}
\put(5,-2){$0\le \tau_{1,\mathrm{IV}}\le \tau_{1,\mathrm{III}}$,}
\put(5,-2.5){$0\le \tau_{2,\mathrm{V}}\le \tau_{2,\mathrm{III}}\le 
\tau_{2,\mathrm{II}}\le w_3$.}
\end{picture}
\caption{A rigid cascade contributing to $m_3^\emptyset$. The arrows
indicate intersections that match via the flow $\phi_{H_\rho}$.}
\label{fig:cascade0}
\end{figure}
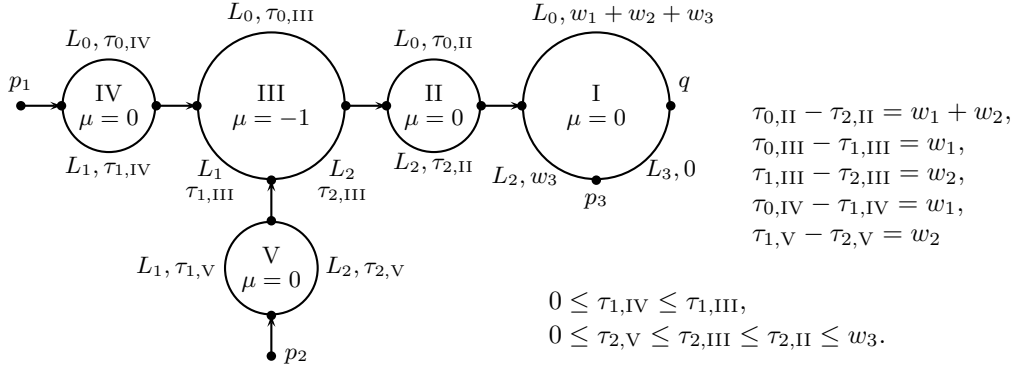

\begin{exm}
By the above discussion, the cascades which contribute to $m_1^\emptyset$
consist of a single index 1 holomorphic disc, so $m_1^\emptyset$ equals the usual Floer
differential $\delta$. The situation is more interesting for $\ell\ge 2$;
for instance,
Figure \ref{fig:cascade0} depicts a rigid plain cascade that contributes to
$m_3^\emptyset$.
\end{exm}

\begin{lem}\label{l:mempty}
The operations $m^\emptyset_\ell$ satisfy the $A_\infty$-relations, i.e.\
$$\sum_{i,k} \,(-1)^{*}\, m_{\ell-k+1}^\emptyset(p_\ell,\dots,p_{i+k+1},m^\emptyset_k(
p_{i+k},\dots,p_{i+1}),p_i,\dots,p_1)=0.$$
\end{lem}

\noindent (Since we work with $\Z/2$ coefficients, we will not worry about
orientations or signs.)

\proof[Sketch of proof]
The argument relies as usual on an analysis of 1-dimensional moduli spaces of cascades.
These moduli spaces are composed of various pieces,
depending on the combinatorial type of the tree $\Gamma$ and the classes
represented by the individual components. At
interior points, exactly one of the components admits a one-parameter
family of deformations, while the others are rigid. 

With one exception, the
portions of the boundary where the non-rigid component breaks into a pair of
$J$-holomorphic disks match exactly with those where the inequality $\tau_{j,v^-}\le
\tau_{j,v^+}$ in condition \ref{def:data}(2)(c) becomes an equality for some directed edge $e=(v^-,v^+)$ 
connecting two interior vertices of $\Gamma$ (one of them carrying the non-rigid
component) and separating regions $i<j$. Accordingly, we glue the various
moduli spaces together along these common boundary strata. 

The exceptional case is when the root component breaks into a pair of rigid discs,
one carrying the $\ell$-th input and the other carrying the output.
In that case we create an edge $e=(v^-,v^+)$ in $\Gamma$ to record the combinatorics of
the breaking, and then split $\Gamma$ along $e$ to obtain a pair of
planar graphs $\Gamma'$, whose root vertex $v^-$ carries the $\ell$-th
input, and $\Gamma''$, whose root vertex $v^+$ carries the original output
(this case has to be treated separately because
$\tau_{\ell,v^-}=\tau_{\ell,v^+}=0$).
One easily checks that the cascade now decomposes into the union of two cascades with
underlying graphs $\Gamma'$ and $\Gamma''$.

The remaining portions of the boundary correspond to the cases where the
inequality $\tau_{j,v^-}\le \tau_{j,v^+}$ becomes an equality
for a directed edge $e=(v^-,v^+)$ that
connects an input leaf to an interior vertex of $\Gamma$. In that case,
we have $\tau_{j,v^+}=0$, and $j$ is necessarily the largest index among
all the regions of $D^2\setminus\Gamma$ adjacent to the vertex $v^+$.
Splitting $\Gamma$ along the outgoing edge from the vertex $v^+$
(and creating a pair of leaves) yields a pair of planar graphs $\Gamma'$
(with root vertex $v^+$) and $\Gamma''$ (with the same root vertex as
$\Gamma$); it is
then easy to check that the cascade decomposes into the union of two
cascades with underlying graphs $\Gamma'$ and $\Gamma''$. 

Conversely,
two cascades such that the outgoing intersection point of one matches
with one of the inputs of the other can be glued to obtain one of the
boundary configurations described above. Thus, the boundary of the moduli space of cascades
can be identified with a union of fibered products of smaller moduli 
spaces of cascades, and the $A_\infty$-relations follow. 
\endproof

We are now ready to define the more general cascades which determine
the operation $m_\ell^F$ for an arbitrary subset $F$ of $\{1,\dots,\ell\}$.

\begin{defn} \label{def:cascadeF}
A {\em cascade} of holomorphic
discs for the boundary
data $(\underline{L};\underline{w},F;\underline{\smash p},q)$
consists of:
\begin{itemize}
\item a labelled planar tree $(\Gamma,\{\tau_{i,v}\},\{p_{e,v}\})$ for the
boundary data, such that for every vertex $v$ of\/ $\Gamma$, the region of greatest index
$j$ among those adjacent to $v$ satisfies $\tau_{j,v}=0$;
\medskip
\item for each interior vertex $v$ of\/ $\Gamma$, an element of
$\mathcal{M}_{\ell_v}^\emptyset(\underline{L}_v;\underline{w}_v;
\underline{\smash p}^-_v,q_v;\varphi_v)$, i.e.\ a plain cascade representing some homotopy class
$\varphi_v$, where $\underline{w}_v$ is the collection of widths
$w_{e,v}$ for the incoming edges at the vertex $v$, and
$p^-_{e,v}=\phi_{\tau_{j,v}H_\rho}^{-1}(p_{e,v})$ for an incoming edge
separating regions $i$ and $j$, $i<j$.
\end{itemize}
We denote by
$\mathcal{M}_\ell^F(\underline{L};\underline{w};\underline{\smash p},
q;\varphi)$
the moduli space of cascades representing a total homotopy class
$\sum \varphi_v=\varphi$.
\end{defn}

The transversality condition \ref{ass:ass}(3) implies that, when the
$w_i$ are positive integers, the moduli space $\mathcal{M}_\ell^F(
\underline{L};\underline{w};\underline{\smash p},q;\varphi)$ is smooth
and of the expected dimension, i.e.\ $\mu(\varphi)+\ell-2+|F|$. The
coefficient of $q$ in $m_\ell^F(p_\ell,\dots,p_1)$ is then defined as a
count of cascades in the moduli spaces $\mathcal{M}_\ell^F(
\underline{L};\underline{w};\underline{\smash p},q;\varphi)$ for which
$\mu(\varphi)=2-\ell-|F|$. Note that the operation $m_1^{\{1\}}$ includes
the empty cascade (where $\Gamma$ has no interior vertices).

Given an interior vertex $v$ of $\Gamma$,
the width parameter $w_{e,v}$ associated to an incoming edge $e$ is 
{\it a priori} free to vary if and only if $e$ can be reached by a directed
path that starts at some input leaf $v_{in,i}$, $i\in F$. Thus, denoting
by $f_v$ the number of such incoming edges at $v$ and by $\ell_v$ the
total number of incoming edges, the dimension of the parametrized moduli
space attached to the vertex $v$ is $\mu(\varphi_v)+\ell_v-2+f_v$.
Hence, the rigid cascades which contribute to $m_\ell^F$ consist of trees 
such that the equality 
\begin{equation}\label{eq:rigidcascade}
\mu(\varphi_v)=2-\ell_v-f_v\end{equation} holds for each interior 
vertex $v$ of $\Gamma$. 

When $F\neq \emptyset$, generic cascades
have the property that $f_v\ge 1$ for every interior vertex $v$, 
i.e.\ each vertex can be reached by a directed path from some input 
leaf $v_{in,i}$, $i\in F$; for otherwise the sum of the individual
dimensions $\mu(\varphi_v)+\ell_v-2+f_v$ turns out to be strictly less
than $\mu(\varphi)+\ell-2+|F|$. (In the case $F=\emptyset$ the same
argument implies that for generic cascades
$\Gamma$ has a single interior vertex,
i.e.\ we are reduced to Definition \ref{def:cascade0}).

\begin{figure}[t]
\setlength{\unitlength}{9mm}
\begin{picture}(14.2,5.2)(-2.8,-2.7)
\psset{unit=\unitlength}
\pscircle(-1.5,1){0.7}
\pscircle(0.9,1){1.1}
\pscircle(3.3,1){0.7}
\pscircle(5.3,1){0.7}
\pscircle(7.7,1){1.1}
\pscircle(10.1,1){0.7}
\pscircle(0.9,-1.4){0.7}
\pscircle(7.7,-1.4){0.7}
\pscircle*(-2.2,1){0.07} \pscircle*(-0.8,1){0.07}
\pscircle*(-0.2,1){0.07} \pscircle*(2,1){0.07} \pscircle*(0.9,-0.1){0.07}
\pscircle*(0.9,-0.7){0.07} \pscircle*(0.9,-2.1){0.07}
\pscircle*(2.6,1){0.07} \pscircle*(4,1){0.07}
\pscircle*(4.6,1){0.07} \pscircle*(6,1){0.07} 
\pscircle*(6.6,1){0.07} \pscircle*(8.8,1){0.07} \pscircle*(7.7,-0.1){0.07}
\pscircle*(7.7,2.1){0.07}
\pscircle*(7.7,-0.7){0.07} \pscircle*(7.2,-1.9){0.07} \pscircle*(8.2,-1.9){0.07}
\pscircle*(9.4,1){0.07} \pscircle*(10.8,1){0.07}
\pscircle*(-2.8,1){0.07} \pscircle*(11.4,1){0.07}
\pscircle*(0.9,-2.7){0.07} \pscircle*(6.5,-2.3){0.07}
\put(10.1,1.8){\makebox(0,0)[cb]{\small $L_0,\tau_{0,\mathrm{I}}$}}
\put(8.6,2.4){\makebox(0,0)[cb]{\small $L_0$}}
\put(8.6,2){\makebox(0,0)[cb]{\small $\tau_{0,\mathrm{II}}$}}
\put(6.8,2.4){\makebox(0,0)[cb]{\small $L_1$}}
\put(6.8,2){\makebox(0,0)[cb]{\small $\tau_{1,\mathrm{II}}$}}
\put(5.3,1.8){\makebox(0,0)[cb]{\small $L_1,\tau_{1,\mathrm{III}}$}}
\put(3.3,1.8){\makebox(0,0)[cb]{\small $L_1,\tau_{1,\mathrm{IV}}$}}
\put(0.9,2.2){\makebox(0,0)[cb]{\small $L_1,\tau_{1,\mathrm{V}}$}}
\put(-1.5,1.8){\makebox(0,0)[cb]{\small $L_1,\tau_{1,\mathrm{VI}}$}}
\put(-1.5,0.2){\makebox(0,0)[ct]{\small $L_2,0$}}
\put(0,0.2){\makebox(0,0)[ct]{\small $L_2$}}
\put(0,-0.2){\makebox(0,0)[ct]{\small $\tau_{2,\mathrm{V}}$}}
\put(0.1,-1.4){\makebox(0,0)[rc]{\small $L_2,\tau_{2,\mathrm{VII}}$}}
\put(1.7,-1.4){\makebox(0,0)[lc]{\small $L_3,0$}}
\put(1.95,0.05){\makebox(0,0)[ct]{\small $L_3,0$}}
\put(3.3,0.2){\makebox(0,0)[ct]{\small $L_3,0$}}
\put(5.3,0.2){\makebox(0,0)[ct]{\small $L_3,0$}}
\put(6.8,0.2){\makebox(0,0)[ct]{\small $L_3$}}
\put(6.8,-0.2){\makebox(0,0)[ct]{\small $\tau_{3,\mathrm{II}}$}}
\put(6.9,-1.4){\makebox(0,0)[rc]{\small $L_3,\tau_{3,\mathrm{VIII}}$}}
\put(7.7,-2.2){\makebox(0,0)[ct]{\small $L_4,w_5$}}
\put(8.5,-1.4){\makebox(0,0)[lc]{\small $L_5,0$}}
\put(8.75,0.05){\makebox(0,0)[ct]{\small $L_5,0$}}
\put(10.1,0.2){\makebox(0,0)[ct]{\small $L_5,0$}}
\put(-1.5,1.2){\makebox(0,0)[cc]{\small VI}}
\put(-1.5,0.8){\makebox(0,0)[cc]{\small $\mu=0$}}
\put(0.9,1.2){\makebox(0,0)[cc]{\small V}}
\put(0.9,0.8){\makebox(0,0)[cc]{\small $\mu=-2$}}
\put(0.9,-1.2){\makebox(0,0)[cc]{\small VII}}
\put(0.9,-1.6){\makebox(0,0)[cc]{\small $\mu=0$}}
\put(3.3,1.2){\makebox(0,0)[cc]{\small IV}}
\put(3.3,0.8){\makebox(0,0)[cc]{\small $\mu=0$}}
\put(5.3,1.2){\makebox(0,0)[cc]{\small III}}
\put(5.3,0.8){\makebox(0,0)[cc]{\small $\mu=0$}}
\put(7.7,1.2){\makebox(0,0)[cc]{\small II}}
\put(7.7,0.8){\makebox(0,0)[cc]{\small $\mu=-3$}}
\put(10.1,1.2){\makebox(0,0)[cc]{\small I}}
\put(10.1,0.8){\makebox(0,0)[cc]{\small $\mu=0$}}
\put(7.7,-1.2){\makebox(0,0)[cc]{\small VIII}}
\put(7.7,-1.6){\makebox(0,0)[cc]{\small $\mu\!=\!-1$}}
\psline[linestyle=dashed]{->}(-2.8,1)(-2.2,1)
\psline[linestyle=dashed]{->}(-0.8,1)(-0.2,1)
\psline[linestyle=dashed]{->}(2,1)(2.6,1)
\psline[linestyle=dashed]{->}(4,1)(4.6,1)
\psline[linestyle=dashed]{->}(6,1)(6.6,1)
\psline[linestyle=dashed]{->}(8.8,1)(9.4,1)
\psline[linestyle=dashed]{->}(10.8,1)(11.4,1)
\psline[linestyle=dashed]{->}(0.9,-2.7)(0.9,-2.1)
\psline[linestyle=dashed]{->}(0.9,-0.7)(0.9,-0.1)
\psline[linestyle=dashed]{->}(7.7,-0.7)(7.7,-0.1)
\psline[linestyle=dashed]{->}(6.5,-2.3)(7.2,-1.9)
\put(-2.8,1.2){\makebox(0,0)[cb]{\small $p_2$}}
\put(7.7,2.3){\makebox(0,0)[cb]{\small $p_1$}}
\put(1.1,-2.7){\makebox(0,0)[lc]{\small $p_3$}}
\put(6.3,-2.3){\makebox(0,0)[rc]{\small $p_4$}}
\put(8.4,-2){\makebox(0,0)[lc]{\small $p_5$}}
\put(11.4,1.2){\makebox(0,0)[cb]{\small $q$}}
\end{picture}
\caption{A rigid cascade contributing to $m_5^{\{2,3,4\}}$. The dotted
lines indicate intersections that match via the maps $\vartheta_w^{w'}$
from \ref{ass:ass}(2).}
\label{fig:cascade}
\end{figure}
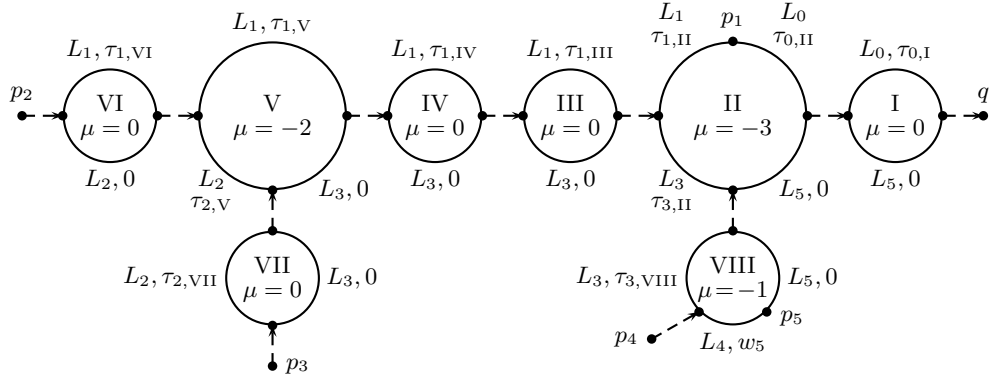

\begin{exm}
Figure \ref{fig:cascade} depicts a rigid cascade that contributes to
$m_5^{\{2,3,4\}}$. Each circle represents either a single holomorphic
disc, or more generally a plain cascade as in Definition \ref{def:cascade0}.
The times $\tau_{i,v}$ satisfy:
\begin{itemize}
\item $\tau_{0,\mathrm{II}}-\tau_{1,\mathrm{II}}=w_1$;
\item $w_2\le \tau_{1,\mathrm{VI}}\le \tau_{1,\mathrm{V}}-\tau_{2,\mathrm{V}}\le
w_2+1$; 
\item $w_3\le \tau_{2,\mathrm{VII}}\le \tau_{2,\mathrm{V}}\le w_3+1$;
\item $\tau_{1,\mathrm{V}}\le \tau_{1,\mathrm{IV}}\le \tau_{1,\mathrm{III}}
\le \tau_{1,\mathrm{II}}-\tau_{3,\mathrm{II}}\le w_2+w_3+2$;
\item $w_4+w_5\le \tau_{3,\mathrm{VIII}}\le \tau_{3,\mathrm{II}}\le w_4+w_5+1$;
\item $\tau_{0,\mathrm{II}}\le \tau_{0,\mathrm{I}}\le w_{out}$.
\end{itemize}
\end{exm}

In order to state the algebraic relation satisfied by the $m_\ell^F$'s,
we first recall the notion of ``admissible cut'' introduced by
Abouzaid and Seidel (cf.\ Section 3.6 of \cite{AS}).

\begin{defn}[\cite{AS}, Definition 3.8]\label{def:admcut}
An {\em admissible cut} of $F\subseteq \{1,\dots,\ell\}$ consists of 
$\ell_+,\ell_-\ge 1$ such that $\ell_-+\ell_+=\ell+1$, a number $i\in
\{1,\dots,\ell_+\}$, and subsets $F_\pm\subseteq \{1,\dots,\ell_\pm\}$
satisfying $|F_-|+|F_+|=|F|$, and with the following property:
$F$ contains all $k\in F_+$ satisfying $k<i$, the numbers $k+\ell_--1$ for
all $k\in F_+$ with $k>i$, and all the numbers
$k+i-1$ for $k\in F_-$. If $i\not\in F_+$ those are all the elements of
$F$, otherwise $F$ has one more element, which lies in the range
$\{i,\dots,i+\ell_--1\}$.
\end{defn}

An admissible cut arises when a cascade decomposes into a pair of cascades
by splitting the graph $\Gamma$ along some edge to obtain a pair of planar
graphs $\Gamma^-$ (carrying input vertices $i$ to $i+\ell_--1$) and
$\Gamma^+$ (carrying input vertices $1$ to $i-1$ and $i+\ell_-$ to $\ell$,
plus a new input vertex arising from the edge that was cut). The elements
of $F$ then decompose in the obvious manner; however we allow ourselves
to delete one element from \hbox{$F\cap \{i,\dots,i+\ell_--1\}$} when forming $F_-$,
in which case $i$ becomes an element of~$F_+$. The width associated to the
cut (i.e., to the output leaf of $\Gamma^-$ and to the $i$-th input leaf
of~$\Gamma^+$) is $w_{cut}=\sum_{k=i}^{i+\ell_--1} w_k+|F_-|$; naturally,
the cut is legal only if the required inequalities
\ref{def:data}(2)(c) hold on either side of the cut, i.e.\ the edge
$e=(v^-,v^+)$ along which $\Gamma$ is split should satisfy $w_{e,v^-}\le
w_{cut}\le w_{e,v^+}$.

The same relation as in \cite[equation (61)]{AS} then holds:

\begin{prop}\label{prop:cuts}
$\sum (-1)^* m_{\ell_+}^{F_+}(p_\ell,\dots,p_{i+\ell_-},
m_{\ell_-}^{F_-}(p_{i+\ell_--1},\dots,p_i),p_{i-1},\dots,p_1)=0$,
where the sum ranges over all admissible cuts.
\end{prop}

\proof[Sketch of proof]
The argument again relies on the study of 1-dimensional moduli spaces
of cascades. As before, these are composed of various pieces (according
to the combinatorial type of the tree $\Gamma$) glued together along
part of their boundaries. For a cascade in a 1-dimensional moduli space, all 
but one of the interior vertices of $\Gamma$ satisfy the equality
$\mu(\varphi_v)=2-\ell_v-f_v$ (i.e., the corresponding plain cascade
is rigid); the remaining interior vertex $v_0$ is associated to a one-dimensional
parametrized moduli space of plain cascades. There are various boundary
strata, corresponding to the following possibilities:
\begin{enumerate}
\item the plain cascade at the vertex $v_0$ breaks up into a pair
of plain cascades, as in the proof of Lemma \ref{l:mempty}; the limiting
cascade is described by a tree with one more vertex;
\item the inequality $w_{e,v^-}\le w_{e,v^+}$ becomes an equality for
some directed edge $e=(v^-,v^+)$ connecting two interior vertices
(one of which is $v_0$);
\item the inequality $w_{e,v^-}\le w_{e,v^+}$ becomes an equality for
some directed edge $e=(v^-,v^+)$ connecting an input leaf $v^-=v_{in,i}$
($i\in F$) to an interior vertex (necessarily $v^+=v_0$);
\item the inequality $w_{e,v^+}\le \sum_{i<k\le j} w_k'$ becomes an equality
for some directed edge $e=(v^-,v^+)$ separating regions $i$ and $j$
(necessarily $v^+=v_0$).
\end{enumerate}
We do not consider the case where $w_{e,v^-}\le w_{e,v^+}$
becomes an equality for $v^-=v_0$ and $v^+=v_{out}$ the outgoing leaf,
since it is a subcase of (4). Moreover, no boundary strata arise from
the inequality $\tau_{j,v^-}\le \tau_{j,v^+}$ (where $e=(v^-,v^+)$ is
a directed edge separating regions $i$ and $j$) becoming an equality:
indeed, Definition \ref{def:cascadeF} implies that $\tau_{j,v^-}=0$, whereas
$\tau_{j,v^+}$ is always zero if $j$ is the greatest index among all regions
adjacent to $v^+$, and always positive and bounded from below otherwise
(due to the positivity of the input width $w_{j+1}$).

Next, we observe that cases (1) and (2) match up exactly, i.e.\ they
correspond to strata along which different moduli spaces are glued
together. (Here it is worth noting that the plain cascades at vertices
$v^-$ and $v^+$ can be glued together to form a single plain cascade
precisely when the widths $w_{e,v^-}$ and $w_{e,v^+}$ are equal, regardless
of whether the times $\tau_{j,v^-}$ and $\tau_{j,v^+}$ match or not.)
We are left with cases (3) and (4), which correspond precisely to the
two types of admissible cuts.

In case (4), we split the tree $\Gamma$ along the edge $e$ to obtain
two trees, $\Gamma^-$ with root vertex $v^-$ and a new output leaf
with width $w_{cut}=\sum_{i<k\le j} w'_k$, and $\Gamma^+$ with
a new input leaf with width $w_{cut}\,(=w_{e,v^+})$.
We obtain a pair of rigid cascades subordinate to an admissible 
cut (with $i\not\in F_+$). (Note: since $i\not\in F_+$, the cut
decreases $f_{v^+}$ by one, which makes the plain cascade at the vertex
$v^+$ rigid after splitting.)

In case (3), namely when $w_{e,v^+}$ becomes equal to $w_i$ for a
directed edge $e$ connecting the $i$-th input leaf to the vertex
$v^+=v_0$, we will find an admissible cut such that the $i$-th input leaf and the
vertex $v_0$ lie within the tree $\Gamma^-$, and the element $i$ is
deleted from $F_-$. Namely, denote by $\hat{e}=(\hat{v}^-,\hat{v}^+)$
the first directed edge encountered along the path from $v_0$ to
the output leaf with the property that $$w_{\hat{e},\hat{v}^-}\le
\sum_{\hat{\textit{\i}}<k\le\hat{\textit{\j}}} w'_k-1\le
w_{\hat{e},\hat{v}^+}$$ where $\hat{\textit{\i}}$ and $\hat{\textit{\j}}$ 
are the labels of the regions separated by $\hat{e}$. Then we split
$\Gamma$ along the edge $\hat{e}$ to obtain two trees: $\Gamma^-$, with
root vertex $\hat{v}^-$ and a new output leaf with width $w_{cut}=
\sum_{\hat{\textit{\i}}<k\le\hat{\textit{\j}}} w'_k-1$, and $\Gamma^+$
with a new input leaf with the same width $w_{cut}$. The cascade then
splits into a pair of rigid cascades subordinate to the relevant cut,
where the label associated to the input leaf $v_{in,i}$ is deleted from
$F_-$ (i.e., $i-\hat{\textit{\i}}+1\not\in F_-$ after relabelling),
whereas the new input is added to $F_+$ (i.e., $\hat{\textit{\i}}\in F_+$).

Finally, each pair of cascades which contributes to the sum in the statement
of the proposition arises precisely once from the splitting of some 
configuration at the boundary of a 1-dimensional moduli space in the 
manner we have described; the result follows.
\endproof

Proposition \ref{prop:cuts} allows us to construct partially wrapped
Fukaya categories (under the assumptions of Definition \ref{ass:ass})
in the same manner as Abouzaid and Seidel \cite{AS}, except we substitute
cascades for popsicles.

We end with the following useful observation:

\begin{lem}\label{l:trivcont}
Let $\{L_i,\ i\in I\}$ be a transverse collection of exact Lagrangian
submanifolds of $\hat{M}$, with the following additional properties:
\begin{enumerate}
\item for all $i,j\in I$, $\phi_{wH_\rho}(L_i)$ is transverse to $L_j$ for
all large enough $w$ ($w\ge m$ for some integer $m=m_{i,j}$), without any intersections being 
created or cancelled;\medskip
\item given any boundary data
$(\underline{L};\underline{w},F;\underline{\smash{p}},q)$ where 
$\underline{L}=(L_{i_0},\dots,L_{i_\ell})$ is a sequence of exact
Lagrangians chosen among the $L_i$, and the widths $w_k$ are large
enough $(w_k\ge m_{i_k,i_{k+1}})$, and given
$\underline{\tau}=(\tau_0,\dots,\tau_\ell)\in \R_{\ge 0}^{\ell+1}$ with 
$\tau_k-\tau_{k+1}=w_k$ and a nontrivial relative class $\varphi$
such that $\mu(\varphi)<2-\ell$, the Lagrangian submanifolds
$\phi_{\tau_kH_\rho}(L_{i_k})$
do not bound any holomorphic disc in the relative class $\varphi$, i.e.\
$\mathcal{M}_\ell^{\mathrm{hol}}(\underline{L};\underline{\tau};
\underline{\smash p}, q;\varphi)=\emptyset$.
\end{enumerate}
Then the operations $m_\ell^F$ are identically zero for $F\neq\emptyset$,
except $m_1^{\{1\}}=\kappa$ which is the natural isomorphism between
Floer complexes induced by identifying intersection points via
the map $\vartheta_w^{w+1}$. 
Thus, up to quasi-isomorphism we can replace the infinitely generated 
complex $\hom(L_i,L_j)$ by $CF^*(\phi_{wH_\rho}(L_i),L_j)$ 
(for any $w\ge m_{i,j}$).
Moreover, the $\ell$-fold product operation
$$m_\ell^\emptyset:CF^*(\phi_{w_\ell H_\rho}(L_{\ell-1}),L_\ell)\otimes
\dots\otimes CF^*(\phi_{w_1H_\rho}(L_0),L_1)\to 
CF^*(\phi_{w_{out}H_\rho}(L_0),L_\ell)$$
simply counts rigid $J$-holomorphic discs in the moduli
spaces $\mathcal{M}_\ell^{\mathrm{hol}}(\underline{L};\underline{\tau};
\underline{\smash p}, q;\varphi)$, where $\tau_i=\sum_{j>i} w_j$ and
we identify the generators of
$CF^*(\phi_{w_iH_\rho}(L_{i-1}),L_i)$ with those of
$CF^*(\phi_{\tau_{i-1}H_\rho}(L_{i-1}),
\phi_{\tau_iH_\rho}(L_i))$ in the obvious manner.
\end{lem}

\proof Recall from the discussion after Definition \ref{def:cascade0}
that rigid plain cascades consist of trees of holomorphic discs
in which the root component is rigid and the other components have index
$1-\ell_v$ where $\ell_v$ is the number of inputs.
However the assuptions give a lower bound
by $2-\ell_v$ on the Maslov index of any nontrivial holomorphic disc.
Thus, rigid plain cascades (those of index $2-\ell$) consist of a single
holomorphic disc, and there are no ``exceptional'' plain cascades
(of index less than $2-\ell$).

Likewise, consider a rigid cascade contributing to $m_\ell^F$ and modelled
after a planar tree~$\Gamma$. Recall from the discussion after Definition
\ref{def:cascadeF} that for each interior vertex $v$ we have a plain cascade
of Maslov index $2-\ell_v-f_v$, where $\ell_v$ is the number of incoming
edges at $v$ and $f_v$ is the number of incoming edges which can be reached
from an input leaf tagged by an element of $F$. Thus the non-existence
of plain cascades of index less than $2-\ell_v$ implies that either
$F=\emptyset$ or $\Gamma$ has no interior vertices (the latter case
corresponds to the empty cascade, which contributes to
$m_1^{\{1\}}=\kappa$).
\endproof

\subsection{Hamiltonian perturbations}\label{ss:cascadeham}
We now modify the above setup by introducing auxiliary Hamiltonian
perturbations in order to make it easier to achieve transversality even
with a degenerate Hamiltonian $H_\rho$. Given two exact Lagrangians
$L_1,L_2$, we fix a family of Hamiltonians $H'_{L_1,L_2,w}$, with the property that
$\phi_{wH_\rho+H'_{L_1,L_2,w}}(L_1)$ is transverse to $L_2$ for large enough $w$, and
we now define $$\hom(L_1,L_2)=\bigoplus_{w=1}^\infty
CF^*(\phi_{wH_\rho+H'_{L_1,L_2,w}}(L_1),L_2)[q].$$ The differential is defined
in terms of linear cascades, exactly as in the unperturbed case. In order 
to define products and higher-order operations
on these complexes, we need to fix homotopies between the relevant
Hamiltonian perturbations, and incorporate them into the definition of
plain cascades (general cascades are then built out of plain cascades as
in the unperturbed case). 

To avoid a lengthy discussion of consistent homotopies between
Hamiltonians, we will restrict ourselves to the case where
the perturbation can be chosen independent of the second Lagrangian, i.e.\
$H'_{L_1,L_2,w}=H'_{L_1,w}$. Thus,
we pick for every Lagrangian $L$ a family of Hamiltonians
$\{H'_{L,\tau}\}_{\tau\ge 0}$, depending smoothly on $\tau$, and with
$H'_{L,0}=0$.
We then replace $\phi_{\tau H_\rho}(L)$ by
$\phi_{\tau H_\rho+H'_{L,\tau}}(L)$ in the construction of plain cascades.

To be more precise, the changes are the following. To start with, we modify Definition \ref{def:data} in the obvious manner, 
so that boundary data now consists of:
\begin{itemize}
\item a collection of exact Lagrangians $\underline{L}=(L_0,\dots,L_\ell)$;
\item positive real numbers $\underline{w}=(w_1,\dots,w_\ell)\in \R^\ell_+$;
\item a subset $F$ of $\{1,\dots,\ell\}$;
\item transverse intersection points $\underline{p}=(p_1,\dots,p_\ell)$ and
$q$, where $$p_i\in
\phi_{w_iH_\rho+H'_{L_{i-1},w_i}}(L_{i-1})\cap L_i\quad\mathrm{and}
\quad q\in \phi_{w_{out}H_\rho+H'_{L_0,w_{out}}}(L_0)\cap L_\ell.$$ 
\end{itemize}

The notion of transversality (Definition \ref{ass:ass})
is modified as follows:

\begin{itemize} \item In condition (1), we now require
$\phi_{(\tau+w)H_\rho+H'_{L_i,\tau+w}}(L_i)$ and $\phi_{\tau
H_\rho+H'_{L_j,\tau}}(L_j)$ to intersect transversely
for all large enough integer values of $w$ and for all $\tau\ge 0$. 
\item Condition (2) again says that, as $w$ increases, new intersections
may be created ``at infinity'', but may not be the outgoing ends of $J$-holomorphic discs. 
\item Condition (3) now requires all relevant moduli spaces of holomorphic discs
with boundaries on the Lagrangians $\phi_{\tau_jH_\rho+H'_{L_j,\tau_j}}(L_j)$ 
to be regular.
\end{itemize}

\noindent
Plain cascades are again built out of $J$-holomorphic discs, taking
the additional Hamiltonian perturbations $H'_{L_j,\tau_j}$ into account.
Given a transverse collection $\underline{L}=(L_0,\dots,L_\ell)$ of
exact Lagrangians, a tuple of real numbers
$\underline\tau=(\tau_0,\dots,\tau_\ell)\in \R^{\ell+1}$,
intersection points $\underline{\smash p}=(p_1,\dots,p_\ell)$ where $p_i
\in \phi_{\tau_{i-1}H_\rho+H'_{L_{i-1},\tau_{i-1}}}(L_{i-1})\cap 
\phi_{\tau_iH_\rho+H'_{L_i,\tau_i}}(L_i)$ and
$q\in \phi_{\tau_0H_\rho+H'_{L_0,\tau_0}}(L_0)\cap \phi_{\tau_\ell
H_\rho+H'_{L_\ell,\tau_\ell}}(L_\ell)$, and
a relative homotopy class $\varphi$,
we now denote by 
$\mathcal{M}_\ell^\mathrm{hol}(\underline{L};\underline{\tau};\underline{\smash p},q;\varphi)$
the moduli space of $J$-holomorphic maps from the disc with $\ell+1$ (ordered)
boundary marked points to $\hat{M}$, with the boundary arcs mapping to the Lagrangian 
submanifolds $\phi_{\tau_i H_\rho+H'_{L_i,\tau_i}}(L_i)$ and the marked points mapping 
to $p_1,\dots,p_\ell,q$, representing the class $\varphi$.

With this change of notation understood, plain cascades are built out of
holomorphic discs exactly as in Definition \ref{def:cascade0}, and
general cascades are defined in terms of plain cascades as in Definition
\ref{def:cascadeF}. With the obvious adaptations, Lemma \ref{l:mempty}, 
Proposition~\ref{prop:cuts} and Lemma \ref{l:trivcont} still hold in this
setting. In particular, we now restate Lemma
\ref{l:trivcont} in the form needed for our purposes:

\begin{lem}\label{l:trivcontham}
Let $\{L_i,\ i\in I\}$ be a transverse collection of exact Lagrangian
submanifolds of $\hat{M}$, with the following additional properties:
\begin{enumerate}
\item for all $i,j\in I$, $\phi_{(\tau+w)H_\rho+H'_{L_i,\tau+w}}(L_i)$ 
is transverse to $\phi_{\tau H_\rho+H'_{L_j,\tau}}(L_j)$ for
all large enough $w$ ($w\ge m=m_{i,j}$) and all $\tau\ge 0$, without any intersections being 
created or cancelled;\medskip
\item given any boundary data
$(\underline{L};\underline{w},F;\underline{\smash{p}},q)$ where 
$\underline{L}=(L_{i_0},\dots,L_{i_\ell})$ is a sequence of exact
Lagrangians chosen among the $L_i$, and the widths $w_k$ are large
enough $(w_k\ge m_{i_k,i_{k+1}})$, and given
$\underline{\tau}=(\tau_0,\dots,\tau_\ell)\in \R_{\ge 0}^{\ell+1}$ with 
$\tau_k-\tau_{k+1}=w_k$ and a nontrivial relative class $\varphi$
such that $\mu(\varphi)<2-\ell$, the Lagrangian submanifolds
$\phi_{\tau_kH_\rho+H'_{L_{i_k},\tau_k}}(L_{i_k})$
do not bound any holomorphic disc in the relative class $\varphi$, i.e.\
$\mathcal{M}_\ell^{\mathrm{hol}}(\underline{L};\underline{\tau};
\underline{\smash p}, q;\varphi)=\emptyset$.
\end{enumerate}
Then the operations $m_\ell^F$ are identically zero for $F\neq\emptyset$,
except $m_1^{\{1\}}=\kappa$ which is the natural isomorphism between
Floer complexes induced by the isotopy.
Thus, up to quasi-isomorphism we can replace the infinitely generated 
complex $\hom(L_i,L_j)$ by $CF^*(\phi_{wH_\rho+H'_{L_i,w}}(L_i),L_j)$ 
(for any $w\ge m_{i,j}$).
Moreover, the $\ell$-fold product operation
\begin{multline*}
m_\ell^\emptyset:CF^*(\phi_{w_\ell H_\rho+H'_{L_{\ell-1},w_\ell}}(L_{\ell-1}),L_\ell)\otimes
\dots\otimes CF^*(\phi_{w_1H_\rho+H'_{L_0,w_1}}(L_0),L_1)\to\\ 
\to CF^*(\phi_{w_{out}H_\rho+H'_{L_0,w_{out}}}(L_0),L_\ell)
\end{multline*}
simply counts rigid $J$-holomorphic discs in the moduli
spaces $\mathcal{M}_\ell^{\mathrm{hol}}(\underline{L};\underline{\tau};
\underline{\smash p}, q;\varphi)$, where $\tau_i=\sum_{j>i} w_j$ and
we identify the generators of
$CF^*(\phi_{w_iH_\rho+H'_{L_{i-1},w_i}}(L_{i-1}),L_i)$ with those of
$CF^*(\phi_{\tau_{i-1}H_\rho+H'_{L_{i-1},\tau_{i-1}}}(L_{i-1}),
\phi_{\tau_iH_\rho+H'_{L_i,\tau_i}}(L_i))$ in the natural manner.
\end{lem}


\begin{thebibliography}{99}
\bibitem{Ab1}
   M. Abouzaid,
   {\sl Homogeneous coordinate rings and mirror symmetry for toric varieties},
   Geom. Topol. {\bf 10} (2006), 1097--1157.
\bibitem{AS}
   M. Abouzaid, P. Seidel, {\sl An open string analogue of Viterbo
   functoriality}, Geom.\ Topol.\ {\bf 14} (2010), 627--718.
\bibitem{AuICM}
   D. Auroux, {\sl Fukaya categories and bordered Heegaard-Floer homology},
   to appear in Proc.\ ICM 2010, arXiv:1003.2962.
\bibitem{Lekili}
   Y.~Lekili, {\sl Heegaard Floer homology of broken fibrations over the
   circle}, arXiv:0903.1773.
\bibitem{LP}
   Y.~Lekili, T.~Perutz, in preparation.
\bibitem{lipshitz}
   R.~Lipshitz, {\sl A cylindrical reformulation of Heegaard-Floer
   homology}, Geom.\ Topol.\ {\bf 10} (2006), 955--1096.
\bibitem{LMW}
   R.~Lipshitz, C.~Manolescu, J.~Wang, {\sl Combinatorial cobordism maps
   in hat Heegaard Floer theory}, Duke Math.\ J. {\bf 145} (2008), 207--247.
\bibitem{LOT} 
   R.~Lipshitz, P.~Ozsv\'ath, D.~Thurston, {\sl Bordered Heegaard Floer
   homology: invariance and pairing}, arXiv:0810.0687.
\bibitem{LOTnew}
   R.~Lipshitz, P.~Ozsv\'ath, D.~Thurston, {\sl Heegaard Floer
   homology as morphism spaces}, arXiv: 1005.1248.
\bibitem{MWW}
   S. Ma'u, K. Wehrheim, C. Woodward, {\sl $A_\infty$-functors for
   Lagrangian correspondences}, preprint.
\bibitem{Perutz}
   T.~Perutz, {\sl Lagrangian matching invariants for fibred four-manifolds:
   I},\, Geom.\ Topol.\ {\bf 11} (2007), 759--828.
\bibitem{PerHH}
   T.~Perutz, {\sl Hamiltonian handleslides for Heegaard Floer homology},
   Proc.\ 14th G\"okova Geometry-Topology Conference (2007), 
   G\"okova, 2008, 15--35, arXiv:0801.0564.
\bibitem{Sarkar}
   S.~Sarkar, {\sl Maslov index of holomorphic triangles},
   arXiv:math.GT/0609673.
\bibitem{SW}
   S.~Sarkar, J.~Wang, {\sl An algorithm for computing some Heegaard Floer
   homologies}, to appear in Ann.\ Math., arXiv:math.GT/0607777.
\bibitem{SeVCM}
   P. Seidel, {\sl Vanishing cycles and mutation},
   Proc. 3rd European Congress of Mathe\-ma\-tics (Barcelona, 2000), Vol. II,
   Progr.\ Math.\ {\bf 202}, Birkh\"auser, Basel, 2001, pp.\ 65--85,
   arXiv: math.SG/0007115.
\bibitem{SeBook}
   P. Seidel,
   {\it Fukaya categories and Picard-Lefschetz theory},
   Zurich Lect.\ in Adv.\ Math., European Math.\ Soc., Z\"urich, 2008.
\bibitem{WW}
   K. Wehrheim, C. Woodward,
   {\it Functoriality for Lagrangian correspondences in Floer theory},
   to appear in Quantum Topology, arXiv:0708.2851.
\bibitem{WW2}
   K. Wehrheim, C. Woodward,
   {\it Quilted Floer cohomology}, Geom.\ Topol.\ {\bf 14} (2010), 833--902.
\bibitem{WWsequence}
   K. Wehrheim, C. Woodward,
   {\it Exact triangle for fibered Dehn twists}, preprint.
\end{thebibliography}
\end{document}